# n- LINEAR ALGEBRA OF TYPE I
# AND ITS APPLICATIONS

**W. B. Vasantha Kandasamy**
e-mail: vasanthakandasamy@gmail.com
web: http://mat.iitm.ac.in/~wbv
www.vasantha.net

**Florentin Smarandache**
e-mail: smarand@unm.edu

# n- LINEAR ALGEBRA OF TYPE I
# AND ITS APPLICATIONS

**W. B. Vasantha Kandasamy**
**Florentin Smarandache**

**2008**

# CONTENTS

# PREFACE

With the advent of computers one needs algebraic structures that can simultaneously work with bulk data. One such algebraic structure namely n-linear algebras of type I are introduced in this book and its applications to n-Markov chains and n-Leontief models are given. These structures can be thought of as the generalization of bilinear algebras and bivector spaces. Several interesting n-linear algebra properties are proved.

This book has four chapters. The first chapter just introduces n-group which is essential for the definition of n-vector spaces and n-linear algebras of type I. Chapter two gives the notion of n-vector spaces and several related results which are analogues of the classical linear algebra theorems. In case of n-vector spaces we can define several types of linear transformations.

The notion of n-best approximations can be used for error correction in coding theory. The notion of n-eigen values can be used in deterministic modal superposition principle for undamped structures, which can find its applications in finite element analysis of mechanical structures with uncertain parameters. Further it is suggested that the concept of n-matrices can be used in real world problems which adopts fuzzy models like Fuzzy Cognitive Maps, Fuzzy Relational Equations and Bidirectional Associative Memories. The applications of these algebraic structures are given in Chapter 3. Chapter four gives some problem to make the subject easily understandable.

The authors deeply acknowledge the unflinching support of Dr.K.Kandasamy, Meena and Kama.

W.B.VASANTHA KANDASAMY
FLORENTIN SMARANDACHE

Chapter One

# BASIC CONCEPTS

In this chapter we introduce the notion of n-field, n-groups ($n \geq 2$) and illustrate them by examples. Throughout this book F will denote a field, Q the field of rationals, R the field of reals, C the field of complex numbers and $Z_p$, p a prime, the finite field of characteristic p. The fields Q, R and C are fields of zero characteristic.

Now we proceed on to define the concept of n-groups.

**DEFINITION 1.1:** *Let $G = G_1 \cup G_2 \cup \ldots \cup G_n$ ($n \geq 2$) where each $(G_i, *_i, e_i)$ is a group with $*_i$ the binary operation and $e_i$ the identity element, such that $G_i \neq G_j$, if $i \neq j$, $1 \leq j$, $i \leq n$. Further $G_i \not\subset G_j$ or $G_j \not\subset G_i$ if $i \neq j$. Any element $x \in G$ would be represented as $x = x_1 \cup x_2 \cup \ldots \cup x_n$; where $x_i \in G_i$, $i = 1, 2, \ldots, n$. Now the operations on G is described so that G becomes a group. For $x, y \in G$, where $x = x_1 \cup x_2 \cup \ldots \cup x_n$ and $y = y_1 \cup y_2 \cup \ldots \cup y_n$; with $x_i, y_i \in G_i$, $i = 1, 2, \ldots, n$.*

$$x * y = (x_1 \cup x_2 \cup \ldots \cup x_n) * (y_1 \cup y_2 \cup \ldots \cup y_n)$$
$$= (x_1 *_1 y_1 \cup x_2 *_2 y_2 \cup \ldots \cup x_n *_n y_n).$$

*Since each* $x_i *_i y_i \in G_i$ *we see* $x * y = (p_1 \cup p_2 \cup \ldots \cup p_n)$ *where* $x_i *_i y_i = p_i$ *for* $i = 1, 2, \ldots, n$. *Thus G is closed under the binary operation* $*$.

*Now let* $e = (e_1 \cup e_2 \cup \ldots \cup e_n)$ *where* $e_i \in G_i$ *the identity of* $G_i$ *with respect to the binary operation,* $*_i$, $i = 1, 2, \ldots, n$ *we see* $e * x = x * e = x$ *for all* $x \in G$. *e will be known as the identity element of G under the operation* $*$.

*Further for every* $x = x_1 \cup x_2 \cup \ldots \cup x_n \in G$; *we have* $x_1^{-1} \cup x_2^{-1} \cup \ldots \cup x_n^{-1}$ *in G such that,*

$$x * x^{-1} = (x_1 \cup x_2 \cup \ldots \cup x_n) * (x_1^{-1} \cup x_2^{-1} \cup \ldots \cup x_n^{-1})$$

$$= \quad x_1 *_1 x_1^{-1} \cup x_2 *_2 x_2^{-1} \cup \ldots \cup x_n *_n x_n^{-1}$$

$$= \quad x^{-1} * x$$

$(e_1 \cup e_2 \cup \ldots \cup e_n) = e.$

$$x^{-1} = x_1^{-1} \cup x_2^{-1} \cup \ldots \cup x_n^{-1}$$

*is known as the inverse of* $x = x_1 \cup x_2 \cup \ldots \cup x_n$. *We define* $(G, *, e)$ *to be the n-group* $(n \geq 2)$. *When* $n = 1$ *we see it is the group.* $n = 2$ *gives us the bigroup described in [37-38] when* $n > 2$ *we have the n-group.*

Now we illustrate this by examples before we proceed on to recall more properties about them.

***Example 1.1:*** Let $G = G_1 \cup G_2 \cup G_3 \cup G_4 \cup G_5$ where $G_1 = S_3$ the symmetric group of degree 3 with

$$e_1 = \begin{pmatrix} 1 & 2 & 3 \\ 1 & 2 & 3 \end{pmatrix},$$

$G_2 = \langle g \mid g^6 = e_2 \rangle$, the cyclic group of order 6, $G_3 = Z_5$, the group under addition modulo 5 with $e_3 = 0$, $G_4 = D_8 = \{a, b \mid a^2 = b^8 = 1; bab = a\}$, the dihedral group of order 8, $e_4 = 1$ is the identity element of $G_4$ and $G_5 = A_4$ the alternating subgroup of $S_4$ with

$$e_4 = \begin{pmatrix} 1 & 2 & 3 & 4 \\ 1 & 2 & 3 & 4 \end{pmatrix}.$$

Clearly $G = S_3 \cup G_2 \cup Z_5 \cup D_8 \cup A_4$ is a n-group with n = 5.

Any $x \in G$ would be of the form

$$x = \begin{pmatrix} 1 & 2 & 3 \\ 2 & 1 & 3 \end{pmatrix} \cup g^2 \cup 4 \cup b^3 \cup \begin{pmatrix} 1 & 2 & 3 & 4 \\ 1 & 3 & 4 & 2 \end{pmatrix}.$$

$$x^{-1} = \begin{pmatrix} 1 & 2 & 3 \\ 2 & 1 & 3 \end{pmatrix} \cup g^4 \cup 1 \cup b^5 \cup \begin{pmatrix} 1 & 2 & 3 & 4 \\ 1 & 4 & 2 & 3 \end{pmatrix}.$$

The identity element of G is

$$\begin{pmatrix} 1 & 2 & 3 \\ 1 & 2 & 3 \end{pmatrix} \cup e_2 \cup 0 \cup 1 \cup \begin{pmatrix} 1 & 2 & 3 & 4 \\ 1 & 2 & 3 & 4 \end{pmatrix}$$
$$= e_1 \cup e_2 \cup e_3 \cup e_4 \cup e_5.$$

Thus G is a 5-group. Clearly the order of G is $o(G_1) \times o(G_2) \times o(G_3) \times o(G_4) \times o(G_5) = 6 \times 6 \times 5 \times 16 \times 12 = 34, 560$.

We see $o(G) < \infty$. Thus if in the n-group $G_1 \cup G_2 \cup \ldots \cup G_n$, every group $G_i$ is of finite order then G is of finite order; $1 \leq i \leq n$.

***Example 1.2:*** Let $G = G_1 \cup G_2 \cup G_3$ where $G_1 = Z_{10}$, the group under addition modulo 10, $G_2 = \langle g \mid g^5 = 1\rangle$, the cyclic group of order 5 and $G_3 = Z$ the set of integers under +.

Clearly G is a 3-group. We see G is an infinite group for order of $G_3$ is infinite.

Further it is interesting to observe that every group in the 3-group G is abelian. Thus if $G = G_1 \cup G_2 \cup \ldots \cup G_n$, is a n-group $(n \geq 2)$, we see G is an abelian n-group if each $G_i$ is an abelian group; i = 1, 2, …, n. Even if one of the $G_i$ in G is a non abelian group then we call G to be only a non abelian n-group.

Having seen an example of an abelian and non abelian group we now proceed on to define the notion of n-subgroup. We need all these concepts mainly to define the new notion of linear n-algebra or n-linear algebra and n-vector spaces of type I.

**DEFINITION 1.2:** *Let $G = G_1 \cup G_2 \cup \ldots \cup G_n$, be a n-group, a proper subset $H \subset G$ of the form $H = H_1 \cup H_2 \cup \ldots \cup H_n$ with $H_i \neq G_i$ or $\{e_i\}$ but $\phi \neq H_i \subset G_i$; $i = 1, 2, \ldots, n$, $H_i$ proper subgroup of $G_i$ is defined to be the proper n-subgroup of the n-group G. If some of the $H_i = G_i$ or $H_i = \{e_i\}$ or $H_i = \phi$ for some i, then H will not be called as proper n-subgroup but only as m-subgroup of the n-group, $m < n$ and m of the subgroups $H_j$ in $G_j$ are only proper and the rest are either $\{e_j\}$ or $\phi$, $1 \leq j \leq n$.*

We illustrate both these situations by the following example.

***Examples 1.3:*** Let $G = G_1 \cup G_2 \cup G_3 \cup G_4$ be a 4-group where $G_1 = S_4$, $G_2 = Z_{10}$ group under addition modulo 10, $G_3 = D_{12}$ the dihedral group with order 12 given by the set $\{a, b \mid a^2 = b^6 = 1, bab = a\}$ and $G_4 = Z$ the set of positive and negative integers with zero under +.

Consider $H = H_1 \cup H_2 \cup H_3 \cup H_4$ where $H_1 = A_4$ the alternating subgroup of $S_4$, $H_2 = \{0, 2, 4, 6, 8\}$ a subgroup of order 5 under addition modulo 10. $H_3 = \{1, b, b^2, b^3, b^4, b^5\}$; the subgroup of $D_{12}$ and $H_4 = \{2n \mid n \in Z\}$ a subgroup of Z. Clearly H is a proper 4-subroup of the 4-group G.

Let $K = K_1 \cup K_2 \cup K_3 \cup K_4 \subseteq G$ where $K_1 = A_4$, $K_2 = \{0, 5\}$, $K_3 = D_{12}$ and $K_4 = Z$. Clearly K is not a proper 4-subgroup of the 4-group G but only a improper 4-subgroup of G.

Let $T = T_1 \cup T_2 \cup T_3 \cup T_4 \subseteq G$ where $T_1 = A_4$, $T_2 = \{0\}$, $T_3 = \phi$ and $T_4 = \{2n \mid n \in Z\}$; clearly T is only a 2-subgroup or bisubgroup of the 4-group G.

We mainly need in this book n-groups which are only abelian.

Now in this section we define the notion of n-fields.

**DEFINITION 1.3:** *Let $F = F_1 \cup F_2 \cup \ldots \cup F_n$ ( $n \geq 2$) be such that each $F_i$ is a field and $F_i \neq F_j$, if $i \neq j$ and $F_i \not\subset F_j$ or $F_j \not\subset F_i$, $1 \leq i, j \leq n$. Then we define $(F, +, \times)$ to be a n-field if $(F, +)$ is a* *n-group and* $F_1 \setminus \{0\} \cup F_2 \setminus \{0\} \cup \ldots \cup F_n \setminus \{0\}$ *is a n-group under* $\times$.

*Further*

$$[(a_1 \cup a_2 \cup \ldots \cup a_n) + (b_1 \cup b_2 \cup \ldots \cup b_n)] \times [(c_1 \cup c_2 \cup \ldots \cup c_n)]$$
$$= (a_1 + b_1) \times c_1 \cup (a_2 + b_2) \times c_2 \cup \ldots \cup (a_n + b_n) \times c_n$$

*and*

$$[(c_1 \cup c_2 \cup \ldots \cup c_n)] \times \{[(a_1 \cup a_2 \cup \ldots \cup a_n)] + [(b_1 \cup b_2 \cup \ldots \cup b_n)]\}$$
$$= c_1 \times (a_1 \cup b_1) \cup c_2 \times (a_2 \cup b_2) \cup \ldots \cup c_n \times (a_n \cup b_n)$$

*for all* $a_i, b_i, c_i \in F$, $i = 1, 2, \ldots, n$. *Thus* $(F, +, \times)$ *is a n-field.*

We illustrate this by the following example.

***Example 1.4:*** Let $F = F_1 \cup F_2 \cup F_3 \cup F_4$ where $F_1 = Q$, $F_2 = Z_2$, $F_3 = Z_{17}$ and $F_4 = Z_{11}$; F is a 4-field.

***Example 1.5:*** Let $F = F_1 \cup F_2 \cup F_3 \cup F_4 \cup F_5 \cup F_6$ where $F_1 = Z_2$, $F_2 = Z_3$, $F_3 = Z_{13}$, $F_4 = Z_7$, $F_5 = Z_{19}$ and $F_6 = Z_{31}$, F is a 6-field. Let $F = F_1 \cup F_2 \cup \ldots \cup F_n$, be a n-field where each $F_i$ is a field of characteristic zero, $1 \leq i \leq n$, then F is called as a n-field of characteristic zero.

Let $F = F_1 \cup F_2 \cup \ldots \cup F_m$ $(m \geq 2)$ be a m-field if each field $F_i$ is of finite characteristic then we call F to be a m-field of finite characteristic. Suppose $F = F_1 \cup F_2 \cup \ldots \cup F_n$, $n \geq 2$ where some $F_i$'s are finite characteristic and some $F_j$'s are zero characteristic then alone we say F is a n-field of mixed characteristic.

***Example 1.6:*** Let $F = F_1 \cup F_2 \cup \ldots \cup F_5$ where $F_1 = Q$, $F_2 = Z_7$, $F_3 = Z_{23}$ and $F_4 = Z_{17}$ and $F_5 = Z_2$ be a 5-field. F is a 5-field of mixed characteristic.

***Example 1.7:*** Let $F = F_1 \cup F_2 \cup \ldots \cup F_6 = Z_4 \cup R \cup Z_7 \cup Q \cup R \cup Z_{11}$. Clearly F is not a 6 field as $F_2 = F_5$. We need each field $F_i$ to be distinct, $1 \leq i \leq n$.

Note: Clearly $F_1 \cup F_2 \cup F_3 = Q \cup R \cup Z_2$ is not a 3-field as $Q \subseteq R$. Because we need in this case also as in case of bistructures non containment of one set in another set.

## Chapter Two

# n-VECTOR SPACES OF TYPE I AND THEIR PROPERTIES

In this chapter we introduce the notion of n-vector spaces and describe some of their important properties.

Here we define the concept of n-vector spaces over a field which will be known as the type I n-vector spaces or n-vector spaces of type I. Several interesting properties about them are derived in this chapter.

**DEFINITION 2.1:** *A n-vector space or a n-linear space of type I ($n \geq 2$) consists of the following:*

1. *a field F of scalars*
2. *a set $V = V_1 \cup V_2 \cup \ldots \cup V_n$ of objects called n-vectors*
3. *a rule (or operation) called vector addition; which associates with each pair of n-vectors $\alpha = \alpha_1 \cup \alpha_2 \cup \ldots \cup \alpha_n$, $\beta = \beta_1 \cup \beta_2 \cup \ldots \cup \beta_n \in V = V_1 \cup V_2 \cup \ldots \cup V_n$; $\alpha + \beta = (\alpha_1 \cup \alpha_2 \cup \ldots \cup \alpha_n) + (\beta_1 \cup \beta_2 \cup \ldots \cup \beta_n) = (\alpha_1 + \beta_1 \cup \alpha_2 + \beta_2 \cup \ldots \cup \alpha_n + \beta_n) \in V$ called the sum of $\alpha$ and $\beta$ in such a way*
   a. *$\alpha + \beta = \beta + \alpha$; i.e., addition is commutative ($\alpha, \beta \in V$).*

b. $\alpha + (\beta + \gamma) = (\alpha + \beta) + \gamma$, i.e., addition is associative ($\alpha, \beta, \gamma \in V$).

c. There is a unique n-vector $0_n = 0 \cup 0 \cup \ldots \cup 0 \in V$ such that $\alpha + 0_n = \alpha$ for all $\alpha \in V$, called the zero n-vector of V.

d. For each n-vector $\alpha = \alpha_1 \cup \alpha_2 \cup \ldots \cup \alpha_n \in V$, there exists a unique vector $-\alpha = -\alpha_1 \cup -\alpha_2 \cup \ldots \cup -\alpha_n \in V$ such that $\alpha + (-\alpha) = 0_n$.

e. A rule (or operation) called scalar multiplication which associates with each scalar c in F and a n-vector α in $V = V_1 \cup V_2 \cup \ldots \cup V_n$ a n-vector cα in V called the product of c and α in such a way that

1. $$\begin{aligned} 1.\alpha &= 1.(\alpha_1 \cup \alpha_2 \cup \ldots \cup \alpha_n) \\ &= 1.\alpha_1 \cup 1.\alpha_2 \cup \ldots \cup 1.\alpha_n \\ &= \alpha_1 \cup \alpha_2 \cup \ldots \cup \alpha_n \\ &= \alpha \end{aligned}$$

for every n-vector α in V.

2. $(c_1. c_2).\alpha = c_1.(c_2. \alpha)$ for all $c_1, c_2 \in F$ and $\alpha \in V$ i.e. if $\alpha_1 \cup \alpha_2 \cup \ldots \cup \alpha_n$ is the n-vector in V we have

$$\begin{aligned} (c_1. c_2).\alpha &= (c_1. c_2)(\alpha_1 \cup \alpha_2 \cup \ldots \cup \alpha_n) \\ &= c_1 [c_2((\alpha_1 \cup \alpha_2 \cup \ldots \cup \alpha_n)] \\ &= c_1 [c_2\alpha_1 \cup c_2\alpha_2 \cup \ldots \cup c_2\alpha_n] \\ &= c_1 [c_2\alpha]. \end{aligned}$$

3. $c(\alpha + \beta) = c.\alpha + c.\beta$ for all $\alpha, \beta \in V$ and for all $c \in F$ i.e., if $\alpha_1 \cup \alpha_2 \cup \ldots \cup \alpha_n$ and $\beta_1 \cup \beta_2 \cup \ldots \cup \beta_n$ are n-vectors of V then for any $c \in F$ we have

$$\begin{aligned} c(\alpha + \beta) &= c[(\alpha_1 \cup \alpha_2 \cup \ldots \cup \alpha_n) + (\beta_1 \cup \beta_2 \cup \ldots \cup \beta_n)] \\ &= c[\alpha_1 + \beta_1 \cup \alpha_2 + \beta_2 \cup \ldots \cup \alpha_n + \beta_n] \\ &= (c(\alpha_1 + \beta_1) \cup c(\alpha_2 + \beta_2) \cup \ldots \cup c(\alpha_n + \beta_n)] \\ &= (c\alpha_1 \cup c\alpha_2 \cup \ldots \cup c\alpha_n) + (c\beta_1 \cup c\beta_2 \cup \ldots \cup c\beta_n) \\ &= c\alpha + c\beta. \end{aligned}$$

4. $(c_1 + c_2).\alpha = c_1\alpha + c_2\alpha$ for all $c_1, c_2 \in F$ and $\alpha \in V$.

*Just like a vector space which is a composite algebraic structure containing the field a set of vectors which form a group, the n-vector space of type I is a composite of set of n-vectors or n-group and a field F of scalars. V is a linear n-algebra or n-linear algebra if V has a multiplicative closed binary operation "." which is associative i.e.; if* $\alpha, \beta \in V$, $\alpha.\beta \in V$, *thus if* $\alpha = (\alpha_1 \cup \alpha_2 \cup \ldots \cup \alpha_n)$ *and* $\beta = (\beta_1 \cup \beta_2 \cup \ldots \cup \beta_n) \in V$ *then if*

$$\begin{aligned} \alpha.\beta &= (\alpha_1 \cup \alpha_2 \cup \ldots \cup \alpha_n) . (\beta_1 \cup \beta_2 \cup \ldots \cup \beta_n) \\ &= (\alpha_1.\beta_1 \cup \alpha_2.\beta_2 \cup \ldots \cup \alpha_n.\beta_n) \in V \end{aligned}$$

*then the linear n-vector space of type I becomes a linear n-algebra of type-I.*

Now we make an important mention that all linear n-algebras of type-I are linear n-vector spaces of type-I; however a n-vector space of type-I over F in general need not be a n- linear algebra of type I over F.

We now illustrate this by the following example.

***Example 2.1:*** Let $V = V_1 \cup V_2 \cup V_3 \cup V_4$ where $V_1 = Q[x]$ the vector space of polynomials over Q. $V_2 = Q \times Q$, the vector space of dimension two over Q,

$$V_3 = \left\{ \begin{pmatrix} a & b \\ c & d \end{pmatrix} \middle| a, b, c, d \in Q \right\}$$

the vector space of all $2 \times 2$ matrices with entries from Q and

$$V_4 = \left\{ \begin{pmatrix} a & b & c \\ d & e & f \end{pmatrix} \middle| a, b, c, d, e, f \in R \right\}$$

be the vector space of all $2 \times 3$ matrices with entries from R over Q. Thus V is a linear 4-vector space over Q of type-I. Clearly V is not a linear 4-algebra of type-I over Q.

Now we give yet another example of a linear n-vector space of type-I.

***Example 2.2:*** Let $V = V_1 \cup V_2 \cup V_3 \cup V_4 \cup V_5$ be a 5-vector space over Q of type-I, where $V_1$ = Q[x], the set of all polynomials with coefficients from Q is a vector space over Q. $V_2 = Q \times R \times Q$ is a vector space over Q,

$$V_3 = \left\{ \begin{pmatrix} a & b & c \\ d & e & f \\ g & h & i \end{pmatrix} \middle| \, a, b, c, d, e, f, g, h, i \in Q \right\}$$

is a vector space over Q,

$$V_4 = \left\{ \begin{pmatrix} a & 0 & 0 & 0 \\ 0 & b & 0 & 0 \\ 0 & 0 & c & 0 \\ 0 & 0 & 0 & d \end{pmatrix} \middle| \, a, b, c, d, \in R \right\}$$

is a vector space over Q and $V_5$ = R is a vector space over Q. Clearly $V = V_1 \cup V_2 \cup V_3 \cup V_4 \cup V_5$ is a linear 5-vector space of type-I over Q. Also V is a linear 5-linear algebra over Q. Thus we have seen from example 2.1 that every vector n-space of type-I need not be a linear n-algebra of type-I. Also every linear n-algebra of type-I is a linear n-vector space of type-I.

Now we can also define the notion of n-vector space of type-I in a very different way.

**DEFINITION 2.2:** *Let $V = V_1 \cup V_2 \cup \ldots \cup V_n$ $(n \geq 2)$ where each $V_i$ is a vector space over the same field F and $V_i \neq V_j$, if $i \neq j$ and $V_i \not\subset V_j$ and $V_j \not\subset V_i$ if $i \neq j$, $1 \leq i, j \leq n$, then V is defined to be a n-vector space of type-I over F.*

*If each of the $V_i$'s are linear algebra over F then we call V to be a linear n-algebra of type-I over F.*

Now we proceed on to define the notion of n-subvector space of the n-vector space of type-I.

**DEFINITION 2.3:** *Let $V = V_1 \cup V_2 \cup \ldots \cup V_n$ ($n \geq 2$) be a n-vector space of type I over F. Suppose $W = W_1 \cup W_2 \cup \ldots \cup W_n$ ($n \geq 2$) is a proper subset of V such that each $W_i$ is a proper subspace of the vector space $V_i$ over F with $W_i \neq V_i$, $W_i \neq \phi$ or (0) such that $W_i \neq W_j$ or $W_i \not\subset W_j$ or $W_j \not\subset W_i$ if $i \neq j$, $1 \leq i, j \leq n$, then we define W to be a n-subspace of type-I over F.*

We now illustrate it by the following example.

***Example 2.3:*** Let $V = V_1 \cup V_2 \cup V_3$ where $V_1 = R \times R$, a vector space over R and $V_2 = R[x]$ a vector space over R and

$$V_3 = \left\{ \begin{pmatrix} a & c \\ d & b \end{pmatrix} \middle| a, b, c, d \in R \right\},$$

a vector space over R i.e., V is a 3-vector space of type-I over R. Let $W = W_1 \cup W_2 \cup W_3 \subset V = V_1 \cup V_2 \cup V_3$ where

$$W_1 = R \times \{0\} \subset V_1,$$

$$W_2 = \left\{ \sum_{i=0}^{n} r_i x^{2i} \middle| r_i \in R \right\} \subset V_2,$$

$$W_3 = \left\{ \begin{pmatrix} a & 0 \\ 0 & b \end{pmatrix} \middle| a, b \in R \right\} \subseteq V_3.$$

Clearly W is a 3 subspace of V of type-I. Suppose

$$T = R \times \{0\} \cup R \cup \left\{ \begin{pmatrix} a & 0 \\ 0 & b \end{pmatrix} \middle| a, b \in R \right\} \subseteq V_1 \cup V_2 \cup V_3,$$

then T is not a 3-subspace of type-I as $R \times \{0\}$ and R are same or $R \subseteq R \times \{0\}$.

Now we proceed on to define the notion of n-linear dependence and n-linear independence in the n-vector space V of type-I.

**DEFINITION 2.4:** *Let $V = V_1 \cup V_2 \cup \ldots \cup V_n$ be a n-vector space of type-I over F. Any proper n-subset $S \subseteq V$ would be of the form $S = S_1 \cup S_2 \cup \ldots \cup S_n \subseteq V_1 \cup V_2 \cup \ldots \cup V_n$ where $\phi \neq S_i$ contained in $V_i$, $1 \leq i \leq n$. $S_i$ a proper subset of $V_i$. If each of the subsets $S_i \subseteq V_i$ is a linearly independent set over F for i = 1, 2, …, n then we define S to be a n-linearly independent subset of V. Even if one of the subset $S_k$ of $V_k$ is not a linearly independent subset of $V_k$ for some $1 \leq k \leq n$ then we call the n-subset of V to be a n-linearly dependent subset or a linearly dependent n-subset of V.*

Now we illustrate this situation by the following examples.

***Example 2.4:*** Let $V = V_1 \cup V_2 \cup V_3 \cup V_4$ be a 4- vector space over Q, where $V_1 = Q[x]$, $V_2 = Q \times Q \times Q$; $V_3$ = { the set of all 2 × 2 matrices with entries from Q} and $V_4$ = [the set of all 4 × 2 matrices with entries from Q, are all vector spaces over Q. Let $S = S_1 \cup S_2 \cup S_3 \cup S_4$ be a 4 subset of V,

$$S_1 = \{1, x^2, x^5, x^7, 3x^8\},$$
$$S_2 = \{(7, 0, 2), (0, 5, 1)\},$$
$$S_3 = \left\{\begin{pmatrix} 5 & 1 \\ 0 & 0 \end{pmatrix}, \begin{pmatrix} 0 & 0 \\ 7 & 3 \end{pmatrix}\right\}$$

and

$$S_4 = \left\{\begin{bmatrix} 0 & 2 \\ 1 & 0 \\ 0 & 0 \\ 3 & 0 \end{bmatrix}, \begin{bmatrix} 1 & 0 \\ 0 & 2 \\ 0 & 0 \\ 0 & 1 \end{bmatrix}, \begin{bmatrix} 0 & 0 \\ 0 & 0 \\ 7 & 3 \\ 0 & 1 \end{bmatrix}\right\}.$$

Clearly we see every subset $S_i$ of $V_i$ is a linearly independent subset, for i = 1, 2, 3, 4. Thus S is a 4- linearly independent subset of V.

***Example 2.5:*** Let $V = V_1 \cup V_2 \cup V_3$ be a 3-vector space over R where $V_1 = R[x]$, $V_2$ = {set of all $3 \times 3$ matrices with entries from R} and $V_3 = R \times R \times R \times R$. Clearly $V_1$, $V_2$ and $V_3$ are all vector spaces over R. Let $S = S_1 \cup S_2 \cup S_3 \subseteq V_1 \cup V_2 \cup V_3 = V$ be a proper 3-subset of V; where

$$S_1 = \{x^3, 3x^3 + 7, x^5\},$$

$$S_2 = \left\{ \begin{pmatrix} 6 & 0 & 0 \\ 0 & 0 & 3 \\ 1 & 1 & 0 \end{pmatrix}, \begin{pmatrix} 0 & 1 & -2 \\ 1 & 0 & 1 \\ 0 & 7 & 0 \end{pmatrix} \right\}$$

and

$$S_3 = \{(3\ 1\ 0\ 0), (0\ 7\ 2\ 1), (5\ 1\ 1\ 1), (0\ 8\ 9\ 1), (2\ 1\ 3\ 0)\}.$$

We see $S_1$ is a linearly dependent subset of $V_1$ over R and $S_2$ is a linearly independent subset over R and $S_3$ is a linearly dependent subset of $V_3$ over R. Thus S is a 3-linearly dependent subset of the 3-vector space V over R.

Now we proceed onto define the notion of n-basis of the n-vector space V over a field F.

**DEFINITION 2.5:** *Let $V = V_1 \cup V_2 \cup \ldots \cup V_n$ be a n-vector space over a field F. A proper n-subset $S = S_1 \cup S_2 \cup \ldots \cup S_n$ of V is said to be n-basis of V if S is a n-linearly independent set and each $S_j \subseteq V_j$ generates $V_j$, i.e., $S_j$ is a basis of $V_j$, true for j = 1, 2, …, n. Even if one of the $S_j$ is not a basis of $V_j$ for $1 \leq j \leq n$ then S is not a n-basis of V.*

As in case of vector spaces the n-vector spaces can also have many basis but the number of base elements in each of the n subsets is the same.

Now we illustrate this situation by the following example.

***Example 2.6 :*** Let $V = V_1 \cup V_2 \cup V_3 \cup V_4$ be a 4-vector space over Q. $V_1$ = {all polynomials of degree less than or equal to 5}, $V_2 = Q \times Q \times Q$, $V_3$ = {the set of all 2×2 matrices with entries from Q} and $V_4 = Q \times Q \times Q \times Q \times Q$ are vector spaces over Q. Now let

$$B = B_1 \cup B_2 \cup B_3 \cup B_4$$
$$= \{1, x, x^2, x^3, x^4, x^5\} \cup \{(1\ 0\ 0), (0\ 1\ 0), (0\ 2\ 1)\} \cup$$

$$\left\{\begin{pmatrix} 0 & 0 \\ 1 & 0 \end{pmatrix}, \begin{pmatrix} 1 & 0 \\ 0 & 0 \end{pmatrix}, \begin{pmatrix} 0 & 0 \\ 0 & 1 \end{pmatrix} \begin{pmatrix} 0 & 1 \\ 0 & 0 \end{pmatrix}\right\} \cup \{(0\ 0\ 0\ 0\ 1), (0\ 0\ 0\ 1\ 0), (0\ 0\ 1\ 0\ 0), (0\ 1\ 0\ 0\ 0), (1\ 0\ 0\ 0\ 0)\}$$
$$\subseteq V_1 \cup V_2 \cup V_3 \cup V_4 = V.$$

B is a 4-basis of V as each $B_i$ is a basis of $V_i$ ; i = 1, 2, 3, 4.

***Example 2.7:*** Let $V = V_1 \cup V_2 \cup V_3 \cup V_4 \cup V_5$ be a 5-vector space over Q where $V_1 = R$, $V_2 = Q \times Q$, $V_3 = Q[x]$, $V_4 = R \times R \times R$ and $V_5$ = {set of all 2×2 matrices with entries from Q}. Clearly $V_1$, $V_2$, $V_3$, $V_4$ and $V_5$ are vector spaces over Q. We see some of the vector spaces $V_i$ over Q are finite dimensional i.e., has finite basis and some of the vector spaces $V_j$ have infinite number of elements in the basis set. We find means to define the new notion of finite n-dimensional space and infinite n-dimensional space. To be more specific in this example, $V_1$ is an infinite dimensional vector space over Q, $V_2$ and $V_3$ are finite dimensional vector spaces over Q. $V_4$ is an infinite dimensional vector space over Q and $V_5$ is a finite dimensional vector space over Q.

**DEFINITION 2.6:** *Let $V = V_1 \cup V_2 \cup \ldots \cup V_n$ be a n-vector space of type-I over F. If every vector space $V_i$ in V is finite dimensional over F then we say the n-vector space is finite n-dimensional over F. Even if one of the vector space $V_j$ in V is infinite dimensional then we say V is infinite dimensional over F. We denote the dimension of V by $(n_1, n_2, \ldots, n_n)$; $n_i$ dimension of $V_i$, i = 1, 2, ..., n.*

We illustrate the definition by some examples.

***Example 2.8:*** Let $V = V_1 \cup V_2 \cup V_3$ be a 3-vector space over Q, where $V_1 = Q[x]$, $V_2$ = {set of all 2×2 matrices with entries from Q} and $V_3 = Q$ the one dimensional vector space over Q. Clearly 3-dimension of the 3-vector space over Q is (∞, 4, 1). Thus V is an infinite 3-dimensional space over Q.

***Example 2.9:*** Let $V = V_1 \cup V_2 \cup V_3 \cup V_4$ be a 4-vector space of type-I over Q. Suppose $V_1$ = {set of all 2 × 2, matrices with entries from Q}; $V_2 = Q \times Q \times Q$ a vector space over Q, $V_3$ = {All polynomials of degree less than or equal to 7 with coefficients from Q} and $V_4$ = {the collection of all 5 × 5, matrices with entries from Q}, we see $V_1$, $V_2$, $V_3$ and $V_4$ are vector spaces over Q. The 4-dimension of V is (4, 3, 8, 25), so, V is finite 4-dimension 4 vector space over Q of type-I.

Having seen sub n-spaces, n-basis and n-dimension of n-vector spaces of type-I now we proceed on to define the notion of n-transformation of n-vector space of type-I.

**DEFINITION 2.7:** *Let $V = V_1 \cup V_2 \cup \ldots \cup V_n$ be a n-vector space over a field F of type-I and $W = W_1 \cup W_2 \cup \ldots \cup W_m$ be another m-vector space over the same field F of type I, $(n \leq m)$ $(m \geq 2)$ and $(n \geq 2)$. We call T a n-map if $T = T_1 \cup T_2 \cup \ldots \cup T_n$: $V \to W$ is defined as $T_i : V_i \to W_j$, $1 \leq i \leq n$, $1 \leq j \leq m$ for every i. If each $T_i$ is a linear transformation from $V_i$ to $W_j$, i = 1, 2, …, n, $1 \leq j \leq n$ then we call the n-map to be a n-linear transformation from V to W or linear n-transformation from V to W. No two $V_i$'s are mapped on to the same $W_j$, $1 \leq i \leq n$, $1 \leq j \leq m$. Even if one of the $T_i$ is not a linear transformation from $V_i$ to $W_j$ then T is not a n-linear transformation.*

We will illustrate this by the simple example.

***Example 2.10:*** Let $V = V_1 \cup V_2 \cup V_3$ be a 3-vector space over Q and $W = W_1 \cup W_2 \cup W_3 \cup W_4$ be a 4-vector space over Q. V is of finite (3, 2, 4) dimension and W is of finite (4, 3, 2, 4) dimension. T: V → W be a 3-linear transformation defined by $T = T_1 \cup T_2 \cup T_3$: $V_1 \cup V_2 \cup V_3 \to W_1 \cup W_2 \cup W_3 \cup W_4$ as

$T_1: V_1 \to W_1$ given by

$$T_1(a_1^1, a_2^1, a_3^1) = (a_1^1 + a_2^1, a_2^1, a_3^1 + a_2^1, a_1^1 + a_2^1 + a_3^1)$$

$T_2: V_2 \to W_3$ defined by

$$T_2(a_1^1, a_2^2) = (a_1^2 + a_2^2, a_1^2)$$

and $T_3 : V_3 \to W_4$ defined by

$$T_3(a_1^3, a_2^3, a_3^3, a_4^3) = (a_2^3, a_4^3, a_4^3, a_1^3 + a_2^3),$$

clearly T is a 3-linear transformation or linear 3-transformation or linear 3 transformation of V to W; i.e. from 3-vector space V to 4-vector space W.

It may so happen that we may have a n-vector space over a field F and it would become essential for us to make a linear n-transformation to a m-vector space over F where n>m. In such situation we define a linear n-transformation which we call as shrinking linear n-transformation which is as follows.

**DEFINITION 2.8:** *Let V be a n-vector space over F and W a m-vector space over F n > m. The shrinking n-map T from* $V = V_1 \cup V_2 \cup \ldots \cup V_n$ *to* $W = W_1 \cup W_2 \cup \ldots \cup W_m$ *is defined as a map from V to W as follows* $T = T_1 \cup T_2 \cup \ldots \cup T_n$ *with* $T_i : V_i \to W_j$ *;* $1 \leq i \leq n$ *and* $1 \leq j \leq m$ *with the condition* $T_j : V_j \to W_k$ *where j may be equal to k. i.e. the range space as in case of linear n-map may not be distinct.*

*Now if* $T_i : V_i \to W_j$ *in addition a linear transformation then we call,* $T = T_1 \cup T_2 \cup \ldots \cup T_n$ *the shrinking n-map to be a shrinking linear n-transformation or a shrinking n-linear transformation.*

We illustrate this situation by the following example.

***Example 2.11:*** Let $V = V_1 \cup V_2 \cup V_3 \cup V_4 \cup V_5$ be a 5-vector space defined over Q of 5-dimenion (3, 2, 5, 7, 6) and $W = W_1 \cup W_2 \cup W_3$ be a 3-vector space defined over Q of 3-dimension (5, 3, 6). $T = T_1 \cup T_2 \cup \ldots \cup T_5 : V \to W$ can only be a shrinking 5-linear transformation defined by

$$T_1 : V_1 \to W_3,$$
$$T_2 : V_2 \to W_1,$$
$$T_3 : V_3 \to W_2,$$
$$T_4 : V_4 \to W_3$$

and

$$T_5 : V_5 \to W_1$$

where

$$T_1(x_1^1, x_2^1, x_3^1) = (x_1^1 + x_2^1, x_3^1, x_2^1, x_2^1 + x_3^1, x_1^1 + x_3^1, x_1^1)$$

where $x_1^1, x_2^1, x_3^1 \in V_1$,

$$T_2(x_1^2, x_2^2) = (x_1^2 + x_2^2, x_2^2, x_1^2 + x_2^2, x_1^2, x_2^2),$$

where $x_1^2, x_2^2 \in V_2$,

$$T_3(x_1^3, x_2^3, x_3^3, x_4^3, x_5^3) = (x_1^3 + x_2^3, x_2^3 + x_3^3, x_4^3 + x_5^3)$$

where $x_1^3, x_2^3, x_3^3, x_4^3$ and $x_5^3 \in V_3$,

$$T_4(x_1^4, x_2^4, x_3^4, x_4^4, x_5^4, x_6^4, x_7^4)$$
$$= (x_1^4 + x_2^4, x_2^4 + x_3^4, x_3^4 + x_4^4, x_4^4 + x_5^4, x_5^4 + x_6^4, x_6^4 + x_7^4)$$

for all $x_1^4, x_2^4, \ldots, x_7^4 \in V_4$ and

$$T_5(x_1^5, x_2^5, x_3^5, x_4^5, x_5^5, x_6^5) = (x_1^5, x_2^5, x_3^5 + x_4^5, x_5^5, x_6^5)$$

for $x_1^5, x_2^5, \ldots, x_6^5 \in V_5$.

Clearly T is a shrinking linear 5 transformation.

*Note:* It may be sometimes essential for one to define a linear n-transformation from a n-vector space V into a m-vector space W, m > n where all the n spaces of the m-vector space may not be used only a set of r vector spaces from W may be needed r < n < m, in such cases we call the linear n-transformation as a special shrinking linear n-transformation of V into W.

We illustrate this situation by the following example.

***Example 2.12:*** Let $V = V_1 \cup V_2 \cup V_3$ be a 3-vector space over Q and $W = W_1 \cup W_2 \cup W_3 \cup W_4 \cup W_5$ be a 5-vector space over Q. Suppose V is a finite 3-dimension (3, 5, 4) space and W be a finite 5-dimension (3, 5, 4, 8, 2) space. Let $T = T_1 \cup T_2 \cup T_3 : V \to W$ be defined by $T_1: V_1 \to W_1$, $T_2: V_2 \to W_3$, $T_3: V_3 \to W_1$ as follows;

$$T_1(x_1^1, x_2^1, x_3^1) = (x_1^1 + x_2^1, x_2^1 + x_3^1, x_2^1)$$

for all $x_1^1, x_2^1, x_3^1 \in V_1$;

$$T_2(x_1^2, x_2^2, x_3^2, x_4^2, x_5^2) = (x_2^2, x_1^2, x_3^2 + x_5^2, x_4^2)$$

for all $x_1^2, x_2^2, x_3^2, x_4^2, x_5^2 \in V_2$,

$$T_3(x_1^3, x_2^3, x_3^3, x_4^3) = (x_1^3 + x_2^3, x_4^3 + x_1^3, x_2^3 + x_3^3)$$

for all $x_1^3, x_2^3, x_3^3$ and $x_4^3$ in $V_3$.

Thus $T:V \to W$ is only a special shrinking linear 3-transformation.

**DEFINITION 2.9:** *Let V be a n-vector space over the field F and W be a n-vector space over the same field F. $T = T_1 \cup T_2 \cup \ldots \cup T_n$ is a linear one to one n transformation if each $T_i$ is a transformation from $V_i$ to $W_j$ and for no $V_k$ we have $T_k : V_k \to W_j$ i.e. no two distinct domain space can have the same range space. Then we call T to be a one to one vector space preserving linear n-transformation.*

We just show this by a simple example.

***Example 2.13:*** Let $V = V_1 \cup V_2 \cup V_3 \cup V_4$ be a 4- vector space over Q and $W = W_1 \cup W_2 \cup W_3 \cup W_4$ be another 4-vector space over Q. Let V be of (3, 4, 5, 2) finite 4 dimensional space and W a (2, 5, 6, 3) finite 4-dimensional space. Let $T = T_1 \cup T_2 \cup T_3 \cup T_4: V = V_1 \cup V_2 \cup V_3 \cup V_4 \to W_1 \cup W_2 \cup W_3 \cup W_4$ given by $T_1: V_1 \to W_2$, $T_2: V_2 \to W_3$, $T_3: V_3 \to W_4$ and $T_4: V_4 \to W_1$ where $T_1$, $T_2$, $T_3$ and $T_4$ are linear transformation. Clearly T is a linear one to one 4-transformation.

*Note:* In the definition 2.9 it is interesting and important to note that all $T_i$'s need not be 1-1 linear transformation with dim $V_i$ = dim $W_j$ if $T_i$: $V_i \to W_j$ i.e., $T_i$'s are not vector space isomorphism for i = 1, 2, …, n. Now we give a new name for a n-linear transformation T: V → W where $T = T_1 \cup T_2 \cup \ldots \cup T_n$ with each $T_i$ a vector space isomorphism or $T_i$ is 1-1 and onto linear transformation from $V_i$ to $W_j$, $1 \leq i \leq n$, $1 \leq j \leq n$.

**DEFINITION 2.10:** *Let V and W be n vector spaces defined over a field F. We say V and W are of same n-dimension if and only if n-dimension of V is ($n_1$, …, $n_n$) then the n-dimension of W is just a permutation of ($n_1$, $n_2$, … , $n_n$).*

***Example 2.14:*** Let $V = V_1 \cup V_2 \cup V_3 \cup V_4 \cup V_5$ be a 5-dimension vector space over R of 5-dimension (7, 2, 3, 4, 5). Suppose $W = W_1 \cup W_2 \cup W_3 \cup W_4 \cup W_5$ is a 5-dimension vector space over R of 5-dimension (2, 5, 4, 7, 3) then we say V and W are of same 5-dimension. If $X = X_1 \cup X_2 \cup X_3 \cup X_4 \cup X_5$ is a 5-vector space of 5-dimension (2, 7, 9, 3, 4) then clearly X and V are not 5-vector spaces of same dimension. So for any n-dimensional n-vector space V we have only $|\underline{n}$ number of n-vector spaces of same dimension including V.

We just show this by an example.

***Example 2.15:*** Let $V = V_1 \cup V_2 \cup V_3$ be a 3-vector space of 3-dimension (7, 5, 3). Then W, X, Y, Z and S of 3-dimension (5, 7, 3), (5, 3, 7), (7, 3, 5), (3, 5, 7) and (3, 7, 5) are of same dimension.

In view of this we have the following interesting theorem.

**THEOREM 2.1:** *Let V be a finite n-dimension n-vector space over the field F of n-dimension ($n_1$, $n_2$, …, $n_n$), then their exist $|\underline{n}$ finite n-dimension n-vector spaces of same dimension as that of V including V over F.*

*Proof:* Given V is a finite n-vector space of n-dimension $(n_1, n_2, \ldots, n_n)$ i.e. each $1 \leq n_i \leq \infty$ and $i \neq j$ implies $n_i \neq n_j$ we know two n-vector spaces V and W are of same dimension if and only if the n-dimension of one (say V) can be got from permuting the n-dimension of W, or vice versa. Further from group theory we know for a set (1, 2, …, n) we have $\lfloor\underline{n}$ permutations of the set (1, 2, …, n). Thus we have $\lfloor\underline{n}$ n-vector spaces of dimension $(n_1, n_2, \ldots, n_n)$.

*Note:* If we have a n-vector space of n-dimension $(m_1, m_2, \ldots, m_n)$ with some $m_i \neq n_j$, $1 \leq i \leq n$ then we get another set of $\lfloor\underline{n}$ n-vector spaces of n-dimension $(m_1, m_2, \ldots, m_n)$ and all its permutations. Clearly this set of m-vector spaces with n-dimension $(n_1, n_2, \ldots, n_n)$ are distinct from the n-vector spaces of n-dimension $(m_1, m_2, \ldots, m_n)$. From this one can conclude we have infinite number of n-vector spaces of varying dimensions. Only same n-dimension vector spaces can be n-isomorphic.

**DEFINITION 2.11:** *Let V and W be n-vector spaces of same dimension. Let n-dimension of V be $(n_1, n_2, \ldots, n_n)$ and that of W be $(n_4, n_2, n_n, \ldots, n_5)$ i.e. let $V = V_1 \cup V_2 \cup \ldots \cup V_n$ and $W = W_1 \cup W_2 \cup \ldots \cup W_n$. A linear n-transformation $T = T_1 \cup T_2 \cup \ldots \cup T_n : V \to W$ is defined to be a n-vector space linear n-isomorphism if and only if $T_i : V_i \to W_j$ is such that $\dim V_i = \dim W_j$; $1 \leq i, j \leq n$.*

We illustrate this situation by an example.

***Example 2.16:*** Let $V = V_1 \cup V_2 \cup V_3 \cup V_4 \cup V_5$ and $W = W_1 \cup W_2 \cup W_3 \cup W_4 \cup W_5$ be two 5-vector spaces of same dimension. Let the 5-dimension of V and W be (3, 2, 5, 4, 6) and (4, 2, 5, 3, 6) respectively. Suppose $T = T_1 \cup T_2 \cup T_3 \cup T_4 \cup T_5$: $V \to W$ given by $T_1(V_1) = W_4$, $T_2(V_2) = W_2$, $T_3(V_3) = W_3$, $T_4(V_4) = W_1$ and $T_5(V_5) = W_5$; then T is a one to one n-isomorphic, n-linear transformation of V to W (n = 5). Suppose $P : V \to W$ where $P = P_1 \cup P_2 \cup P_3 \cup P_4 \cup P_5$ given by $P_1$:$V_1 \to W_2$, $P_2$: $V_2 \to W_3$, $P_3$: $V_3 \to W_4$, $P_4$: $V_4 \to W_5$, and $P_5$: $V_5 \to W_1$ the linear transformation so that P is a 5-linear transformation from V to W. Clearly P is not the one to one isomorphic 5-linear transformation of V. P is only a one to one 5-linear transformation of V.

Now having seen different types of linear n-transformation of a n-vector space V to W, W a linear n-space we proceed on to define the notion of n-kernel of T.

**DEFINITION 2.12:** *Let $V = V_1 \cup V_2 \cup \ldots \cup V_n$ be a n-vector space over the field F and $W = W_1 \cup W_2 \cup \ldots \cup W_m$ be a m-vector space over the field F. Let $T = T_1 \cup T_2 \cup \ldots \cup T_n$ be n-linear transformation of T from V to W defined by $T_i: V_i \to W_j$; $1 \leq i \leq n$ and $1 \leq j \leq m$ such that no two domain spaces are mapped on to the same range space. The n-kernel of T denoted by ker T = ker $T_1 \cup$ ker$T_2 \cup \ldots \cup$ ker$T_n$ where*

$$\ker T_i = \{v^i \in V_i \mid T(v^i) = \bar{0}\},\ i = 1, 2, \ldots, n.$$

*Thus*

$$\ker T = \{(v^1, v^2, \ldots, v^n) \in V_1 \cup V_2 \cup \ldots \cup V_n \mid T(v^1, v^2, \ldots, v^n) = T(v^1) \cup T(v^2) \cup \ldots \cup T(v^n) = 0 \cup 0 \cup \ldots \cup 0\}.$$

It is easily verified that Ker T is a proper n-subgroup of V. Further Ker T is a n-subspace of V.

We will illustrate this situation by the following example.

***Example 2.17:*** Let $V = V_1 \cup V_2 \cup V_3$ be a 3-vector space over Q of 3-dimension (3, 2, 4). Let $W = W_1 \cup W_2 \cup W_3 \cup W_4$ be a 4-vector space over Q of 4-dimension (4, 3, 2, 5). Let $T = T_1 \cup T_2 \cup T_3: V \to W$ be a 3-linear transformation given by

$$T_1: V_1 \to W_4,$$
$$T_1(x_1^1, x_2^1, x_3^1) = (x_1^1 + x_2^1, x_3^1, x_1^1, x_1^1 + x_3^1, x_2^1)$$

for all $x_1^1, x_2^1, x_3^1 \in V_1$,

$$\ker T_1 = \{(x_1^1, x_2^1, x_3^1) \mid T(x_1^1, x_2^1, x_3^1) = (0) \text{ i.e. } x_3^1 = 0,\ x_1^1 = 0,\ x_2^1 = 0 \text{ and } x_1^1 + x_2^1 = 0 \text{ and } x_1^1 + x_3^1 = 0\}$$

Thus ker $T_1$ = {(0, 0, 0)} is the trivial subspace of $V_1$,

$$T_2: V_2 \to W_3,\ T_2(x_1^2, x_2^2) = (x_1^2 + x_2^2, x_1^2)$$

for all $x_1^2, x_2^2 \in V_2$.

$$\ker T_2 = \{(x_1^2, x_2^2) \,|\, T(x_1^2, x_2^2) = (0)\}$$

i.e. $x_1^2 + x_2^2 = 0$ and $x_1^2 = 0$ which forces $x_2^2 = 0$. Thus ker $T_2$ = $\{(0\ 0)\}$.

Now

$$T_3: V_3 \to W_1$$

given by

$$T_3(x_1^3, x_2^3, x_3^3, x_4^3) = (x_1^3 + x_2^3, x_3^3, x_4^3, x_3^3 + x_4^3)$$

for all $x_1^3, x_2^3, x_3^3, x_4^3 \in V_3$.

Now ker $T_3$ gives

$$x_1^3 + x_2^3 = 0,\ x_3^3 = 0,\ x_4^3 = 0,\ x_3^3 + x_4^3 = 0.$$

This gives the condition $x_1^3 = -x_2^3$ and $x_3^3 = x_4^3 = 0$. Thus

$$\ker T_3 = \{(x_1^3, -x_1^3, 0, 0)\}.$$

Thus a subspace of $V_3$. Hence we see the 3-kernel of T is a 1-susbspace of V i.e. $\langle\{(0\ 0\ 0\ 0) \cup (0\ 0) \cup (x_1^3, -x_1^3, 0, 0)\}\rangle$. We can define kernel for any n-linear transformation T be it a usual n-linear transformation or a one to one n-linear transformation. It is easily verified that for any n-vector space $V = V_1 \cup V_2 \cup \ldots \cup V_n$ and any m-vector space $W = W_1 \cup W_2 \cup \ldots \cup W_m$ over the same field F. Suppose $T: V \to W$ is any n-linear transformation from V to W then ker T = ker $T_1 \cup$ ker $T_2 \cup \ldots \cup$ ker $T_n$ would be always a t-subspace of V as each ker $T_i$ is a subspace of $V_i$, i = 1, 2, …, n. It may so happen that some of the ker $T_i$ may be the zero space in such case we will call the subspace of V only as a t-subspace of V where $1 \leq t \leq n$. If all the subspaces given by ker $T_i$ is zero then we call ker T to be the n zero subspace of V; i = 1, 2, …, n.

Now we proceed on to give some more results in case of n-vector spaces and their related linear n-transformation.

**DEFINITION 2.13:** *Let $V = V_1 \cup V_2 \cup \ldots \cup V_n$ be a n-vector space over a field F of type-I. Let $T = T_1 \cup T_2 \cup \ldots \cup T_n : V \to$*

*$V$ be a linear n-transformation of $V$ such that each $T_i : V_i \rightarrow V_i$, i = 1, 2, ..., n. i.e., each $T_i$ is a linear operator on $V_i$ then we define $T = T_1 \cup T_2 \cup ... \cup T_n$ to be a n-linear operator on V. Clearly all n-linear transformations need not be n-linear operator on V. Thus T is a n-linear operator on V if and only if each $T_i$ is a linear operator from $V_i$ to $V_i$, $1 \leq i \leq n$.*

This is the marked difference between the linear operator on a vector space and a n-linear operator on a n-vector space. All n-linear transformations from the n-vector space V to the same n-vector space V need not always be a n-linear operator.

We illustrate this situation by the following example.

***Example 2.18 :*** Let $V = V_1 \cup V_2 \cup V_3 \cup V_4$ be a 4-vector space over Q of 4-dimension (5, 4, 2, 3). Let $T = T_1 \cup T_2 \cup T_3 \cup T_4 : V \rightarrow V$ be a 4-linear transformation given by $T_1: V_1 \rightarrow V_2$, $T_2: V_2 \rightarrow V_3$, $T_3: V_3 \rightarrow V_4$ and $T_4: V_4 \rightarrow V_1$. Clearly none of the linear transformation $T_i$'s are linear operators for they have different domain and range spaces; i = 1, 2, 3, 4. So T though is on the same n-vector space V still T is a linear n-transformation and not a linear n-operator on V, where n = 4.

Suppose we define a 4-linear transformation $P = P_1 \cup P_2 \cup P_3 \cup P_4 : V \rightarrow V$ defined by $P_1: V_1 \rightarrow V_1$, $P_2: V_2 \rightarrow V_2$, $P_3: V_3 \rightarrow V_3$, and $P_4: V_4 \rightarrow V_4$, clearly the 4-linear transformation P is a 4-linear operator of V.

The above example shows the reader that in general a n-linear transformation of a n-vector space V need not in general be a n-linear operator on V. But of course trivially every n-linear operator on V is a n-linear transformation on V.

We have the following result in case of finite n-dimensional n-vector spaces over the field F.

**THEOREM 2.2:** *Let $V = V_1 \cup V_2 \cup ... \cup V_n$ and $W = W_1 \cup W_2 \cup ... \cup W_n$ be any two n-vector spaces over the field F. Let* $B = \{(\alpha_1^1, \alpha_2^1, ..., \alpha_{n_1}^1) \cup (\alpha_1^2, \alpha_2^2, ..., \alpha_{n_2}^2) \cup ... \cup (\alpha_1^n, \alpha_2^n, ..., \alpha_{n_n}^n)\}$ *be a* *n-basis of $V_1 \cup V_2 \cup \ldots \cup V_n$; i.e., $(\alpha_1^i, \alpha_2^i, \ldots, \alpha_{n_i}^i)$ is a basis of $V_i$ $i = 1, 2, \ldots, n$. Let*

$$C = \{(\beta_1^1, \beta_2^1, \ldots, \beta_{n_1}^1) \cup (\beta_1^2, \beta_2^2, \ldots, \beta_{n_2}^2) \cup \ldots \cup (\beta_1^n, \beta_2^n, \ldots, \beta_{n_n}^n)\}$$

*be any n-vector in $W = W_1 \cup W_2 \cup \ldots \cup W_n$ then there is precisely only one linear n-transformation $T = T_1 \cup T_2 \cup \ldots \cup T_n$ from V on to W such that $T\alpha_j^i = \beta_j^i$, $j = 1, 2, \ldots, n_i$, $1 \leq i \leq n$.*

*Proof:* To prove that there is some n-linear transformation T with T(B) = C, it is enough if we show for the $T = T_1 \cup T_2 \cup \ldots \cup T_n$ we have $T_i\alpha_j^i = \beta_j^i$, i = 1, 2, …, n and j = 1, 2, …, $n_i$. Given $(\alpha_1^i, \alpha_2^i, \ldots, \alpha_{n_i}^i)$ in $V_i$ there is a unique $n_i$ tuple $(x_1^i, x_2^i, \ldots, x_{n_i}^i)$ such that $\alpha^i = x_1^i\alpha_1^i + x_2^i\alpha_2^i + \ldots + x_{n_i}^i\alpha_{n_i}^i$, for this vector $\alpha^i$ we define

$$T_i\alpha^i = x_1^i\beta_1^i + x_2^i\beta_2^i + \ldots + x_{n_i}^i\beta_{n_i}^i$$

true for each i; i = 1, 2, …, n.

Clearly $T_i$ is a well defined rule for associating with each vector $\alpha^i$ in $V_i$ a vector $T_i\,\alpha^i$ in $W_i$. From the definition it is clear that $T_i\alpha_j^i = \beta_j^i$ for each j. To see that $T_i$ is linear; let $\beta^i$ = $y_1^i\alpha_1^i + y_2^i\alpha_2^i + \ldots + y_{n_i}^i\alpha_{n_i}^i$ be in $V_i$ and $c_i$ be any scalar.

$$T(c_i\alpha^i + \beta^i) = (c_ix_1^i + y_1^i)\beta_1^i + \ldots + (c_ix_{n_i}^i + y_{n_i}^i)\beta_{n_i}^i.$$

On the other hand,

$$c_i(T_i(\alpha^i)) + T_i(\beta^i) = c_i\sum_{j=1}^{n_i} x_j^i\beta_j^i + \sum_{j=1}^{n_i} y_j^i\beta_j^i$$

$$= \sum_{j=1}^{n_i}(c_ix_j^i + y_j^i)\beta_j^i$$

and thus $T_i(c_i\alpha^i + \beta^i) = c_i(T_i\alpha^i) + T_i\,\beta^i$ true for each i; i = 1, 2, …, n. If $U = U_1 \cup U_2 \cup \ldots \cup U_n$ is a linear n-transformation from V on to W with $U_i\alpha^i_j = \beta^i_j$, j = 1, 2, …, $n_i$ and true for each i; i = 1, 2, … , n. then for the vector $\alpha^i = \sum_{j=1}^{n_i} x^i_j \alpha^i_j$ we have

$$U_i\alpha^i = U_i(\sum_{j=1}^{n_i} x^i_j\alpha^i_j)$$
$$= \sum_{j=1}^{n_i} x^i_j U_i \alpha^i_j$$
$$= \sum_{j=1}^{n_i} x^i_j \beta^i_j ,$$

true for each and every i; i = 1, 2, .., n. Thus $U_i$ is exactly the rule $T_i$ , i = 1, 2, … , n hence U is exactly the rule T which we have defined. This shows the n-linear transformation $T = T_1 \cup T_2 \cup \ldots \cup T_n$ is unique.

Having defined n-kernel of a n-linear transformation T we now proceed on to define the n-range of the n-linear transformation $T = T_1 \cup T_2 \cup \ldots \cup T_n$.

**DEFINITION 2.14:** *Let $T = T_1 \cup T_2 \cup \ldots \cup T_n$ be a n-linear transformation from the n-vector space $V = V_1 \cup V_2 \cup \ldots \cup V_n$ in to another m-vector space W, m > n. The range of T is called the n-range of T denoted by $R^n_T$, is a p-subspace of W p < m that is $R^n_T = \{\beta = \beta_1 \cup \beta_2 \cup \ldots \cup \beta_m \in W\}$ such that $\beta = T(\alpha)$ for some $\alpha = \alpha_1 \cup \alpha_2 \cup \ldots \cup \alpha_n$ in V. Clearly if β, γ ∈ $R^n_T$ and c any scalar, then there are n-vectors α, δ in V such that Tα = β and Tβ = γ. Since T is n-linear.*

$$T(c\alpha + \delta) = cT\alpha + T\delta$$
$$= c\beta_1 + \beta_2$$

*which is in $R^n_T$. Now V and W be any two n-vector space and m vector space respectively defined over the field F and let $T = T_1 \cup T_2 \cup \ldots \cup T_n$ be a linear n-transformation from V into W. The n-null space T is a n set of all n-vectors α in V such that*

$$\begin{aligned} T\alpha &= T(\alpha_1 \cup \alpha_2 \cup \ldots \cup \alpha_n) \\ &= (T_1 \cup T_2 \cup \ldots \cup T_n)\,(\alpha_1 \cup \alpha_2 \cup \ldots \cup \alpha_n) \\ &= T_1\alpha_1 \cup T_2\alpha_2 \cup \ldots \cup T_n\alpha_n \\ &= 0 \cup 0 \cup \ldots \cup 0. \end{aligned}$$

*If the n-vector space V is finite n-dimensional, the n-rank of T is the n-dimension of the n-range of T and will vary depending on the nature of the n-linear transformation like if T is a shrinking n-linear transformation, it would be different and so on.*

Now we can prove the most important theorem relating the n-rank of T and n-nullity of T for a n-linear transformation only as for other n-linear transformation like shrinking n-linear transformation the result in general may not be true.

**THEOREM 2.3:** *Let V and W be two n-vector space and m-vector space over the field F, m > n and let T is be linear n-transformation i.e.* $T = T_1 \cup T_2 \cup \ldots \cup T_n$ *from V to W is such that* $T_i: V_i \to W_j$ *and the* $W_j$*'s are distinct spaces for each* $T_i$*, i.e. no two subspaces of V are mapped on to the same subspace in W. Suppose V is* $(n_1, n_2, \ldots, n_n)$ *finite dimensional, then n rank T + n nullity T = n dim V.*

*Proof:* Given $V = V_1 \cup V_2 \cup \ldots \cup V_n$ is a n-vector space over F and $W = W_1 \cup W_2 \cup \ldots \cup W_n$ is a m-vector space over F ($m > n$) of dimensions $(n_1, n_2, \ldots, n_n)$ and $(m_1, m_2, \ldots, m_n)$ respectively. $T = T_1 \cup T_2 \cup \ldots \cup T_n$ is a n-linear transformation such that each $T_i$ is a linear transformation from $V_i$ to a unique $W_j$, i.e. no two vector spaces $V_i$ and $V_k$ can be mapped to same $W_j$, if $i \neq k$; $1 \leq i, k \leq n$ and $1 \leq j \leq m$ . Now n-rank T + n nullity T = n dim W
i.e. n-rank $(T_1 \cup T_2 \cup \ldots \cup T_n)$ + n nullity of $(T_1 \cup T_2 \cup \ldots \cup T_n)$ = n dim $(V_1 \cup V_2 \cup \ldots \cup V_n)$.
i.e. rank $T_1 \cup$ rank $T_2 \cup \ldots \cup$ rank $T_n$ + nullity $T_1 \cup$ nullity $T_2 \cup \ldots \cup$ nullity $T_n$ = (dim $V_1$, dim $V_2$, …, dim $V_n$) = $(n_1, n_2, \ldots, n_n)$.

Suppose $N = N_1 \cup N_2 \cup \ldots \cup N_n$ be the p-null space of the n-space V; $0 \leq p \leq n$. Let

$$\alpha = \left\{\left(\alpha_1^1, \alpha_2^1, \ldots, \alpha_{k_i}^1\right) \cup \left(\alpha_1^2, \alpha_2^2, \ldots, \alpha_{k_2}^2\right) \cup \ldots \cup \left(\alpha_1^n, \alpha_2^n, \ldots, \alpha_{k_n}^n\right)\right\}$$

be a n-basis for N. Here $0 \le k_i \le n_i$ ; i = 1, 2, … , n. If $k_i = 0$ then the corresponding null space is the zero space. Now we show the working for any i; $T_i: V_i \rightarrow W_j$. This result which we would prove is true for all i = 1, 2, …, n.

Let $\left\{\alpha_1^i, \alpha_2^i, ..., \alpha_k^i\right\}$ be a basis for $N_i$, the null space of $T_i$. There are vectors $\alpha_{k_i+1}^i, ..., \alpha_{n_i}^i$ in $V_i$ such that $\left\{\alpha_1^i, \alpha_2^i, ..., \alpha_{n_i}^i\right\}$ is a basis for $V_i$; true for each i; i = 1, 2, …, n. We shall now prove that $\left\{T_i\alpha_{k_i+1}, ..., T_i\alpha_{n_i}\right\}$ is a basis for the range of $T_i$. The vectors $T_i\alpha_1^i, T_i\alpha_2^i, ..., T_i\alpha_n^i$ certainly span the range of $T_i$ and since $T_i\alpha_j^i = 0$ for $j \le k_i$ we see that $T_i\alpha_{k_i+1}, ..., T_i\alpha_{n_i}$ span the range, to see that these vectors are linearly independent, suppose we have scalars $c_j$'s such that

$$\sum_{j=k_i+1}^{n_i} c_j T_i(\alpha_j^i) = 0 .$$

This says that

$$T_i(\sum_{j=k_i+1}^{n_i} c_j \alpha_j^i) = 0$$

and accordingly the vector

$$\alpha^i = \sum_{j=k_i+1}^{n_i} c_j \alpha_j^i$$

is in the null space of $T_i$. Since $\alpha_1^i, \alpha_2^i, ..., \alpha_{n_i}^i$ form a basis for $N_i$ there must be scalars $\beta_1^i, \beta_2^i, ..., \beta_{n_i}^i$ such that $\alpha^i = \sum_{j=1}^{k_i} \beta_j^i \alpha_j^i$ . Thus

$$\sum_{j=1}^{k_i} \beta_j^i \alpha_j^i - \sum_{j=k_i+1}^{n_i} \beta_j^i \alpha_j^i = 0.$$

Since $\alpha_1^i, \alpha_2^i, ..., \alpha_{n_i}^i$ are linearly independent we must have $b_1^i$ = … = $b_{k_i}^i$ = $c_{k_{i+1}}^i$ = … = $c_{n_i}^i$ = 0. If $r_i$ is the rank of $T_i$ the fact that $T_i\ \alpha_{k_{i+1}}^i$ , . . . , $T_i\ \alpha_{n_i}^i$ form a basis, for the range of $T_i$ tells us that $r_i = n_i - k_i$. Since $k_i$ is the nullity of $T_i$ and $n_i$ is the dimension of

$V_i$, we get rank $T_i$ + nullity $T_i$ = dim. $V_i$. This is true for each and every i. That is

(rank $T_1$ + nullity $T_1$) $\cup$ (rank $T_2$ + nullity $T_2$) $\cup \ldots \cup$ (rank $T_n$ + nullity $T_n$)

$$= \dim (V_1 \cup V_2 \cup \ldots \cup V_n)$$

i.e., (rank $T_1 \cup$ rank $T_2 \cup \ldots \cup$ rank $T_n$) + (nullity $T_1 \cup$ nullity $T_2 \cup \ldots \cup$ nullity $T_n$)

$$= \dim (V_1 \cup V_2 \cup \ldots \cup V_n) \text{ that is}$$

rank $(T_1 \cup T_2 \cup \ldots \cup T_n)$ + nullity $(T_1 \cup T_2 \cup \ldots \cup T_n)$

$$= (n_1, n_2, \ldots, n_n).$$

n rank T + n nullity T = n dim V.

Now in the relation

$$\text{n rank } T + \text{n nullity } T = \text{n dim } (V) = (n_1, n_2, \ldots, n_n).$$

We assume the n-linear transformation is such that it is not shrinking it is a n-linear transformation given in definition 2.12. Also we see if nullity $T_i = 0$ for some i in such cases we have rank $T_i = \dim V_i$. Since a p-nullspace in general need not always be a nontrivial subspace we may have the p-nullspace of the n-vector space be such that $p < n$.

Now we proceed on to the algebra of n-linear transformations. Let us assume V and W are two n-vector space and m-vector space respectively defined over the field K.

**THEOREM 2.4:** *Let V and W be any two n-vector space and m-vector space respectively defined over the field F(m > n). Let T and U be n-linear transformations as given in definition from V into W. The n-function (T + U) defined by (T + U)α = Tα + Uα is a n-linear transformation from V into W, if c is any element from F, the function cT defined by (cT)α = cTα is a n-linear transformation from V into W.*

*The set of all n-linear transformations from V into W with addition and scalar multiplication defined above is an n-vector space over the field F.*

*Proof:* Let $V = V_1 \cup V_2 \cup \ldots \cup V_n$ be a n-vector space over F and $W = W_1 \cup W_2 \cup \ldots \cup W_m$ (m>n) a m-vector space over F. $T = T_1 \cup T_2 \cup \ldots \cup T_n$ a n-linear transformation from V to W.

If $U = U_1 \cup U_2 \cup \ldots \cup U_n$ is a n-linear transformation from V into W; define the n-function (T + U) for $\alpha = \alpha_1 \cup \alpha_2 \cup \ldots \cup \alpha_n \in V$ by $(T + U)\ \alpha = T\alpha + U\alpha$ then (T + U) is a n-linear transformation of V into W.

$(T + U)\ (c\alpha + \beta)$

$= [(T_1 \cup T_2 \cup \ldots \cup T_n) + (U_1 \cup U_2 \cup \ldots \cup U_n)]\ [c\alpha + \beta]$

$= [(T_1 + U_1) \cup (T_2 + U_2) \cup \ldots \cup (T_n + U_n)]$
$[c(\alpha_1 \cup \alpha_2 \cup \ldots \cup \alpha_n) + (\beta_1 \cup \beta_2 \cup \ldots \cup \beta_n))]$

$= [(T_1 + U_1) \cup (T_2 + U_2) \cup \ldots \cup (T_n + U_n)]$
$[(c\alpha_1 + \beta_1) \cup (c\alpha_2 + \beta_2) \cup \ldots \cup (c\alpha_n + \beta_n)]$

$= [(T_1 + U_1)\ (c\alpha_1 + \beta_1)] \cup [(T_2 + U_2)\ (c\alpha_2 + \beta_2)] \cup \ldots \cup [(T_n + U_n)\ (c\alpha_n + \beta_n)]$.

Now using the properties of linear transformation on linear vector space we get $(T_i + U_i)\ (c\alpha_i + \beta_i) = c\ (T_i + U_i)\ (\alpha_i) + (T_i + U_i)\ (\beta_i)$ for each i = 1, 2, …, n.

Thus $(T + U)\ (c\alpha + \beta) = \{[c(T_1 + U_1)\ \alpha_1 \cup c(T_2 + U_2)\ \alpha_2 \cup \ldots \cup c(T_n + U_n)\ \alpha_n] + (T_1 + U_1)\ \beta_1 \cup (T_1 + U_1)\ \beta_2 \cup \ldots \cup (T_n + U_n)\ \beta_n\} = c(T + U)\ \alpha + (T + U)\beta$, which shows (T + U) is a n-linear transformation from V into W.

$cT(d\alpha + \beta)$

$= c[(T_1 \cup T_2 \cup \ldots \cup T_n)\ [d(\alpha_1 \cup \alpha_2 \cup \ldots \cup \alpha_n) + (\beta_1 \cup \beta_2 \cup \ldots \cup \beta_n)]$

$= c[T_1 \cup T_2 \cup \ldots \cup T_n]\ [(d\alpha_1 + \beta_1) \cup (d\alpha_2 + \beta_2) \cup \ldots \cup (d\alpha_n + \beta_n)]$

$= c[T_1(d\alpha_1 + \beta_1) \cup T_2(d\alpha_2 + \beta_2) \cup \ldots \cup T_n(d\alpha_n + \beta_n)]$

$= T_1(c(d\alpha_1 + \beta_1)) \cup T_2(c(d\alpha_2 + \beta_2)) \cup \ldots \cup T_n(c(d\alpha_n + \beta_n)$
(since each $cT_i$ is a linear transformation)

$= T[c(d\alpha + \beta)]$

$= c[T(d\alpha + \beta)]$

$= d[(cT)\alpha] + (cT)\beta$.

This shows cT is a n-linear transformation of V into W.

**THEOREM 2.5:** *Let V be a n-vector space of n-dimension ($n_1$, $n_2$, ..., $n_n$) over the field F, and let W be a m dimensional m-vector space over the same field F with m-dimension ($m_1$, $m_2$, ..., $m_n$) (m > n). Then $L^n(V,W)$ is the finite dimensional n-space over F of n-dimension $\left(m_{i_1}n_1, m_{i_2}n_2, \ldots, m_{i_n}n_n\right)$ where $L^n(V,W)$ denotes the space of all n-linear transformations of V into W, $1 \le i_1, i_2, \ldots, i_n \le m$.*

*Proof:* Let

$B = \{(\alpha_1^1, \alpha_2^1, \ldots, \alpha_{n_1}^1) \cup (\alpha_1^2, \alpha_2^2, \ldots, \alpha_{n_2}^2) \cup \ldots \cup (\alpha_1^n, \alpha_2^n, \ldots, \alpha_{n_n}^n)\}$

be a n-basis for $V = V_1 \cup V_2 \cup \ldots \cup V_n$ of the n-vector space of n-dimension ($n_1$, $n_2$, … , $n_n$).

Let

$C = \{(\beta_1^1, \beta_2^1, \ldots, \beta_{m_1}^1) \cup (\beta_1^2, \beta_2^2, \ldots, \beta_{m_2}^2) \cup \ldots \cup (\beta_1^m, \beta_2^m, \ldots, \beta_{m_n}^m)\}$

be a m-basis of the m-vector space $W = W_1 \cup W_2 \cup \ldots \cup W_m$ of m-dimension ($m_1$, $m_2$, …, $m_n$).

Let $L^n(V, W)$ be the set of all n-linear transformation of V into W. For every pair of integers ($p^j$, $q^i$), $1 \le j \le m$ and $1 \le i \le n$, $1 \le p^j \le m_j$ and $1 \le q^i \le n_i$, we define a linear transformation $E^{p^j,q^i}$ ; $1 \le i \le n$ and $k_j \le m$ of $V_i$ into $W_j$ by

$$E^{p^j,q^i}(\alpha_t^i) = \begin{cases} 0 \text{ if } t \ne q^i \\ \beta_{p^j}^j \text{ if } t = q^i \end{cases} = \delta_{tq^i}\beta_{p^j}^j .$$

By the theorem $(T_j\alpha_j^i = \beta_i^j)$ their is a unique linear transformation from $V_i$ into $W_j$. We claim that $m_jn_i$ transformation $E^{p^j,q^i}$ form a basis for $L_i(V_i, W_j)$. This is true for each i, i = 1, 2, …, n and the appropriate j, $1 \le j \le m$ with no two spaces $V_i$ of V mapped into the same $W_j$. Let $T_i$ be a linear transformation from $V_i$ into $W_j$, $1 \le i \le n$, $1 \le j \le m$. For each $k_i \le k \le n_i$. Let $A_{1k}$, …, $A_{m_jk}$ be the coordinates of the vector $T_i\alpha_k^i$ in the ordered basis $(\beta_1^j, \beta_2^j, \ldots, \beta_{m_j}^j)$ the n-basis of W given in C.

$$T_i\, \alpha^i_k = \sum_{p=1}^{m_j} A_{p^k} \beta^j_p \qquad (1)$$

We wish to show that

$$T_i = \sum_{p^j=1}^{m_j} \sum_{q^i=1}^{n_i} A_{p^j q^i} E^{p^j, q^j} \qquad (2)$$

Let $U_i$ be the linear transformation in the right hand member of (2). Then for each k

$$U_i \alpha^i_k = \sum_{p^j=1} \sum_{q^i=1} A_{p^j q^i} E^{p^j, q^i} (\alpha^i_k)$$

$$= \sum_{p^j} \sum_{q^i} A_{p^j q^i} \delta_{k q^i} \beta^j_{p^j}$$

$$= \sum_{p^j=1}^{m_j} A_{p^j k} \beta^j_{p^j}$$

$$= T_i \alpha^i_k \,;$$

and consequently $U_i = T_i$. Now from 2 we see $E^{p^j, q^i}$ spans $L_i(V_i, W_j)$. We must only now show they form a linearly independent set. This is very clear from the fact

$$U_i = \sum_{p^j} \sum_{q^i} A_{p^j q^i} E^{p^j, q^i}$$

is the zero transformation, then $U_i \alpha^i_k = 0$ for each k, so that

$$\sum_{p^j=1}^{m_j} A_{p^j k} \beta^j_{p^j} = 0$$

and thus the independence of the $\beta^j_{p^j}$ implies that $A_{p^j k} = 0$ for every $p^j$ and k. Since this is true of every i, i = 1, 2, … , n. We have

$$L^n(V,W) = L^n_1(V_1, W_{i_1}) \cup L^n_2(V_2, W_{i_2}) \cup \ldots \cup L^n_n(V_n, W_{i_n});$$

where $i_1, i_2, \ldots, i_n$ are distinct elements from the set {1, 2, …,m} and m > n. Hence $L^n(V,W)$ is a n-space of dimension $(m_{i_1} n_1, m_{i_2} n_2, \ldots, m_{i_n} n_n)$ over the same field F. This n-space will be known as the n-space of n-linear transformation of the n-vector space $V = V_1 \cup V_2 \cup \ldots \cup V_n$ of n-dimension $(n_1, n_2, \ldots, n_n)$ into the m-vector space $W = W_1 \cup W_2 \cup \ldots \cup W_m$ of m-dimension $(m_1, m_2, \ldots, m_n)$, $m > n$.

Now having proved that the space of all n-linear transformations of a n-vector space V into a m-vector space W forms a n-vector space over the same field F, we prove another interesting theorem.

**THEOREM 2.6:** *Let V and W be two n-vector spaces of n-dimensions $(n_1, n_2, \ldots, n_n)$ and $(t_1, t_2, \ldots, t_n)$ respectively defined over the field F. Z be a m-vector space defined over the same field F(m > n). Let T be a n-linear transformation of V into W and U be a n linear transformation from W into Z. Then the composed function UT defined by $(UT)(\alpha) = U(T(\alpha))$ is a n-linear transformation from V into Z, $\alpha \in V$.*

*Proof*: Given $V = V_1 \cup V_2 \cup \ldots \cup V_n$ and $W = W_1 \cup W_2 \cup \ldots \cup W_n$ are 2 n-vector spaces over F. $Z = Z_1 \cup Z_2 \cup \ldots \cup Z_m$ is given to be a m-vector space over F, $m > n$. $T: V \to W$ is a n-linear transformation; that is $T = T_1 \cup T_2 \cup \ldots \cup T_n : V \to W$ with $T_i: V_i \to W_j$ and no two vector spaces in V are mapped into the same vector space $W_j$, $i = 1, 2, \ldots, n$ and $1 \leq j \leq n$.
Now $U = U_1 \cup \ldots \cup U_n: W \to Z$ is a n-linear transformation such that $U_j: W_j \to Z_k$, $j = 1, 2, \ldots, n$ and $1 \leq k \leq m$ such that no two subspaces of W are mapped into the same $Z_k$.

Now

$$
\begin{aligned}
(U_j\, T_i)\,(c\alpha^i + \beta^i) &= U_j\,[T_i\,(c\alpha^i + \beta^i)\\
&= U_j\,[T_i\,(c\alpha^i) + T\,(\beta^i)]\\
&= U_j\,[c\, T_i\,(\alpha^i) + T_i\,(\beta^i)]\\
&= U_j\,[c\omega^j + \delta^j\,]
\end{aligned}
$$

(as $T_i : V_i \to W_j\,;\ \omega^j, \delta^j \in W_j$)

$$
\begin{aligned}
&= c\, U_j\,(\omega^j) + U_j\,(\delta^j)\\
&= ca^k + b^k;\ a^k, b^k \in Z_k.
\end{aligned}
$$

Thus $U_jT_i$ is a n-linear transformation from $W_j$ to $Z_k$. Hence the claim; for the result is true for each i and each j. Thus UT is a n-linear transformation from W to Z.

So

$$\begin{aligned} U \text{ o } T &= (U_1 \cup U_2 \cup \ldots \cup U_n) \text{ o } (T_1 \cup T_2 \cup \ldots \cup T_n) \\ &= U_1 T_{i_1} \cup U_2 T_{i_2} \cup \ldots \cup U_n T_{i_n} \end{aligned}$$

$(i_1, i_2, \ldots , i_n)$ is a permutation of 1, 2, 3, … , n. Now we for the notational convenience recall that if $V = V_1 \cup V_2 \cup \ldots \cup V_n$ is a n vector space over a field F then $V_i$'s are called as component subvector spaces of V. $V_i$ is also unknown as the component of V.

Now we proceed on to define the notion of linear n-operator.

**DEFINITION 2.15:** *Let $V = V_1 \cup V_2 \cup \ldots \cup V_n$ be a n-vector space over F, a n-linear operator on V is a n-linear transformation T from V to V, such that $T = T_1 \cup T_2 \cup \ldots \cup T_n$ with $T_i: V_i \to V_i$ for $1 \leq i \leq n$. Thus in the above theorem not only $V = W = Z$ but U and T are such that $T_i$: $V_i \to V_i$; $U_i$: $V_i \to V_i$ so that U and T are n-linear operators on the n space V, we see composition UT is again a n-linear operator on V.*

*Thus the n-space $L^n (V, V)$ has a multiplication defined as composition. In this case the operator TU is also defined. In general $TU \neq UT$ i.e., $UT - TU \neq 0$.*

*Now $L^n (V, V)$ would be only a n-vector space of dimension $(n_1^2, n_2^2, \ldots, n_n^2)$, n-dimension of V is $(n_1, n_2, \ldots , n_n)$.*

$$\begin{aligned} UT &= (U_1 \cup U_2 \cup \ldots \cup U_n) \circ (T_1 \cup T_2 \cup \ldots \cup T_n) \\ &= (U_1T_1 \cup U_2 T_2 \cup \ldots \cup U_nT_n). \\ TU &= (T_1 \cup T_2 \cup \ldots \cup T_n) \circ (U_1 \cup U_2 \cup \ldots \cup U_n) \\ &= T_1U_1 \cup T_2U_2 \cup \ldots \cup T_nU_n. \end{aligned}$$

*Here $T_i : V_i \to V_i$ and $U_i : V_i \to V_i$, i = 1, 2, … , n.*

*Now only in this case $T^2 = TT$ and in general $T^n = TT \ldots T$; n times for n = 1, 2, …, n. We define $T^o = I_1 \cup I_2 \cup \ldots \cup I_n$ = identity n-function of $V = V_1 \cup V_2 \cup \ldots \cup V_n$. It may so happen*

*depending on each $V_i$ we will have different power of $T_i$ to be approaching identity for varying linear transformation. i.e. if T : $T_1 \cup T_2 \cup \ldots \cup T_n$ on $V = V_1 \cup V_2 \cup \ldots \cup V_n$, such that $T_i: V_i \to V_i$ (only) for i = 1, 2, …, n since n-dimension of V is $(n_1, n_2, \ldots, n_n)$. so T o T = $T^2 = (T_1 \cup T_2 \cup \ldots \cup T_n)(T_1 \cup T_2 \cup \ldots \cup T_n) = (T_1^2, T_2^2, \ldots, T_n^2)$. Like wise any power of T. $I = I_1 \cup I_2 \cup \ldots \cup I_n$ is the identity function on V i.e. each $I_i : V_i \to V_i$ is such that $I_i(\mu^i) = \mu^i$ for all $\mu^i \in V_i$ ; i = 1, 2, …, n. Only under these special conditions we define $L^n$ (V, V); elements of Ln(v, v) are called special n-linear operators.*

**LEMMA 2.1:** *Let $V = V_1 \cup V_2 \cup \ldots \cup V_n$ be a n-vector space over the field F; let U, $T_1$ and $T_2$ be n-linear operators on V; let c be an element of F*

a. *IU = UI = U where $I = I_1 \cup I_2 \cup \ldots \cup I_n$ is the n-identity transformation*
b. *$U(T_1 + T_2) = UT_1 + UT_2$*
   *$(T_1 + T_2)U = T_1U + T_2U$*
c. *$C(U T_1) = (C U) T_1 = U (C T_1)$.*

*Proof:* Given $V = V_1 \cup V_2 \cup \ldots \cup V_n$ be a n-vector space over F, F a field. $I = I_1 \cup I_2 \cup \ldots \cup I_n$ be the n-identity transformation of V to V i.e. $I_j: V_j \to V_j$; is the identity transformation of each $V_j$, j = 1, 2, …, n. $U = U_1 \cup U_2 \cup \ldots \cup U_n: V \to V$ such that $U_i: V_i \to V_i$ for i = 1, 2, …, n. $T_i$ . $T_1^i \cup T_2^i \cup \ldots \cup T_n^i : V \to V$ such that $T_j^i : V_j \to V_j$; j = 1, 2, …, n and i = 1, 2.

$$U = (U_1 \cup U_2 \cup \ldots \cup U_n)(I_1 \cup I_2 \cup \ldots \cup I_n)$$
$$= U_1 \circ I_1 \cup U_2 \circ I_2 \cup U_n I_n$$
$$= U_1 \cup \ldots \cup U_n.$$

Further

$$IU = (I_1 \cup I_2 \cup \ldots \cup I_n) \circ (U_1 \cup U_2 \cup \ldots \cup U_n)$$
$$= I_1 \circ U_1 \cup I_2 \circ U_2 \cup I_n \text{ o } U_n$$
$$= U_1 \cup U_2 \cup \ldots \cup U_n.$$

Thus IU=UI.

$$
\begin{aligned}
U(T_1 + T_2) &= UT_1 + UT_2 \\
&= [U_1 \cup U_2 \cup \ldots \cup U_n] [(T_{11} \cup T_{12} \cup \ldots \cup T_{1n}) + (T_{21} \cup T_{22} \cup \ldots \cup T_{2n})] \\
&= [U_1 \cup U_2 \cup \ldots \cup U_n] (T_{11} + T_{21}) \cup (T_{12} + T_{22}) \cup \ldots \cup (T_{1n} + T_{2n}).
\end{aligned}
$$

We know from the results in linear algebra $U(T_1 + T_2) = UT_1 + UT_2$ for any $U, T_1, T_2 \in L(V_1 ; V_1)$ where $V_1$ is a vector space and $L(V_1, V_1)$ is the collection of all linear operators from $V_1$ to $V_1$.

Now in $U_i (T_{1i} + T_{2i})$, $U_i$ $T_{1i}$ and $T_{2i}$ are linear operators from $V_i$ to $V_i$ true for each $i = 1, 2, \ldots, n$. Thus $U(T_1 + T_2) = UT_1 + UT_2$ and $(T_1 + T_2) U = T_1U + T_2U$. Further $C(UT_1) = (CU) T_1 = U(CT_1)$ for all $U, T_1 \in L^n(V,V)$. Let $U = U_1 \cup U_2 \cup \ldots \cup U_n$ and $T_1 = (T_1^1 \cup T_2^1 \cup \ldots \cup T_n^1)$ where $U_i$: $V_i \to V_i$ for each i and $T_i^1 : V_i \to V_i$ for each $i = 1, 2, \ldots, n$.

$$
\begin{aligned}
&C[(U_1 \cup U_2 \cup \ldots \cup U_n)(T_1^1 \cup T_2^1 \cup \ldots \cup T_n^1)] \\
&\quad= C[U_1T_1^1 \cup U_2T_2^1 \cup \ldots \cup U_nT_n^1] \\
&\quad= (CU_1T_1^1 \cup CU_1T_2^1 \cup \ldots \cup CU_nT_n^1).
\end{aligned}
$$

$$(CU_1 \cup CU_2 \cup \ldots \cup CU_n) (T^1{}_1 \cup T^1{}_2 \cup \ldots \cup T^1{}_n) = (CU)T_1.$$

But

$$
\begin{aligned}
C(UT_1) &= (CU_1 \cup CU_2 \cup \ldots \cup CU_n) (T^1{}_1 \cup T^1{}_2 \cup \ldots \cup T^1{}_n) \\
&= (U_1 \cup U_2 \cup \ldots \cup U_n) (CT_1) \\
&= (U_1 \cup U_2 \cup \ldots \cup U_n) (CT_1^1 \cup CT_2^1 \cup \ldots \cup CT_n^1) \\
&= U_1(CT_1^1) \cup U_2(CT_2^1) \cup \ldots \cup U_n(CT_n^1) \\
&= (U_1 \cup U_2 \cup \ldots \cup U_n) (CT_1^1 \cup CT_2^1 \cup \ldots \cup CT_n^1) \\
&= U(CT_1).
\end{aligned}
$$

Let us denote the set of all n-linear transformation from V to V; this will also include the set of all n-linear operator $T = T_1 \cup T_2 \cup \ldots \cup T_n$ with $T_i : V_i \to V_i$, $i = 1, 2, \ldots, n$. Let us denote the n-linear transformation on V by $L_T^n$ (V, V). Clearly $L^n$ (V, V) $\subseteq$ $L_T^n$ (V, V). This is the marked difference between the usual linear operator and n-linear operator. For a n-linear transformation can be n-linear operator or n-linear transformation. But every linear operator from V to V is always a linear transformation.

Let V to W be two n-linear vector spaces of same dimension say $(n_1, n_2, \ldots, n_n)$ and $(n_{i_1}, n_{i_2}, \ldots, n_{i_n})$ where $(i_1, i_2, \ldots, i_n)$ is a permutation of (1, 2, …, n).

Let $T_s$: V → W be a n-linear transformation where $T_s = T_1 \cup T_2 \cup \ldots \cup T_n$; $T_i$: $V_i \to W_j$ where $W_j$ is such that dim $V_i$ = dim $W_j$, this is the way every $V_i$ is matched. This will certainly happen because the n-dimension of both V and W are one and the same. We call such n-linear transformation from same dimensional space V into W satisfying the conditions mentioned by each $T_i$; i = 1, 2, …, n. denoted by $T_s$, for this is a special n-linear transformation.

If each $T_i$ in $T_s$; i = 1, 2, …, n is invertible; then we can find a special n-linear transformation $U_s$ : W → V such that $T_sU_s = U_sT_s$ and is the identity function on W. If $T_s$ is invertible the function $U_s$ is unique and is denoted by $T_s^{-1}$. Further more $T_s$ is 1-1 that is $T_s\alpha = T_s\beta$ implies $\alpha = \beta$ where $\alpha = \alpha_1 \cup \alpha_2 \cup \ldots \cup \alpha_n$ and $\beta = \beta_1 \cup \beta_2 \cup \ldots \cup \beta_n$. $T_s$ is onto, that is the range of $T_s$ is all of W.

**THEOREM 2.7:** *Let V and W be n-vector spaces over the field F of same dimension $(n_1, n_2, \ldots, n_n)$ over the field F. If $T_s$ is a special n-linear transformation from V into W and $T_s$ is invertible then the inverse function $T_s^{-1}$ is a special n-linear transformation from W into V.*

*Proof:* Let $T_s = T_1 \cup T_2 \cup \ldots \cup T_n$ be a special n-linear transformation from the same n-dimensional spaces V into W, where n-dimension of V is $(n_1, n_2, \ldots, n_n)$ and that of W is $\left(n_{i_1}, n_{i_2}, \ldots, n_{i_n}\right)$; $(i_1, i_2, \ldots, i_n)$ a permutation of (1, 2, 3, …, n).

i.e., $T_s$: V $\to$ W; $T_i$: $V_i \to W_j$ where dim $V_i$ = dim $W_j$. $T_s^{-1} = T_1^{-1} \cup T_2^{-1} \cup \ldots \cup T_n^{-1}$ is the inverse of $T_s$.

Let $\beta_1$, $\beta_2$ be vectors in W let C $\in$ F. To show $T_s^{-1}$ ($C\beta_1$ + $\beta_2$) = $C\,T_1^{-1}$ $\beta_1$ + $T_1^{-1}\beta_2$ where $\beta_1 = \beta_1^1 \cup \beta_2^1 \cup \ldots \cup \beta_n^1$ and $\beta_2 = \beta_1^2 \cup \beta_2^2 \cup \ldots \cup \beta_n^2$.

$T_s^{-1}$ $(C\beta_1 + \beta_2)$

$= T_s^{-1}\left(C\beta_1^1 + \beta_1^2 \cup C\beta_2^1 + \beta_2^2 \cup \ldots \cup C\beta_n^1 + \beta_n^2\right)$

$= \left(T_1^{-1}(C\beta_1^1 + \beta_1^2) \cup T_2^{-1}(C\beta_2^1 + \beta_2^2) \cup \ldots \cup T_n^{-1}(C\beta_n^1 + \beta_n^2)\right)$

$= (CT_1^{-1}\beta_1^1 + T_1^{-1}\beta_1^2) \cup (CT_2^{-1}\beta_2^1 + T_2^{-1}\beta_2^2) \cup \ldots \cup (CT_n^{-1}\beta_n^1 + T_n^{-1}\beta_n^2)$

$= C\,T_s^{-1}\beta_1 + T_s^{-1}\beta_2.$

Let $\alpha_i = C\,T_s^{-1}\beta_i$ ; i = 1, 2, that is let $\alpha_i$ be the unique n-vector in the V such that $T_s\alpha_i = \beta_i$. Since $T_s$ is n-linear;

$$T_s(C\alpha_1 + \alpha_2) = CT_s\alpha_1 + T_s\alpha_2 = C\beta_1 + \beta_2.$$

Thus $C\alpha_1 + \alpha_2$ is the unique n-vector in V which is sent by $T_s$ into $C\beta_1 + \beta_2$ and so

$$T_s^{-1}\ (C\beta_1 + \beta_2) = C\alpha_1 + \alpha_2 = C(T_s^{-1}\beta_1) + T_s^{-1}\beta_2$$

and $T_s^{-1}$ is n-linear, the proof is similar to the earlier one using $T_s = T_1 \cup T_2 \cup \ldots \cup T_n$ and $T_s^{-1} = T_1^{-1} \cup T_2^{-1} \cup \ldots \cup T_n^{-1}$ and $\alpha_1 = \alpha_1^1 \cup \alpha_2^1 \cup \ldots \cup \alpha_n^1$ and $\beta_1 = \beta_1^1 \cup \beta_2^1 \cup \ldots \cup \beta_n^1$.

**THEOREM 2.8:** *Let $T = T_1 \cup T_2 \cup \ldots \cup T_n$ be a n-linear transformation from $V = V_1 \cup V_2 \cup \ldots \cup V_n$ and $W = W_1 \cup W_2 \cup \ldots \cup W_n$ where dim V = $(n_1, n_2, \ldots, n_n)$ and dim W = $\left(n_{i_1}, n_{i_2}, \ldots, n_{i_n}\right)$ where $i_1, i_2, \ldots, i_n$ is a permutation of (1, 2, …, n.). Then $T_s$ is non singular if and only if $T_s$ carries each n-linearly independent n-subset of V into a n-linearly independent n-subset of W.*

*Proof:* Suppose first we assume $T_s$ is non singular. Let $S = S_1 \cup S_2 \cup \ldots \cup S_n$ be a n-linearly independent n-subset of $V = V_1 \cup V_2 \cup \ldots \cup V_n$ i.e. $S_i \subseteq V_i$ is a linearly independent subset of $V_i$, i = 1, 2, …, n. Let

$S = \{ \left(\alpha_1^1, \alpha_2^1, \ldots, \alpha_{k_1}^1\right) \cup \left(\alpha_1^2, \alpha_2^2, \ldots, \alpha_{k_2}^2\right) \cup \ldots \cup \left(\alpha_1^n, \alpha_2^n, \ldots, \alpha_{k_n}^n\right) \}$

$= S_1 \cup S_2 \cup \ldots \cup S_n \subset V_1 \cup V_2 \cup \ldots \cup V_n$. Given $T_s = T_1 \cup T_2 \cup \ldots \cup T_n$. Here $T_i: V_i \to W_j$, $\left(T_i\alpha_1^i, T_i\alpha_2^i, \ldots, T_i\alpha_{k_i}^i\right)$ are linearly independent for each i for if

$$C_1^i(T_i\alpha_1^i) + \ldots + C_{k_i}^i(T_i\alpha_{k_i}^i) = 0,$$

then $T_i\ (C_1^i\alpha_1^i + \ldots + C_{k_i}^i\alpha_{k_i}^i) = 0$ and since $T_i$ is non singular $(C_1^i\alpha_1^i + \ldots + C_{k_i}^i\alpha_{k_i}^i) = 0$ from which it follows each $C_j^i = 0$, j = 1, 2, …, $k_i$, because $S_i$ is an independent set. This is true of each i, i.e. $S = S_1 \cup S_2 \cup \ldots \cup S_n$ is an independent n-set. This shows the image of S under $T_s$ is independent. Suppose $T_s$ carries independent n-sets onto independent n sets. Let $\alpha = \alpha_1 \cup \alpha_2 \cup \ldots \cup \alpha_n$ be a non zero n vector of V.

Then if $S = S_1 \cup S_2 \cup \ldots \cup S_n = \alpha_1 \cup \alpha_2 \cup \ldots \cup \alpha_n$ with $S_i = \{\alpha_i\}$; i = 1, 2, …, n; is independent. The image n-set of S is the n-row vector $T_1\alpha_1 \cup T_2\alpha_2 \cup \ldots \cup T_n\alpha_n$ and this set is independent. Hence $T_s(\alpha) = T_1\alpha_1 \cup T_2\alpha_2 \cup \ldots \cup T_n\alpha_n \neq 0$ because the set consisting of the zero n-vector alone is dependent. Thus null space of $T_s$ is $0 \cup 0 \cup \ldots \cup 0$.

The following concept of non singular n-linear transformation is little different.

**DEFINITION 2.16:** *Let V and W be two same n-dimension spaces over F i.e. dim V = $(n_1, n_2, \ldots, n_n)$ and dim W = $\left(n_{i_1}, n_{i_2}, \ldots, n_{i_n}\right)$ where $(i_1, i_2, \ldots, i_n)$ is a permutation of ((1, 2, 3, …, n). If $T = T_1 \cup T_2 \cup \ldots \cup T_n$ is a special n-linear transformation of V into W i.e. if $T_i: V_i \to W_j$ then dim $V_i$ = dim $W_j = n_i$ for every i. Then T is n-non singular if each $T_i$ is non singular.*

In view of this the reader is expected to prove the following theorem.

**THEOREM 2.9:** *Let V and W be two n-vector spaces of same dimension defined over the same field F. T is special n-linear transformation from V into W. Then T is n-invertible if and only if T is non-singular.*

*Note:* We say the special n-linear transformation is n-invertible if and only if each $T_i = T_1 \cup T_2 \cup \ldots \cup T_n$ is invertible for i = 1, 2, …, n.

Now we proceed on to define the n-representation of n-transformations by n-matrices. ($n \geq 2$).

Let V be a n-vector space of n dimension $(n_1, n_2, \ldots, n_n)$ and W be a m-vector space of m-dimension $(m_1, m_2, \ldots, m_n)$ defined over the same field F. Let

$B = \{(\alpha_1^1, \alpha_2^1, \ldots, \alpha_{n_1}^1) \cup (\alpha_1^2, \alpha_2^2, \ldots, \alpha_{n_2}^2) \cup \ldots \cup (\alpha_1^n, \alpha_2^n, \ldots, \alpha_{n_n}^n)\}$

be a n-ordered n-basis of V. We say the n-basis is an n-ordered n-basis if each of the basis $(\alpha_1^i, \alpha_2^i, \ldots, \alpha_{n_i}^i)$ of $V_i$ is an ordered basis for i = 1, 2, …, n and

$$B^1 = \{(\beta_1^1, \beta_2^1, \ldots, \beta_{m_1}^1) \cup (\beta_1^2, \beta_2^2, \ldots, \beta_{m_2}^2) \cup \ldots \cup (\beta_1^m, \beta_2^m, \ldots, \beta_{m_m}^m)\}$$

be a m-ordered m basis of W. If T is any n-linear transformation from V into W i.e. $T = T_1 \cup T_2 \cup \ldots \cup T_n$ then each $T_i: V_i \to W_k$ is determined by its action on the vector $\alpha_j^i$; $1 \leq k \leq m$ ; true for each i = 1, 2, …, n and $i \leq j \leq n_i$. Each of the $n_i$ vector $T_i\alpha_j$ is uniquely expressible as a linear combination

$$T_i \alpha_j^i = \sum_{i=1}^{m_k} A_{ij} \beta_i^k \qquad (1)$$

$1 \leq k \leq m$ and $\beta_i^k \in W_k$, the scalars $A_{1j}, A_{2j}, \ldots, A_{m_k j}$ being coordinates of $T_i\alpha_j^i$ in the m-ordered m-basis $B^1$. Accordingly the transformation $T_i$ is determined by the $m_k n_i$ scalars; $A_{ij}$ via equation (1). The $m_k \times n_i$ matrix $A_i^k$ defined by $A_{i_j}$ is called the submatrix relative to the n-linear transformation $T = T_1 \cup T_2 \cup \ldots \cup T_i \cup \ldots \cup T_n$ of the pair of ordered basis

$$\{(\alpha_1^i, \alpha_2^i, ..., \alpha_{n_i}^i)\} \text{ and } \{(\beta_1^k, \beta_2^k, ..., \beta_{m_k}^k)\}$$

of $V_i$ and $W_k$ respectively. This is true for each i and k, $1 \leq i \leq n$ and $1 \leq k \leq m$ i.e. the m-matrix of T is given by

$$A_{n_1}^{m_{i_1}} \cup A_{n_2}^{m_{i_2}} \cup \ldots \cup A_{n_n}^{m_{i_n}}$$

= A (m, n) here $\left(m_{i_1}, m_{i_2}, ..., m_{i_n}\right) \subset (m_1, m_2, \ldots, m_{m_n})$. Clearly A is only a n-linear transformation map $V_i \to W_j$ and no two $V_i$'s are mapped onto same $W_j$, $1 \leq i \leq n$ and $1 \leq j \leq m$. Thus if $\alpha^i$ is a vector in $V_i$ then $\alpha^i = x_1^i a_1^i + ... + x_{n_i}^i a_{n_i}^i$ is a vector in $V_i$ then

$$T_i\alpha^i = T_i \left(\sum_{j=1}^{n_i} x_j^i \alpha_j^i\right)$$

$$= \left(\sum_{j=1}^{n_i} x_j^i\right).(T\alpha_j^i)$$

$$= \sum_{j=1}^{n_i} x_j^i \sum_{i=1}^{m_k} A_{ij}\beta_j^k$$

$$= \sum_{j=1}^{m_k}\left(\sum_{i=1}^{n_i} A_{ij} x_j^i\right)\beta_i^k .$$

This is true for each i, i = 1, 2, ..., n. If $X = X^1 \cup X^2 \cup \ldots \cup X^n$ is the coordinate n-matrix of α in the n-basis B then the computation above shows that

$AX = (A_{n_1}^{m_{i_1}} \cup A_{n_2}^{m_{i_2}} \cup \ldots \cup A_{n_n}^{m_{i_n}})(X^1 \cup X^2 \cup \ldots \cup X^n)$ is the coordinate n-matrix of the n-vector Tα in the ordered basis $B^1$ because the scalars

$$\sum_{j=1}^{n_1} A_{ij}^{m_{i_1}} x_j^1 \cup \sum_{j=1}^{n_2} A_{ij}^{m_{i_2}} x_j^2 \cup \ldots \cup \sum_{j=1}^{n_n} A_{ij}^{m_{i_n}} x_j^n$$

is the entry of the $i^{th}$ n-row of the n-column matrix AX. Let us observe that A is given by the $m_i \times n_j$, n-matrices over the field F, then

$$T_1\left(\sum_{j=1}^{n_1} x_j^1\alpha_j^1\right) \cup T_2\left(\sum_{j=1}^{n_2} x_j^2\alpha_j^2\right) \cup \ldots \cup T_n\left(\sum_{j=1}^{n_n} x_j^n\alpha_j^n\right)$$

$$= \sum_{i=1}^{m_{i_1}} \left( \sum_{j=1}^{n_1} A_{ij}^{m_{i_1}} x_j^1 \right) \beta_i^{i_1} \cup \sum_{i=1}^{m_{i_2}} \left( \sum_{j=1}^{n_2} A_{ij}^{m_{i_2}} x_j^2 \right) \beta_i^{i_2} \cup \ldots \cup$$

$$\sum_{i=1}^{m_{i_n}} \left( \sum_{j=1}^{n_n} A_{ij}^{m_{i_n}} x_j^n \right) \beta_i^{i_n} ,$$

where $(i_1, i_2, \ldots, i_n) \subset \{1, 2, \ldots, m\}$ taken in some order, defines a n-linear transformation T from V into W, the n-matrix of A relative to the n-basis B and m-basis $B^1$ which is stated by the following theorem.

**THEOREM 2.10:** *Let $V = V_1 \cup V_2 \cup \ldots \cup V_n$ be a finite n-dimensional i.e., $(n_1, n_2, \ldots, n_n)$ n-vector space over the field F and $W = W_1 \cup W_2 \cup \ldots \cup W_m$, an m-dimensional $(m_1, m_2, \ldots, m_n)$ vector space over the same field F, $(m > n)$. For each n-linear transformation T from V into W there is a n-mixed rectangular matrices A of orders $(m_1 \times n_1, m_2 \times n_n, \ldots, m_n \times n_n)$ with entries in F such that $[T\alpha]_{B^1} = A[\alpha]_B$ for every $\alpha \in V$. $T \to A$ is a one to one correspondence between the set of all n-linear transformations from V into W and the set of all $m_i \times n_i$, mixed rectangular n-matrices, $i = 1, 2, \ldots, n$ over the field F. The matrix $A = A_{n_1}^{m_{i_1}} \cup A_{n_2}^{m_{i_2}} \cup \ldots \cup A_{n_n}^{m_{i_n}}$ is the associated n-matrix with T; the n-linear transformation of V into W relative to the basis B and $B^1$.*

Several interesting results true for the usual vector spaces can be derived in case of n-vector spaces $n \geq 2$ with appropriate modifications.

Now we give the definition of n-inner product on a n-vector space V.

**DEFINITION 2.17:** *Let F be a field of reals or complex numbers and $V = V_1 \cup V_2 \cup \ldots \cup V_n$ a n-vector space over F. An n-inner product on V is a n-function which assigns to each ordered pair of n-vectors $\alpha = \alpha_1 \cup \alpha_2 \cup \ldots \cup \alpha_n$ and $\beta = \beta_1 \cup \beta_2 \cup \ldots \cup \beta_n$ in the n-vector space V a scalar n-tuple from F. $\langle \alpha \mid \beta \rangle = \langle \alpha_1 \cup \alpha_2 \cup \ldots \cup \alpha_n \mid \beta_1 \cup \beta_2 \cup \ldots \cup \beta_n \rangle = (\langle \alpha_1 \mid \beta_1 \rangle, \langle \alpha_2 \mid \beta_2 \rangle, \ldots, \langle \alpha_n \mid \beta_n \rangle)$,* *where* $\langle \alpha_i | \beta_i \rangle$ *is a inner product on* $V_i$ *as* $\alpha_i, \beta_i \in V_i$*, this is true for each i, i = 1, 2, …, n; satisfying the following conditions:*

a. $\langle \alpha + \beta \mid \gamma \rangle = \langle \alpha \mid \gamma \rangle + \langle \beta \mid \gamma \rangle$ *(where* $\alpha = \alpha_1 \cup \alpha_2 \cup \ldots \cup \alpha_n$, $\beta = \beta_1 \cup \beta_2 \cup \ldots \cup \beta_n$ *and* $\gamma = \gamma_1 \cup \gamma_2 \cup \ldots \cup \gamma_n$ *where* $\alpha_i, \beta_i, \gamma_i \in V_i$ *for each i = 1, 2, …, n.)* $= (\langle \alpha_1 \mid \gamma_1 \rangle, \langle \alpha_2 \mid \gamma_2 \rangle, \ldots, \langle \alpha_n \mid \gamma_n \rangle) + (\langle \beta_1 \mid \gamma_1 \rangle \langle \beta_2 \mid \gamma_2 \rangle, \ldots, \langle \beta_n \mid \gamma_n \rangle) = (\langle \alpha_1 \mid \gamma_1 \rangle + \langle \beta_1 \mid \gamma_1 \rangle, \langle \alpha_2 \mid \gamma_2 \rangle + \langle \beta_2 \mid \gamma_2 \rangle, \ldots, \langle \alpha_n \mid \gamma_n \rangle + \langle \beta_n \mid \gamma_n \rangle)$
b. $\langle C\alpha \mid \beta \rangle = C\langle \alpha \mid \beta \rangle = (C_1\langle \alpha_1 \mid \beta_1 \rangle, C_2\langle \alpha_2 \mid \beta_2 \rangle, \ldots, C_n\langle \alpha_n \mid \beta_n \rangle)$
c. $\langle \beta \mid \alpha \rangle = \overline{\langle \alpha \mid \beta \rangle}$*, the bar denoting the complex conjugation.*
d. $\langle \alpha \mid \alpha \rangle > (0, 0, \ldots, 0)$ *if* $\alpha \neq 0$ *i.e.,* $(\langle \alpha_1 \mid \alpha_1 \rangle, \langle \alpha_2 \mid \alpha_2 \rangle, \ldots \langle \alpha_n \mid \alpha_n \rangle) > (0, 0, \ldots, 0)$ *each* $\alpha_i \neq 0$ *in* $\alpha = \alpha_1 \cup \alpha_2 \cup \ldots \cup \alpha_n$, *i = 1, 2, …, n.*

*On* $F = F^{n_1} \cup F^{n_2} \cup \ldots \cup F^{n_n}$ *there is a n-inner product which we call the n-standard inner product. It is defined on*

$$\alpha = (x_1^1, x_2^1, \ldots, x_{n_1}^1) \cup (x_1^2, x_2^2, \ldots, x_{n_2}^2) \cup \ldots \cup (x_1^n, x_2^n, \ldots, x_{n_n}^n)$$

*and*

$$\beta = (y_1^1, y_2^1, \ldots, y_{n_1}^1) \cup (y_1^2, y_2^2, \ldots, y_{n_2}^2) \cup \ldots \cup (y_1^n, y_2^n, \ldots, y_{n_n}^n) \in P$$

*by*

$$\langle \alpha \mid \beta \rangle = (\sum_{j=1}^{n_1} x_j^1 y_j^1, \sum_{j=1}^{n_2} x_j^2 y_j^2, \ldots, \sum_{j=1}^{n_n} x_j^n y_j^n)$$

*if F is a real field. If F is the field of complex numbers then*

$$\langle \alpha \mid \beta \rangle = (\sum_{j=1}^{n_1} x_j^1 \overline{y}_j^1, \sum_{j=1}^{n_2} x_j^2 \overline{y}_j^2, \ldots, \sum_{j=1}^{n_n} x_j^n \overline{y}_j^n).$$

The reader is expected to work out the properties related with n-inner products on the n-vector spaces over the field F.

Now we proceed on to define n-orthogonal sets.

**DEFINITION 2.18:** *Let* $V = V_1 \cup V_2 \cup \ldots \cup V_n$ *be a n-vector space over the field F. We say V is a n-inner product space if on V is defined an n-inner product. Let* $\alpha = (\alpha_1 \cup \alpha_2 \cup \ldots \cup \alpha_n)$ *and*

*$\beta = (\beta_1 \cup \beta_2 \cup \ldots \cup \beta_n) \in V$ with $\alpha_i$, $\beta_i \in V_i$, $i = 1, 2, \ldots, n$. We say $\alpha$ is n-orthogonal to $\beta$ if $\langle\alpha|\beta\rangle = (0, 0, \ldots, 0) = (\langle\alpha_1|\beta_1\rangle, \ldots, \langle\alpha_n|\beta_n\rangle)$ i.e. if each $\alpha_i$ is orthogonal to $\beta_i \in V_i$ i.e. $\langle\alpha_i|\beta_i\rangle = 0$ for $i = 1, 2, \ldots, n$. This equivalently implies $\beta$ is n-orthogonal to $\alpha$. Hence we simply say $\alpha$ and $\beta$ are orthogonal.*

*If $S = S_1 \cup S_2 \cup \ldots \cup S_n \subseteq V_1 \cup V_2 \cup \ldots \cup V_n = V$ be a n-set of n-vectors in V. S is called an n-orthogonal set provided all pairs of distinct n-vectors in S are orthogonal. An n-orthogonal set is called an n-orthonormal set if $||\alpha|| = (1, 1, \ldots, 1)$ for every $\alpha$ in $S = S_1 \cup S_2 \cup \ldots \cup S_n$.*

*We denote $\langle\alpha \mid \beta\rangle$ also by $(\alpha|\beta)$.*

**THEOREM 2.11:** *Let $V = V_1 \cup V_2 \cup \ldots \cup V_n$ be a n-vector space which is a n-inner product space defined over the field. Let $S = S_1 \cup S_2 \cup \ldots \cup S_n$ be an n-orthogonal set in V. The set of non zero vectors in S are n-linearly independent.*

*Proof:* Let $V = V_1 \cup V_2 \cup \ldots \cup V_n$ be a n-vector space over F. Let $S = S_1 \cup S_2 \cup \ldots \cup S_n \subset V = V_1 \cup V_2 \cup \ldots \cup V_n$ be a orthogonal n-set of V. To show the elements in the n-sets are n-orthogonal. Let $\alpha_1^i, \alpha_2^i, \ldots, \alpha_{m_i}^i \in S_i$ for $i = 1, 2, \ldots, n$. i.e. $\alpha_1^1, \alpha_2^1, \ldots, \alpha_{m_1}^1 \in S_1$, $\alpha_1^2, \alpha_2^2, \ldots, \alpha_{m_2}^2 \in S_2$ and so on. $\alpha_1^n, \alpha_2^n, \ldots, \alpha_{m_n}^n \in S_n$. Let $\alpha_1^i, \alpha_2^i, \ldots, \alpha_{m_i}^i$ be the distinct set of n-vectors in $S_i$ and that $\beta^i = c_1^i\alpha_1^i + c_2^i\alpha_2^i + \ldots + c_{m_i}^i\alpha_{m_i}^i$. Then

$$(\beta^i \mid \alpha_k^i) = \left(\sum_{j=1}^{m_i} c_j^i \alpha_j^i \mid \alpha_k^i\right)$$

$$= \sum_{j=1} c_j^i (\alpha_j^i \mid \alpha_k^i)$$

$$= c_k^i \; (\alpha_k^i \mid \alpha_k^i).$$

Since $(\alpha_k^i \mid \alpha_k^i) \neq 0$, $c_k^i \neq 0$.

Thus when $\beta^i = 0$ then each $c_k^i = 0$ for each i. So each $S_i$ is an independent set. Hence $S = S_1 \cup S_2 \cup \ldots \cup S_n$ is a n-independent set.

Several interesting results including Gram-Schmidt n-orthogonalization process can be derived.

We now proceed onto define the notion of n-best approximation to the n-vector β relative to a n-sub-vector space W.

**DEFINITION 2.19:** *Let $V = V_1 \cup V_2 \cup \ldots \cup V_n$ be a n-inner product n-vector space over the field F. $W = W_1 \cup W_2 \cup \ldots \cup W_n$ be a n-subspace of V. Let $\beta = \beta_1^1 \cup \beta_2^1 \cup \ldots \cup \beta_n^n \in V$ the n-best approximation to β by n-vectors in W is a n-vector $\alpha = \alpha_1^1 \cup \alpha_2^1 \cup \ldots \cup \alpha_n^n$ in W such that $||\beta - \alpha|| \leq ||\beta - \gamma||$ for every n-vector γ in W i.e. $||\beta^i{}_i - \alpha^i{}_i|| \leq || \beta^i{}_i - \gamma^i{}_i||$ for every $\gamma^i{}_i \in W_i$ and this is true for each i; i = 1, 2, …, n. We know if $\beta = \beta_1^1 \cup \beta_2^1 \cup \ldots \cup \beta_n^n$ and if β is a n-linear combination of an n-orthogonal sequence of non zero-vectors $\alpha_1, \alpha_2, \ldots, \alpha_m$ where each $\alpha_i = \alpha^i{}_1 \cup \alpha^i{}_2 \cup \ldots \cup \alpha^i{}_n$, i = 1,2,…, m, then*

$$\beta = \left( \sum_{k=1}^{m_1} \frac{\left(\beta_1^1 \middle| \alpha_k^1\right)\alpha_k^1}{\left\|\alpha_k^1\right\|^2} \cup \sum_{k=1}^{m_2} \frac{\left(\beta_2^2 \middle| \alpha_k^2\right)\alpha_k^2}{\left\|\alpha_k^2\right\|^2} \cup \ldots \cup \sum_{k=1}^{m_n} \frac{\left(\beta_n^n \middle| \alpha_k^n\right)\alpha_k^n}{\left\|\alpha_k^n\right\|^2} \right).$$

The following theorem is left as an exercise for the reader to prove.

**THEOREM 2.12:** *Let $W = W_1 \cup W_2 \cup \ldots \cup W_n$ be a n-subspace of an n-inner product space $V = V_1 \cup V_2 \cup \ldots \cup V_n$ and $\beta = \beta_1^1 \cup \beta_2^1 \cup \ldots \cup \beta_n^n$ be a n-vector in V*

1. *The n-vector $\alpha = \alpha_1^1 \cup \alpha_2^1 \cup \ldots \cup \alpha_n^n$ in W is a n-best approximation to β by the n-vector in W if and only if β-α is n-orthogonal to every n-vector in W.*
2. *If a n-best approximation to β by n-vectors in W exists, it is unique.*

3. *If W is n-finite n-dimensional and* $\{\alpha_1^1 \cup \alpha_2^1 \cup \ldots \cup \alpha_n^n\}$ *is any n-orthonormal n-basis for W then the n-vector*

$$\alpha = \sum_k \frac{(\beta_1^1 \mid \alpha_k^1)\alpha_k^1}{\|\alpha_k^1\|^2} \cup \sum_k \frac{(\beta_2^2 \mid \alpha_k^2)\alpha_k^2}{\|\alpha_k^2\|^2} \cup \ldots \cup \sum_k \frac{(\beta_n^n \mid \alpha_k^n)\alpha_k^n}{\|\alpha_k^n\|^2}$$

*is the unique n-best approximation to β by n-vector in W.*

Now we proceed on to define the notion of n-orthogonal complement.

**DEFINITION 2.20:** *Let V be a n-inner product n-space and S any n-set of n vectors in V. The n-orthogonal complement of S is the n-set* $S^{\perp}$ *of all n-vectors in V which are n-orthogonal to every n-vector in S; where* $S = S_1 \cup S_2 \cup \ldots \cup S_n \subseteq V = V_1 \cup V_2 \cup \ldots \cup V_n$ *and* $S^{\perp} = S_1^{\perp} \cup S_2^{\perp} \cup \ldots \cup S_n^{\perp} \subseteq V$*. i.e. each* $S_i^{\perp}$ *is the orthogonal complement of* $S_i$ *for every i, i = 1, 2, ..., n. We call α to be the n-orthogonal projection of β on W. If every n-vector in V has an n-orthogonal projection on W, the n-mapping that assigns to each n-vector in V its n-orthogonal projection on W is called the n-orthogonal projection of V on W.*

The reader is expected to prove the following theorems.

**THEOREM 2.13:** *Let* $V = V_1 \cup V_2 \cup \ldots \cup V_n$ *be a n-inner product n-vector space defined over the field F. W a finite dimensional n-subspace of V and E the n-orthogonal projection of V on W, Then the n-mapping* $\beta \rightarrow (\beta - E\beta)$ *is the n-orthogonal projection of V on* $W^{\perp}$*.*

**THEOREM 2.14:** *Let* $W = W_1 \cup W_2 \cup \ldots \cup W_n \subset V$ *be a finite dimensional n-subspace of the n inner product space* $V = V_1 \cup V_2 \cup \ldots \cup V_n$ *and let* $E = E_1 \cup E_2 \cup \ldots \cup E_n$ *be the n-orthogonal projection of V on W. Then E is an n-idempotent n-linear transformation of V onto W and* $W^{\perp}$ *is the n-null space of E and*

*$V = W \oplus W^{\perp}$ i.e. if $V = V_1 \cup V_2 \cup \ldots \cup V_n$ and $W = W_1 \cup W_2 \cup \ldots \cup W_n$ and $W^{\perp} = W_1^{\perp} \cup W_2^{\perp} \cup \ldots \cup W_n^{\perp}$ then $V = (W \oplus W^{\perp}) = (W_1 \oplus W_1^{\perp}) \cup (W_2 \oplus W_2^{\perp}) \cup \ldots \cup (W_n \oplus W_n^{\perp})$.*

**THEOREM 2.15:** *Under the conditions of the above theorems, I-E is the n-orthogonal projection of V on $W^{\perp}$. It is an n-idempotent linear n-transformation of V onto $W^{\perp}$, with null space W.*

**THEOREM 2.16:** *Let $\{\alpha_1^1 \cup \alpha_2^1 \cup \ldots \cup \alpha_n^n\}$ be an orthogonal n-set of non-zero vectors in an n-inner product space V. If β is any vector in V then $\sum_k |(\beta|\, \alpha_k^k)|^2 \,/\, \|\alpha_k^k\|^2 \leq \|\beta\|^2$ where $\beta = \beta_1^1 \cup \beta_2^1 \cup \ldots \cup \beta_n^n \in V$.*

It is pertinent to mention here that the notion of linear functional dual space or adjoints cannot be extended in an analogous way in case of n-vector spaces of type I.

Now we proceed on to define the notion of n-unitary operators on n-inner product n vector spaces V over the field F.

**DEFINITION 2.21:** *Let V and W be n-inner product n-vector space and m vector space over the same field F respectively. Let T be a n-linear transformation from V into W. We say that T preserves n inner products if $(T\alpha \mid T\beta) = (\alpha \mid \beta)$ for all $\alpha,\beta \in V$ i.e. if $V = V_1 \cup V_2 \cup \ldots \cup V_n$ and $W = W_1 \cup W_2 \cup \ldots \cup W_n$ and $T = T_1 \cup T_2 \cup \ldots \cup T_n$ with $\alpha = \alpha_1 \cup \alpha_2 \cup \ldots \cup \alpha_n$ and $\beta = \beta_1 \cup \beta_2 \cup \ldots \cup \beta_n \in V$. $T_i: V_i \to W_j$. with no two $V_i$ mapped on to the same $W_j$, then $T_i\alpha_i, T_i\beta_i \in W_j$ and $(T_i\alpha_i \mid T_i\beta_i) = (\alpha_i \mid \beta_i)$ for every i, i = 1, 2, …, n. An n-isomorphism of V into W is a n-vector space isomorphism T of V onto W which also preserves n-inner products.*

**THEOREM 2.17:** *Let V and W be n-finite dimensional n-inner vector spaces of same n-dimension i.e. dim $V = (n_1, n_2, \ldots, n_n)$ and dim $W = \left(n_{i_1}, n_{i_2}, \ldots, n_{i_n}\right)$ where $(i_1, i_2, \ldots, i_n)$ is a* *permutation of (1, 2, …, n) defined over the same field T. If $T = T_1 \cup T_2 \cup \ldots \cup T_n$ is a n-linear transformation from V into W the following are equivalent*

1. *T preserves inner products i.e., each $T_i$ in T preserves inner product i.e. $T_i$: $V_i \to W_j$; $1 \leq i, j \leq n$.*
2. *T is an n-inner product n- isomorphism*
3. *T carries every n-orthonormal n-basis for V onto an n-orthogonal n-basis for W.*
4. *T carries some n-orthogonal n-basis for V onto an n-orthonormal basis for W i.e. $T_i$ carries some orthogonal basis of $V_i$ into an orthogonal basis for $W_j$.*

The reader is expected to prove the following theorems.

**THEOREM 2.18:** *Let V and W be n-dimensional finite inner product n-spaces over the same field F. Then $V = V_1 \cup V_2 \cup \ldots \cup V_n$ is n-isomorphic with $W = W_1 \cup W_2 \cup \ldots \cup W_n$ i.e. each $T_i$: $V_i \to W_j$ is an isomorphism for i = 1, 2, …, n if V and W are of same n-dimension.*

**THEOREM 2.19:** *Let V and W be two n-inner product spaces over the same field F. Let $T = T_1 \cup T_2 \cup \ldots \cup T_n$ be a n-linear transformation from V into W. Then T preserves n-inner product if and only if $||T\alpha|| = ||\alpha||$ i.e. $||(T_1 \cup T_2 \cup \ldots \cup T_n)(\alpha_1^1 \cup \alpha_2^1 \cup \ldots \cup \alpha_n^n)|| = ||T_1(\alpha_1^1) \cup T_2(\alpha_2^2) \cup \ldots \cup T_n(\alpha_n^n)|| = (||\alpha_1^1||, ||\alpha_2^2||, \ldots, ||\alpha_n^n||)$ for every $\alpha \in V$ i.e. for every $\alpha_i \in V_i$, i = 1, 2, …, n.*

We define the notion of n unitary operator of a n-vector space V over the field F.

**DEFINITION 2.22:** *A n-unitary operator on an n-inner product space V is a n-isomorphism of V onto itself.*

**DEFINITION 2.23:** *If T is a n-linear operator on an n-inner product space $V = V_1 \cup V_2 \cup \ldots \cup V_n$, then we say $T = T_1 \cup T_2 \cup \ldots \cup T_n$ has an n-adjoint on V if there exists a n-linear* *operator* $T^* = T_1^* \cup T_2^* \cup \ldots \cup T_n^*$ *on V such that* $(T\alpha \mid \beta) = (\alpha \mid T^*\beta)$ *for all* $\alpha = \alpha_1^1 \cup \alpha_2^1 \cup \ldots \cup \alpha_n^n$ , $\beta = \beta_1^1 \cup \beta_2^1 \cup \ldots \cup \beta_n^n$ *in* $V = V_1 \cup V_2 \cup \ldots \cup V_n$ *i.e.* $T_i(\alpha_i^i \mid \beta_i^i) = (\alpha_i^i \mid T^* \beta_i^i)$ *for each* $i = 1, 2, \ldots, n$.

It is easily verified as in case of adjoints the n-adjoints of T not only depends on T but also on the n-inner product on V.

Interesting results in this direction can be derived for any reader. The following theorems are also left as an exercise for the reader.

**THEOREM 2.20:** *Let* $V = V_1 \cup V_2 \cup \ldots \cup V_n$ *be a finite n-dimensional n-inner product n-space defined over the field F. If T and U are n-linear operators on V and c is a scalar, then*

1. $(T + U)^* = T^* + U^*$ *i.e. if* $T = T_1 \cup T_2 \cup \ldots \cup T_n$ *and* $U = U_1 \cup U_2 \cup \ldots \cup U_n$ *then in* $(T + U)^*$ *we have for each* $i$, $(T_i + U_i)^* = T_i^* + U_i^*$, $i = 1, 2, \ldots, n$.

2. $(cT)^* = cT^*$

3. $(TU)^* = T^*U^*$, *here also* $(T_iU_i)^* = U_i^* T_i^*$ *for* $i = 1, 2, \ldots, n$. *i.e.* $(TU)^* = (T_1U_1)^* \cup (T_2U_2)^* \cup \ldots \cup (T_nU_n)^* = U_1^* T_1^* \cup U_2^* T_2^* \cup \ldots \cup U_n^* T_n^*$

4. $(T^*)^* = T$ *since* $(T_i^*)^* = T_i$, *for each* $i = 1, 2, \ldots, n$.

**THEOREM 2.21:** *Let U be a n-linear operator on an n-inner product space V, defined over the field F. Then U is n-unitary if and only the n-adjoint, U* of U exists and* $UU^* = U^*U = I$.

**THEOREM 2.22:** *Let* $V = V_1 \cup V_2 \cup \ldots \cup V_n$ *be a n-vector space of a n-inner vector space of finite dimension and U be a n-linear operator on V. Then U is n-unitary if and only if the n-matrix related with U in some ordered n-orthonormal n-basis is also a n-unitary matrix i.e. if* $A = A_1 \cup A_2 \cup \ldots \cup A_n$ *is the n-matrix* *each $A_i$ in A is unitary i.e. $A_i * A_i = I$ for each i, i.e. $A * A = A_1{*}A_1 \cup A_2{*}A_2 \cup \ldots \cup A_n{*}A_n = I_1 \cup \ldots \cup I_n$.*

Several interesting results can be obtained using appropriate and analogous proper modifications.

Now we proceed on to define the notion of n-normal operator or normal n-operators on a n-vector space V. The principle objective for doing this is we can obtain some interesting properties about the n-orthonormal n-basis of $V = V_1 \cup V_2 \cup \ldots \cup V_n$.

*Let the n-orthonormal n-basis of V be denoted by $B = \{(\alpha_1^1, \alpha_2^1, \ldots, \alpha_{n_1}^1) \cup (\alpha_1^2, \alpha_2^2, \ldots, \alpha_{n_2}^2) \cup \ldots \cup (\alpha_1^n, \alpha_2^n, \ldots, \alpha_{n_n}^n)\}$ where each $(\alpha_1^i, \alpha_2^i, \ldots, \alpha_{n_i}^i)$ is a orthogonal basis of $V_i$ for i = 1, 2, …, n. Let $T = T_1 \cup T_2 \cup \ldots \cup T_n$ the n-linear operator on V be defined by $T_i \alpha_j^i = c_j^i \alpha_j^i$ for j = 1, 2, …, $n_i$ and for each $T_i$, i = 1, 2, …, n. This simply implies that the n-matrix of T (consequently each matrix of $T_i$ in the ordered basis $(\alpha_1^i, \alpha_2^i, \ldots, \alpha_{n_i}^i)$ is a diagonal matrix with the diagonal entries $(c_1^i, c_2^i, \ldots, c_{n_i}^i)$ is a n-diagonal n matrix given by*

$$D = \begin{bmatrix} c_1^1 & & & 0 \\ & c_2^1 & & \\ & & \ddots & \\ 0 & & & c_{n_1}^1 \end{bmatrix} \cup \begin{bmatrix} c_1^2 & & & 0 \\ & c_2^2 & & \\ & & \ddots & \\ 0 & & & c_{n_2}^2 \end{bmatrix} \cup \ldots \cup \begin{bmatrix} c_1^n & & & 0 \\ & c_2^n & & \\ & & \ddots & \\ 0 & & & c_{n_n}^n \end{bmatrix}.$$

*The n-adjoint operator $T^* = T^*_1 \cup T^*_2 \cup \ldots \cup T^*_n$ of $T = T_1 \cup T_2 \cup \ldots \cup T_n$ is represented by the n-conjugate transpose n matrix i.e. once again a n-diagonal n matrix with diagonal entries $c_1^{-1}, c_2^{-1}, \ldots, c_{n_i}^{-i}$ ; i = 1, 2, ..., n. If V is a real n-vector space over the real field F then of course we have T = T**

**DEFINITION 2.24:** *If $V = V_1 \cup V_2 \cup \ldots \cup V_n$ be a n-dimensional n-inner product n-vector space and T a n-linear operator on V be say T is n-normal if it commutes with its n-adjoint T* of T i.e. $TT^* = T^*T$.*

Now in order to define some more properties we now proceed onto define the notion of n-characteristic values or n-eigen values of a n-vector space V and so on.

**DEFINITION 2.25:** *Let $V = V_1 \cup V_2 \cup \ldots \cup V_n$ be a n-vector space over the field F of type I. Let $T = T_1 \cup T_2 \cup \ldots \cup T_n$ be a n-linear operator on V. A n-characteristic value (or equivalently characteristic n-value) of T is a n-tuple of scalars $c_1^1 \cup c_2^2 \cup \ldots \cup c_n^n$ such that their exists a non zero n vector $\alpha = \alpha_1^1 \cup \alpha_2^2 \cup \ldots \cup \alpha_n^n$ in V with $T\alpha = c\alpha$. i.e. $(c_1^1 \cup c_2^2 \cup \ldots \cup c_n^n)$ $(\alpha_1^1 \cup \alpha_2^2 \cup \ldots \cup \alpha_n^n) = c_1^1\alpha_1^1 \cup c_2^2\alpha_2^2 \cup \ldots \cup c_n^n\alpha_n^n = T_1\alpha_1^1 \cup T_2\alpha_2^2 \cup \ldots \cup T_n\alpha_n^n$ ; If $c = c_1^1 \cup c_2^2 \cup \ldots \cup c_n^n$ is the n-characteristic value of T then*

a. *any $\alpha = \alpha_1^1 \cup \alpha_2^2 \cup \ldots \cup \alpha_n^n$ such that $T\alpha = c\alpha$ is called the n-characteristic n-vector of T associated with the n-characteristic value $c = c_1^1 \cup c_2^2 \cup \ldots \cup c_n^n$.*

b. *The collection of all $\alpha = \alpha_1^1 \cup \alpha_2^2 \cup \ldots \cup \alpha_n^n$ such that $T\alpha = c\alpha$ is called the n-characteristic space associated with c. n-characteristic values will also be known as n-eigen values or n-spectral values.*

*If T is any n-linear operator on the n-vector space V and c any n scalar the set of n-vector α in $V = V_1 \cup V_2 \cup \ldots \cup V_n$ such that*

*$T\alpha = c\alpha$ is a n-subspace of V. It is the n-null space of the n-linear transformation $(T - cI) = (T_1 \cup T_2 \cup \ldots \cup T_n) - (c_1^1 \cup c_2^2 \cup \ldots \cup c_n^n)\ (I_1 \cup I_2 \cup \ldots \cup I_n)) = (T_1 - c_1^1 I_1) \cup (T_2 - c_2^2 I_2) \cup \ldots \cup (T_n - c_n^n I_n)$. We call c the n-characteristic value of T and if this n-subspace is different from the zero subspace i.e. if $(T - cI) = (T_1 - c_1^1 I_1) \cup \ldots \cup (T_n - c_n^n I_n)$ fails to be one to one i.e. each $T_i - c_i^i I_i$ fails to be one to one. If V is a finite n-dimension n-vector space, $(T - cI)$ fails to be one to one. Only when the n determinant i.e. $\det(T - cI) = \det(T_1 - c_1^1 I_1) \cup \det(T_2 - c_2^2 I_2) \cup \ldots \cup \det(T_n - c_n^n I_n) \neq (0 \cup 0 \cup \ldots \cup 0)$ i.e. each $\det(T_i - c_i^i I_i) \neq 0$ for $i = 1, 2, \ldots, n$.*

This is made into the following nice theorem

**THEOREM 2.23:** *Let T be a n-linear operator on a finite n-dimensional n-vector space $V = V_1 \cup V_2 \cup \ldots \cup V_n$ and $c = c_1^1 \cup c_2^2 \cup \ldots \cup c_n^n$ be a n scalar then the following are equivalent*

a. *$c = c_1^1 \cup c_2^2 \cup \ldots \cup c_n^n$ is a n-characteristic value of $T = T_1 \cup T_2 \cup \ldots \cup T_n$ i.e. each $c^i_i$ is a characteristic value of $T_i$ ; $i = 1, 2, \ldots, n$.*
b. *The n-operator $(T - cI) = (T_1 - c_1^1 I_1) \cup \ldots \cup (T_n - c_n^n I_n)$ is non singular (i.e. non invertible) i.e. each $(T_i - c_i^i I_i)$ is non invertible i.e. non singular for each n-vector spaces, $i = 1, 2, \ldots, n$.*
c. *$\det(T - cI) = (0 \cup 0 \cup \ldots \cup 0)$ i.e. det $(T_i - c_i^i I_i) = 0$ for each $i = 1, 2, \ldots, n$.*

Now we give the analogous for n-matrix.

**DEFINITION 2.26:** *Let $A = A_1 \cup A_2 \cup \ldots \cup A_n$ be a n-square matrix where each matrix $A_i$ is $n_i \times n_i$ matrix $i = 1, 2, \ldots, n$; if $i \neq j$ then $n_i \neq n_j$, $1 \leq i, j \leq n$ over the field F, a n-characteristic value of A in F is a n scalar $C = C_1^1 \cup C_2^2 \cup \ldots \cup C_n^n$ ; $C_i^i \in F$, i =*

*1, 2, 3, ... , n such that the n-matrix (A – CI) =* $(A_1 - C_1^1 I_1) \cup (A_2 - C_2^2 I_2) \cup \ldots \cup (A_n - C_n^n I_n)$ *is singular, i.e. C is a n-characteristic value of A if and only if det (A – CI) = 0* $\cup$ *0* $\cup$ *...* $\cup$ *0 or equivalently det (CI – A) = 0* $\cup$ *0* $\cup$ *...* $\cup$ *0, i.e. if* $\det(A_i - C_i^i I_i) = 0$ *for each and every i, i = 1, 2, ..., n we form the matrix (xI – A) where* $x = x_1 \cup x_2 \cup \ldots \cup x_n$ *with polynomial entries and consider the n-polynomial* $f = \det(xI - A) = \det(x_1 I_1 - A_1) \cup \det(x_2 I_2 - A_2) \cup \ldots \cup \det(x_n I_n - A_n) = f_1 \cup f_2 \cup \ldots \cup f_n$ *in n variables* $x_1, x_2, \ldots, x_n$*. Clearly the n-characteristic values of A in F are just the n-tuple scalars C =* $C_1^1 \cup C_2^2 \cup \ldots \cup C_n^n$ *in F such that*

$$\begin{aligned} f(C) &= 0 \cup 0 \cup \ldots \cup 0 \\ &= f_1(C_1^1) \cup f_2(C_2^2) \cup \ldots \cup f_n(C_n^n). \end{aligned}$$

*For this reason f is called the n-characteristic n-polynomial of A.*

It is important to note that f is a n-monic polymonial which has degree exactly $(n_1, n_2, \ldots, n_n)$ is the n-degree of the n-monic polynomial $f = f_1 \cup f_2 \cup \ldots \cup f_n$.

We can prove the following simple lemma.

**LEMMA 2.2**: *Similar n-matrices have same n-characteristic polynomial.*

*Proof:* We just recall if A and B are the mixed square n-matrices of dimension $(n_1, n_2, \ldots, n_n)$ and $(n_1, n_2, \ldots, n_n)$ i.e. same dimension i.e. identity permutation of $(n_1, n_2, \ldots, n_n)$. We say A is similar to B or B is similar to A if their exists a invertible n matrix P of dimension $(n_1, n_2, \ldots, n_n)$ such that if $A = A_1 \cup A_2 \cup \ldots \cup A_n$, $B = B_1 \cup B_2 \cup \ldots \cup B_n$ and $P = P_1 \cup P_2 \cup \ldots \cup P_n$ then $B = P^{-1}AP$ i.e. $B = B_1 \cup B_2 \cup \ldots \cup B_n = P_1^{-1} A_1 P_1 \cup P_2^{-1} A_2 P_2 \cup \ldots \cup P_n^{-1} A_n P_n$; i.e. each $A_i$ is similar $B_i$ for i = 1, 2, …, n.

Suppose A and B are similar n-mixed square matrices of identical dimension i.e. order $A_i$ = order $B_i$ for i = 1, 2, …, n then $B = P^{-1}AP$ the

$$\begin{aligned}
\text{n-det}(xI - B) &= \text{n-det}(xI - P^{-1} A P) \\
&= \text{n-det}(P^{-1}(xI - A) P) \\
&= \text{n-det } P^{-1} \det(xI - A) .\det P \\
&= \det P_1^{-1} \det(x_1^1 I_1 - A_1) \det P_1 \cup \\
&\quad \det P_2^{-1} \det(x_2^2 I_2 - A_2) \det P_2 \cup \ldots \cup \\
&\quad \det P_n^{-1} \det(x_n^n I_n - A_n) \det P_n \\
&= \det(xI - B) \\
&= \det(x_1^1 I_1 - B_1) \cup \det(x_2^2 I_2 - B_2) \cup \ldots \cup \det(x_n^n I_n - B_n)
\end{aligned}$$

i.e. n-det (xI – B) = n-det (xI – A) .

**DEFINITION 2.27**: *Let $T = T_1 \cup T_2 \cup \ldots \cup T_n$ be a special n-linear operator on a finite dimension n-vector space $V = V_1 \cup V_2 \cup \ldots \cup V_n$. We say T is n diagonalizable if there is a n-basis for V each n-vector of which is a n-characteristic n-vector of T.*

The following two lemmas are left as an exercise to the reader.

**LEMMA 2.3**: *Suppose $T\alpha = C\alpha$ where $T = T_1 \cup T_2 \cup \ldots \cup T_n$, $\alpha = \alpha_1^1 \cup \alpha_2^2 \cup \ldots \cup \alpha_n^n$ and $C = C_1^1 \cup C_2^2 \cup \ldots \cup C_n^n$. If $F = F_1 \cup F_2 \cup \ldots \cup F_n$ is any n-polynomial then $f(T)\,\alpha = f(C)\,\alpha$ i.e.,*

$$f_1(T_1)\alpha_1^1 \cup f_2(T_2)\alpha_2^2 \cup \ldots \cup f_n(T_n)\alpha_n^n =$$
$$f_1(C_1^1)\alpha_1^1 \cup f_2(C_2^2)\alpha_2^2 \cup \ldots \cup f_n(C_n^n)\alpha_n^n .$$

**LEMMA 2.4**: *Let $T = T_1 \cup T_2 \cup \ldots \cup T_n$ be a n-linear operator on a finite $(n_1, n_2, \ldots, n_n)$ dimensional n-vector space $V = V_1 \cup V_2 \cup \ldots \cup V_n$. Let $\{(C_1^1, C_2^1, \ldots, C_{k_1}^1) \cup (C_1^2, C_2^2, \ldots, C_{k_2}^2) \cup \ldots \cup (C_1^n, C_2^n, \ldots, C_{k_n}^n)\}$ be distinct n-characteristic values of $T_1 \cup T_2 \cup \ldots \cup T_n$ and let $W_1^1, W_2^1, \ldots, W_{k_1}^1$ be the subspaces of $V_1$ associated with characteristic values $C_1^1, C_2^1, \ldots, C_{k_1}^1$ respectively, $W_1^2, W_2^2, \ldots, W_{k_2}^2$ be the subspaces of $V_2$ with associated* *characteristic values* $C_1^2, C_2^2, \ldots, C_{k_2}^2$ *respectively; and so on and let* $W_1^n, W_2^n, \ldots, W_{k_n}^n$ *be the subspaces of* $V_n$ *with associated characteristic values* $C_1^n, C_2^n, \ldots, C_{k_n}^n$ *and if*

$$W^1 = W_1^1 + W_2^1 + \ldots + W_{k_1}^1$$

*and if*

$$W^2 = W_1^2 + W_2^2 + \ldots + W_{k_2}^2, \ \ldots, \ W^n = W_1^n + W_2^n + \ldots + W_{k_n}^n$$

*and if* $W = W^1 \cup W^2 \cup \ldots \cup W^n$ *the n-dim W =* $(\dim W^1, \dim W^2, \ldots, \dim W^n)$ *with* $\dim W^j = \dim W_1^j + \dim W_2^j + \ldots + \dim W_{k_j}^j$ *for each j = 1, 2, …, n; and if* $B_i^i$ *is an ordered basis of* $W^i$*, i = 1, 2, …, k; then* $(B_1^1, B_2^2, \ldots, B_n^n)$ *is the n-ordered n-basis of* $W^1\ W^2 \ldots W^k$.

Using these lemmas the reader is expected to prove the following theorem.

**THEOREM 2.24:** *Let* $T = T_1 \cup T_2 \cup \ldots \cup T_n$ *be a n-linear operator of the finite n-dimensional n vector space* $V = V_1 \cup V_2 \cup \ldots \cup V_n$*. Let* $\{(C_1^1, C_2^1, \ldots, C_{k_1}^1), (C_1^2, C_2^2, \ldots, C_{k_2}^2) \cup \ldots \cup (C_1^n, C_2^n, \ldots, C_{k_n}^n)\}$ *be the distinct n-characteristic n-values of T and let* $(W^1 \cup W^2 \cup \ldots \cup W^n)$ *be the null n-subspace of* $(T - CI)$ *i.e.* $W^i$ *is a subspace of* $(T_i - C_i^i I)$ *for i = 1, 2, …, n. Then the following are equivalent*

1. *T is n-diagonalizable*
2. *The n-characteristic polynomial for T is*

$f = f_1 \cup f_2 \cup \ldots \cup f_n$ *where* $f_i = (x_i - C_1^i)^{d_1^i}(x_i - C_2^i)^{d_2^i} \ldots (x_i - C_{k_i}^i)^{d_{k_i}^i}$ *for every i = 1, 2, … , n. and* $\dim W^i = d^i$ *where* $d^i = d_1^i + d_2^i + \ldots + d_{k_i}^i$ *for every i = 1, 2, … , n.* $\dim W^1 + \dim W^2 + \ldots + \dim W^k = \dim V = (n_1, n_2, \ldots n_n)$ *i.e.,* $\dim W^1 = \dim W_1^1 + \dim W_2^1 + \ldots + \dim W_{k_1}^1 = n_1$, $\dim W^2 =$

$dim W_1^2 + dim\ W_2^2 + ... + dim\ W_{k_2}^1 = n_2$ *and so on and* $dim W^n =$ $dim W_1^n + dim\ W_2^n + ... + dim\ W_{k_n}^n = n_n$ .

The proof left as an exercise for the reader.

Now we proceed on to define the notion of n-annihilating polynominals.

Let T = $T_1 \cup T_2 \cup \ldots \cup T_n$ be a n-linear operator on a n-vector space V over the field F. If p(x) = $p_1(x) \cup p_2(x) \cup \ldots \cup p_n(x)$ be a n-polynominal in x with coefficients from F then p(T) = $p_1(T_1) \cup p_1(T_1) \cup \ldots \cup p_n(T_n)$ is again a n-linear operator on V. If q(x) = $q_1(x) \cup q_2(x) \cup \ldots \cup q_n(x)$ is another n-polynomial over F then

$$(p + q) (T) = p(T) + q(T) .$$
$$pq(T) = p(T) q(T)$$
$$= p_1(T) q_1(T) \cup p_2(T) q_2(T) \cup \ldots \cup p_n(T) q_n(T) .$$

Therefore the collection of n-polynomials p(x) which n-annihilate T in the sense that p(T) = 0, is a n ideal in the n-polynomial algebra F[x]. Clearly $L^n$ (V, V) is a n-linear space of dimension $(n_1^2, n_2^2 ..., n_n^2)$ where $n_i$ is the dimension of the vector space $V_i$ in V = $V_1 \cup V_2 \cup \ldots \cup V_n$. If we take in the n-linear operator T = $T_1 \cup T_2 \cup \ldots \cup T_n$, for each $T_i$ a $n_i^2 + 1$ power of $T_i$ for i = 1, 2, …, n then $C_0^i + C_1^i T_i + C_2^i T_i^2 + ... + C_{n_i^2}^i T_i^{n_i^2} = 0$ for some scalars $C_j^i$ not all zero, $1 \leq j \leq n_i^2$.

So the n-ideal of polynomials which n-annihilate T contains a non zero n-polynomial of n-degree $(n_1^2, n_2^2, ..., n_n^2)$ or less.

Now we define the notion of n-minimal polynomial for T = $T_1 \cup T_2 \cup \ldots \cup T_n$.

**DEFINITION 2.28:** *If T is a n-linear operator on a finite dimensional n-vector space V over the field F. The n-minimal polynomial for T =* $T_1 \cup T_2 \cup ... \cup T_n$ *is the unique n-monic generator of the n-ideals of polynomial over F which n-annihilate T, i.e., the n-monic generator of the n-ideals of polynomials over F which annihilate each* $T_i$ *for i = 1, 2, ..., n.*

*The term n-minimal comes from the fact that the n-generator of a polynomial n-ideal is characterized by being the n-monic polynomials each of minimum degree that is every ideal in the n-ideals; that is the n-minimal polynomial $p = p_1 \cup p_2 \cup \ldots \cup p_n$ for the n-linear operator T is uniquely determined by these three properties. In $p = p_1 \cup p_2 \cup \ldots \cup p_n$, each $p_i$ is a monic polynomial over the scalar field F, which we shortly call as the n-monic polynomial over F. $p(T) = 0$ implies $p_i(T_i) = 0$ for each i, $i = 1, 2, \ldots, n$; i.e., $p_1(T_1) \cup p_2(T_2) \cup \ldots \cup p_n(T_n) = 0 \cup 0 \cup \ldots \cup 0$. No n-polynomial over F which n-annihilates T has smaller degree than p, i.e., polynomial over F which annihilates $T_i$ has smaller degree than $p_i$ for each $i = 1, 2, \ldots, n$. If A is a n-mixed square matrix over F i.e., $A = A_1 \cup A_2 \cup \ldots \cup A_n$ is a n-mixed matrix where each $A_i$ is a $n_i \times n_i$ matrix over F, we define the n-minimal polynomial for A in an analogous way as unique n-monic generator ideal of all n-polynomials over F which n-annihilate A or annihilates $A_i$ for each i, $i = 1, 2, \ldots, n$.*

Similar results which hold good in case of linear vector spaces can be analogously extended to the case of n-vector spaces with proper and appropriate modifications.

The proof of the following interesting theorem can be obtained by any interested reader.

**THEOREM 2.25:** *Let $T = T_1 \cup T_2 \cup \ldots \cup T_n$ be a n-linear operator on a $(n_1, n_2, \ldots, n_n)$ finite dimensional n-vector space [or let $A = A_1 \cup A_2 \cup \ldots \cup A_n$, a n-mixed square matrix where each $A_i$ is a $n_i \times n_i$ matrix, $i = 1, 2, \ldots, n$] then n-characteristic and n-minimal polynomial for T[for A] have the same n-roots except for multiplicities.*

The Cayley-Hamilton theorem for n-linear operator T on the n-vector space V is stated, the proof is also left as an exercise for the reader.

**THEOREM 2.26:** (CAYLEY HAMILTON THEOREM FOR n-VECTOR SPACES) *Let $T = T_1 \cup T_2 \cup \ldots \cup T_n$ be a n-linear operator on a finite $(n_1, n_2, \ldots, n_n)$ dimensional n-vector space $V = V_1 \cup V_2 \cup$*

*... $\cup V_n$ over a field F. If $f = f_1 \cup f_2 \cup ... \cup f_n$ is the n-characteristic polynomial for $T = T_1 \cup T_2 \cup ... \cup T_n$ then $f(T) = 0 \cup 0 \cup ... \cup 0$ i.e., $f_1(T_1) \cup f_2(T_2) \cup ... \cup f_n(T_n) = 0 \cup 0 \cup ... \cup 0$; in other words the n-minimal polynomial divides the n-characteristic polynomial for T.*

We just give an hint of the proof.

*Hint:* Choose a n-ordered n-basis $\{(\alpha_1^1, \alpha_2^1, ..., \alpha_{n_1}^1) \cup (\alpha_1^2, \alpha_2^2, ..., \alpha_{n_2}^2) \cup \ldots \cup (\alpha_1^n, \alpha_2^n, ..., \alpha_{n_n}^n)\}$ for $V = V_1 \cup V_2 \cup \ldots \cup V_n$ and let $A = A_1^1 \cup A_1^2 \cup \ldots \cup A_n^n$ be the n matrix which represents $T = T_1 \cup T_2 \cup \ldots \cup T_n$ in the given n-basis. Then $T_i\alpha_k^i = \sum_{j=1}^{n_i} A_{jk}^i \alpha_j^i$; $1 \leq j \leq n_i$. This is true for each i; i.e., true for each $T_i$. Thus $p_k = \sum_{j=1}^{n_k} (\delta_{ij} T_k - A_{ji}^k I_k)\alpha_j^k = 0$, this equation being true for k = 1, 2, … , n, i.e., $P = P_1 \cup P_2 \cup \ldots \cup P_n$. Suppose $K = K_1 \cup K_2 \cup \ldots \cup K_n$ be a commutative n-ring with identity consisting of all n polynomials in $T = T_1 \cup T_2 \cup \ldots \cup T_n$. Let $B^1 \cup B^2 \cup \ldots \cup B^n$ be an element of

$$K^{n\times n} = K^{n_1\times n_1} \cup K^{n_2\times n_2} \cup ... \cup K^{n_n\times n_n}$$

with entries $B_{ij}^k = \delta_{ij} T_k - A_{ji}^k I_k$, k = 1, 2, …, n. We can show f(T) = det B i.e., $f_1(T_1) \cup f_2(T_2) \cup \ldots \cup f_n(T_n) = \det B^1 \cup \det B^2 \cup \ldots \cup \det B^n$ .

Using this hint the interested reader can prove the result.

Now we proceed on to define the notion of a n-subspace W of V to be n-invariant under T.

**DEFINITION 2.29:** *Let $V = V_1 \cup V_2 \cup ... \cup V_n$ be a n-vector space over F. $T = T_1 \cup T_2 \cup ... \cup T_n$ be a n-linear operator on V. If $W = W_1 \cup W_2 \cup ... \cup W_n$ is a n-subspace of V, we say that W is n-invariant under T if for each vector $\alpha = \alpha^1 \cup \alpha^2 \cup ... \cup \alpha^n$ in W the vector $T(\alpha)$ is in W i.e., each $T_i(\alpha^i) \in W_i$ for i = 1, 2,*

*… , n. i.e., T(W) is contained in W or this is the same as $T_i(W_i)$ is contained in $W_i$ for i = 1, 2, …, n.*

**LEMMA 2.5:** *Let V be a finite ($n_1$, $n_2$, …, $n_n$) dimensional n-vector space over the field F. Let T = $T_1 \cup T_2 \cup \ldots \cup T_n$ be a n-linear operator on V such that the n-minimal polynomial for T is a product of linear n-factors p = $p_1 \cup p_2 \cup \ldots \cup p_n$ where $p_i$ = $\left(x - C_1^i\right)^{r_1} \ldots \left(x - C_{k_i}^i\right)^{r_{k_i}}$ ; $C_j^i \in F$ , 1 ≤ j ≤ $k_i$, for i = 1, 2, …, m. Let W = $W^1 \cup W^2 \cup \ldots \cup W^n$ be a proper (W ≠ V) subspace of V where each $W^i = W_1^i + \ldots + W_{k_i}^i$ , i = 1, 2, … , n. which is n-invariant under T. There exists a vector α = $\alpha_1 \cup \alpha_2 \cup \ldots \cup \alpha_n$ in V such that α is not in W;*

$$(T - CI)\ \alpha = (T_1 - C^1I_1)\ \alpha_1 \cup (T_2 - C^2I_2)\alpha_2 \cup \ldots \cup (T_n - C^nI_n)\alpha_n$$

*is in W for some m-characteristic values $C^1 = (C_1^1, C_2^1 \ldots C_{k_1}^1)$, 1 ≤ $k_1 \le n_1$; $C^2 = (C_1^2, C_2^2, \ldots, C_{k_2}^2)$ and so on.*

The proof can be derived without much difficulty; infact very straight forward, using the working for each $T_i$: $V_i \to V_i$ and $W^i = W_1^i + \ldots + W_{k_i}^i$ , 1 ≤ $k_i$ ≤ $n_i$. When the result holds for every component of V and T it is true for the n-vector space and its n-linear operator T which is defined on V.

The following theorem on the n-diagonalizablily of the n-linear operator T on V is given below.

**THEOREM 2.27:** *Let V = $V_1 \cup V_2 \cup \ldots \cup V_n$ be a finite ($n_1$, $n_2$, … , $n_n$) dimensional n-vector space over the field F and let T = $T_1 \cup T_2 \cup \ldots \cup T_n$ be a n-linear operator on V. Then T is n diagonalizable if and only if the n-minimal polynomial for T has the form,*

$$p = \left\{\left(x - C_1^1\right)\ldots\left(x - C_{k_1}^1\right)\right\} \cup$$

$$\left\{\left(x-C_1^2\right)\left(x-C_2^2\right)...\left(x-C_{k_2}^2\right)\right\} \cup ... \cup$$
$$\left\{\left(x-C_1^n\right)\left(x-C_2^n\right)...\left(x-C_{k_n}^n\right)\right\}$$

*where $C_{j_i}^i$ are distinct elements of F (i.e., $C_1^1, C_2^1, ..., C_{k_1}^1$, forms a distinct set in F, $C_1^2, C_2^2, ..., C_{k_2}^2$ forms a distinct set in F, so on $C_1^n, C_2^n, ..., C_{k_n}^n$ forms a distinct set of F) .*

*Proof:* We know if $T = T_1 \cup T_2 \cup \ldots \cup T_n$ is n-diagonalizable its n-minimal polynomial is a n-product of distinct linear factors i.e., each $T_i: V_i \to V_i$ (where $V_i$ is a component of the n-vector space $V = V_1 \cup V_2 \cup \ldots \cup V_n$ and $T_i$ is a linear operator of $V_i$ and a component of T).

So we can say if $p_i = \left(x - C_1^i\right)\left(x - C_2^i\right)...\left(x - C_{k_i}^i\right)$ the minimal polynomial associated with the diagonalizable operator $T_i$ then the $p_i$ is a product of distinct linear factors. This is true for each i; i = 1, 2, …, n, Hence the claim. So to prove the converse, let $W = W^1 \cup W^2 \cup \ldots \cup W^n$ be the n-subspace spanned by all the n-characteristic n-vectors of T and suppose $W \neq V$ that is; each $W^i \neq V^i$ for i = 1, 2, …, n.

This implies we have a n-vector $\alpha = \alpha^1 \cup \alpha^2 \cup \ldots \cup \alpha^n$ not in W (i.e., each $\alpha^i \notin W^i$ for i = 1, 2, …, n.) and a n-characteristic value $C = C^1 \cup C^2 \cup \ldots \cup C^n$ of T such that the vector $\beta = (T - CI)\,\alpha$ lies in W i.e., $\beta = \beta^1 \cup \beta^2 \cup \ldots \cup \beta^n$ then $\beta^i = \left(T_i - C_j^i I_i\right)\alpha^i$ lies in $W^i$ ($1 \leq j \leq k_i$) this is true for each i, i = 1, 2, …, n. Since $\beta^i \in W^i$ we have $\beta^i = \beta_1^i + \beta_2^i + ... + \beta_{k_i}^i$ (true for each i, i = 1, 2, …, n) where $T_i\beta_j^i = C_j^i\beta_j^i$; $1 \leq j \leq k_i$ and i = 1, 2, …, n and hence the vector in $W^i$. $h^i(T_i)\beta^i = h^i(C_1^i)\beta_1^i + ... + h^i(C_{k_i}^i)\beta_{k_i}^i$ is in $W_i$ for every polynomial $h^i$; this is true for each i, i = 1, 2, …, n.

Now $p_i = (x - C_j^i)\ q_i$ for some polynomial $q_i$ also $q_i - q_i(C_j^i) = x - (C_j^i)h^i$ (this is true for each i, i = 1, 2, …, n).

We have $q_i(T_i)\alpha^i - q_i(C^i_j)\alpha^i = h^i(T_i)(T_i - C^i_j I_i)\alpha^i = h^i(T_i)\beta^i$, 1 $\leq$ i $\leq$ n. But $h^i(T_i)\beta^i$ is in $W_i$ (for each i) and since 0 = $p_i(T_i)$ $\alpha^i$ = $\left(T_i - C^i_j I_i\right) q_i(T_i)\alpha^i$, the vector $q_i(T_i)\alpha^i$ is in $W^i$.

Therefore $q_i(C^i_j)\alpha^i$ is in $W^i$. Since $\alpha^i$ is not $W^i$ we have $q_i(C^i_j)$ = 0 true for every i = 1, 2, …, n. This contradicts the fact $p_i$ has distinct roots for i = 1, 2, …, n. Hence the claim.

How ever we give an illustration of this theorem so that the reader can understand how it is applied in general.

***Example 2.19:*** Let V = $V_1 \cup V_2 \cup V_3$ where $V_1 = Q \times Q$, $V_2 = Q \times Q \times Q \times Q$ and $V_3 = Q \times Q \times Q$ i.e., V a 3-vector space over Q of finite dimension and 3-dimension (2, 4, 3) .Define T**:** V $\to$ V by T = $T_1 \cup T_2 \cup T_3 : V_1 \cup V_2 \cup V_3 \to V_1 \cup V_2 \cup V_3$ by $T_1$**:** $V_1 \to V_1$ defined by the related matrix

$$A_1 = \begin{bmatrix} 1 & 2 \\ 0 & 2 \end{bmatrix}.$$

$T_2$**:**$V_2 \to V_2$ defined by the related matrix

$$A_2 = \begin{bmatrix} 2 & 1 & 1 & 3 \\ 0 & 1 & 2 & 1 \\ 0 & 0 & 3 & 5 \\ 0 & 0 & 0 & 4 \end{bmatrix}$$

and $T_3$**:**$V_3 \to V_3$ defined by the related matrix

$$A_3 = \begin{bmatrix} 5 & -6 & -6 \\ -1 & 4 & 2 \\ 3 & -6 & -4 \end{bmatrix}.$$

The 3-matrix associated with T is given by

$$\begin{bmatrix} 1 & 2 \\ 0 & 2 \end{bmatrix} \cup \begin{bmatrix} 2 & 1 & 1 & 3 \\ 0 & 1 & 2 & 1 \\ 0 & 0 & 3 & 5 \\ 0 & 0 & 0 & 4 \end{bmatrix} \cup \begin{bmatrix} 5 & -6 & -6 \\ -1 & 4 & 2 \\ 3 & -6 & -4 \end{bmatrix}$$

that is 3-characterstic polynomial associated with T is given by

$$\begin{aligned} C &= (x-1)(x-2) \cup (x-2)(x-1)(x-3)(x-4) \cup (x-2)^2 (x-1) \\ &= C_1 \cup C_2 \cup C_3. \end{aligned}$$

The 3-minimal polynomial p is given by

$$\begin{aligned} p &= p_1 \cup p_2 \cup p_3 \\ &= (x-1)(x-2) \cup (x-2)(x-1)(x-3)(x-4) \cup (x-1)(x-2). \end{aligned}$$

Hence T is a 3-diagonalizable operator and the 3-diagonal 3-matrix associated with T is given by

$$D = \begin{bmatrix} 1 & 0 \\ 0 & 2 \end{bmatrix} \cup \begin{bmatrix} 2 & 0 & 0 & 0 \\ 0 & 1 & 0 & 0 \\ 0 & 0 & 3 & 0 \\ 0 & 0 & 0 & 4 \end{bmatrix} \cup \begin{bmatrix} 1 & 0 & 0 \\ 0 & 2 & 0 \\ 0 & 0 & 2 \end{bmatrix}.$$

Now we proceed on to describe the n-linear operator which is n-diagonalizable in the language of n-invariant direct sum decomposition.

**DEFINITION 2.30:** *Let V be a n-vector space over F, a n-projection of* $V = V_1 \cup V_2 \cup \ldots \cup V_n$ *is a n-linear operator* $E = E^1 \cup E^2 \cup \ldots \cup E^n$ *on V such that* $E^2 = E$ *i.e.,* $E^2 = (E^1)^2 \cup (E^2)^2 \cup \ldots \cup (E^n)^2$. *All properties associated with linear operators as projection can be analogously derived. Clearly if* $V = V_1 \cup V_2 \cup \ldots \cup V_n$ *and* $V = (W_1^1 \oplus \ldots \oplus W_{k_1}^1) \cup$ $(W_1^2 \oplus ... \oplus W_{k_2}^2) \cup ... \cup (W_1^n \oplus ... \oplus W_{k_n}^n)$ *then for each space* $W_j^i$*;* $1 \leq i \leq n$ *and* $1 \leq j \leq k_j$ *we can define* $E_j^i$ *an operator on* $V_i$ *such that if* $\alpha^i \in V_i$ *is of the form* $\alpha^i = \alpha_1^i + \alpha_2^i + ... + \alpha_{k_i}^i$ *with* $\alpha_j^i \in W_j^i$ *define* $E_j^i(\alpha^i) = \alpha_j^i$, $E_j^i$ *is a well defined rule; this is true for each i and j so*

$$E_1^1 + E_2^1 + ... + E_{k_1}^1 = E^1,$$

$$E_1^2 + E_2^2 + ... + E_{k_2}^2 = E^2, ...,$$

$$E_1^n + E_2^n + ... + E_{k_n}^n = E^n;$$

$E = E^1 \cup E^2 \cup ... \cup E^n$.

Now as in case of linear vector space we can in case of n-vector spaces derive the properties of projections.

**Theorem 2.28:** *Let* $V = V_1 \cup V_2 \cup ... \cup V_n$ *be a n-vector space over the field F. Suppose each* $V_i = W_1^i \oplus ... \oplus W_{k_i}^i$ *for i = 1, 2, ..., n i.e.,* $V = (W_1^1 \oplus ... \oplus W_{k_1}^1) \cup (W_1^2 \oplus W_2^2 \oplus ... \oplus W_{k_2}^2) \cup ... \cup (W_1^n \oplus W_2^n \oplus ... \oplus W_{k_n}^n)$ *then there exists* $(k_1 + k_2 + ... + k_n)$ *linear operators* $E_1^1, E_2^1, ..., E_{k_1}^1$, $E_1^2, E_2^2, ..., E_{k_2}^2$, ..., $E_1^n, E_2^n, ..., E_{k_n}^n$ *on the n-vector space V such that*

1. *Each* $E_{j_i}^i$ *is a projection, i = 1, 2, ..., n;* $1 \leq j_i \leq k_i$
2. $E_{j_i}^i . E_{j_k}^i = 0$ *if* $j_i \neq j_k$
3. $I^i = E_1^i + ... + E_{k_i}^i$ *i = 1, 2,..., n i.e.,* $I = I^1 \cup I^2 \cup ... \cup I^n$.
4. *range of* $E_j^i = W_j^i$ *for i = 1, 2, ..., n,* $1 \leq j \leq k_i$

*Conversely if* $E_1^1, E_2^1, ..., E_{k_1}^1$, $E_1^2, E_2^2, ..., E_{k_2}^2$, ..., $E_1^n, ..., E_{k_n}^n$ *are* $k_1 + k_2 + ... + k_n$ *linear operators on V which satisfy the condition 1, 2 and 3 and if we let* $W_j^i$ *be the range of* $E_j^i$ *then* $V = (W_1^1 \oplus ... \oplus W_{k_1}^1) \cup (W_2^2 \oplus ... \oplus W_{k_2}^2) \cup ... \cup (W_1^n \oplus ... \oplus W_{k_n}^n)$.

***Proof:*** Now to prove the converse statement we proceed as follows; from the basic definition and properties the condition 1 to 4 are true, which can be easily verified.

Suppose we have $E = E^1 \cup E^2 \cup \ldots \cup E^n$ where each $E^i$ is a $E_1^i, E_2^i, \ldots, E_{k_i}^i$ $k_i$ number of linear operators of $V_i$, $V_i$ a component of the n-vector space $V = V_1 \cup V_2 \cup \ldots \cup V_n$ and what we prove for i, is true for i = 1, 2, … , n.

Given E satisfies all the three conditions given in (1) (2) and (3) and if we let $W_j^i$ to be range of $E_j^i$ then certainly $V = W^1 \cup \ldots \cup W^n$ where $W^i = W_1^i \oplus \ldots \oplus W_{k_i}^i$ by condition (3) we have for $\alpha = \alpha^1 \cup \alpha^2 \cup \ldots \cup \alpha^n$,

$$\alpha = (E_1^1\alpha_1^1 + E_2^1\alpha_2^1 + \ldots + E_{k_1}^1\alpha_{k_1}^1) \cup$$
$$(E_1^2\alpha_1^2 + E_2^2\alpha_2^2 + \ldots + E_{k_2}^2\alpha_{k_2}^2) \cup \ldots \cup$$
$$(E_1^n\alpha_1^n + E_2^n\alpha_2^n + \ldots + E_{k_n}^n\alpha_{k_n}^n)$$

for each $\alpha \in V_j$ where each $I_i = E_1^i + \ldots + E_{k_i}^i$, i = 1, 2, … , n and $E_{j_p}^i . E_{j_k}^i = 0$ if $p \neq k$; $1 \leq j \leq k_i$ and $\alpha^i = \alpha_1^i + \ldots + \alpha_{k_i}^i$ true for i = 1, 2, …, n. This is true for each $\alpha^i \in V_i$ and hence for each $\alpha \in V$ and $E_j^i\alpha_j^i$ in $W^i$. This expression for each $\alpha^i$ is unique and hence each $\alpha$ is unique, because if $\alpha = (\alpha_1^1 + \ldots + \alpha_{k_1}^1) \cup (\alpha_1^2 + \alpha_2^2 + \ldots + \alpha_{k_2}^2) \cup \ldots \cup (\alpha_1^n + \alpha_2^n + \ldots + \alpha_{k_n}^n)$ is unique with each $\alpha^i \in W^i$, i.e., $\alpha_j^i \in W_j^i$. Suppose $\alpha_j^i = E_j^i\beta_j^i$ then from (1) and (2) we have

$$E_j^i\alpha^i = \sum_{k=1}^{k_i} E_j^i\alpha_{jk}^i$$
$$= \sum_{k=1}^{k_i} E_j^i E_k^i \beta_{jk}^i = (E_j^i)^2\beta_j^i$$
$$= E_j^i\beta_j^i = \alpha_j^i.$$

This is true for every i, i = 1, 2, …, n and every j, j = 1, 2, …, $k_i$. This proves each $V_i$ is direct sum of $W^i$, hence V is a direct sum $W_1^1,...,W_{k_1}^1$, … , $W_1^n, W_2^n, ..., W_{k_n}^n$. Hence the result.

Now we give a sketch of proof of the following theorem. However reader is expected to prove the theorem.

**THEOREM 2.29:** *Let $T = T_1 \cup T_2 \cup ... \cup T_n$ be a n-linear operator on the n-space $V = V_1 \cup V_2 \cup ... \cup V_n$ and let $W^1$, ..., $W^n$ and $E^1$, $E^2$, ..., $E^n$ be as in the above theorem. Then a necessary and sufficient condition that each n-subspace $W^i$ to be n-invariant under T (i.e., each $W_j^i$ invariant under $T_i$) is that T commutes with each of the projections $E^i$ i.e., $TE^i = E^iT$ for i = 1, 2, ..., m (i.e., each $T_i$ commutes with $E_j^i$ i.e., $T_i E_j^i = E_j^i T_i$, i = 1, 2, ..., n and j = 1, 2, ..., $k_i$).*

*Proof:* Suppose T commutes with each $E_j^i$ i.e., $T_i$ commutes with $E_j^i$ for j = 1, 2, …, $k_i$. This is true for each $T_i$ also. Let $\alpha = \alpha^1 \cup \alpha^2 \cup \ldots \cup \alpha^n$ with $\alpha_j^i \in W_j^i$, then $E_j^i \alpha_j^i = \alpha_j^i$ and for $T_i \alpha_j^i = T_i(E_j^i \alpha_j^i) = E_j^i T_i \alpha_j^i$ (since $T_i$ commutes with $E_j^i$ for j = 1, 2, …, $k_i$ and i = 1, 2, …, n) .This shows that $T_i \alpha_j^i$ is in the range of $E_j^i$ i.e., $W_j^i$ is invariant under $T_i$.

Assume now that each $W_j^i$ is invariant under $T_i$, 1 < j < $k_i$; i = 1, 2, …, n; we shall show that $T_i E_j^i = E_j^i T_i$ for every i, $1 \le i \le n$ and j = 1, 2, …, $k_i$. Let

$$\alpha^i \in V_i$$

$$\alpha^i = E_1^i \alpha^i + ... + E_{k_i}^i \alpha^i$$

$$T\alpha^i = TE_1^i \alpha^i + ... + TE_{k_i}^i \alpha^i.$$

Since $E_j^i \alpha^i$ is in $W_j^i$ which is invariant under $T_i$ we must have $T_i(E_j^i \alpha^i) = E_j^i \beta_j^i$ for some $\beta_j^i$.

Then $E_j^i T_i E_k^i \alpha^i = E_j^i E_k^i \beta_k^i$

$$= \begin{cases} 0 \text{ if } & k \neq j \\ E_j^i \beta_j^i \text{ if } & k = j. \end{cases}$$

$$E_j^i T_i \alpha^i = E_j^i T_i . E_1^i \alpha^i + \ldots + E_j^i T_i E_{k_i}^i \alpha^i$$

$$= E_j^i \beta_j^i$$

$$= T_i E_j^i \alpha^i .$$

This is true for each $\alpha^i \in V_i$ so $E_j^i T_i = T_i E_j^i$. This result is true for each i, i = 1, 2, … , n.

We now prove the main theorem which describes n-diagonalization of a n-linear operator.

**THEOREM 2.30:** *Let $T = T_1 \cup T_2 \cup \ldots \cup T_n$ be a n-linear operator on a finite dimensional n-vector space $V = V_1 \cup V_2 \cup \ldots \cup V_n$. If T is n-diagonalizable and if $(C_1^1, C_2^1, \ldots, C_{k_1}^1) \cup (C_1^2, C_2^2, \ldots, C_{k_2}^2) \cup \ldots \cup (C_1^n, C_2^n, \ldots, C_{k_n}^n)$ are n-characteristic values such that for each i, $C_1^i, C_2^i, \ldots, C_{k_i}^i$ are distinct characteristic values of $T_i$ for i = 1, 2, … , n, then their exists n-linear operators*

$$(E_1^1, E_2^1, \ldots, E_{k_1}^1), (E_1^2, E_2^2, \ldots, E_{k_2}^2), \ldots, (E_1^n, E_2^n, \ldots, E_{k_n}^n)$$

*on V such that*

1. $T = (C_1^1 E_1^1 + \ldots + C_{k_1}^1 E_{k_1}^1) \cup (C_1^2 E_1^2 + \ldots + C_{k_2}^2 E_{k_2}^2) \cup \ldots \cup (C_1^n E_1^n + \ldots + C_{k_n}^n E_{k_n}^n)$

2. $I = (E_1^1 + E_2^1 + \ldots + E_{k_1}^1) \cup \ldots \cup (E_1^n + \ldots + E_{k_n}^n) = I_1 \cup \ldots \cup I_n$

3. $E_k^i E_j^i = 0, j \neq k.$

4. $E_j^i E_j^i = E_j^i$

5. *The range of each $E_j^i$ is the characteristic space for $T_i$ associated with $C_j^i$.*

*Conversely if there exists ($k_1$, $k_2$, …, $k_n$) set of $k_i$ distinct n-scalars $C_1^i, C_2^i, ..., C_{k_i}^i$, i = 1, 2, …, n and $k_i$ distinct linear operators $E_1^i, E_2^i, ..., E_{k_i}^i$; i = 1, 2, …, n which satisfy conditions (1), (2) and (3) then $T_i$ is diagonalizable; hence $T = T_1 \cup T_2 \cup \ldots \cup T_n$ is n-diagonalizable. $C_1^i, C_2^i, ..., C_{k_i}^i$ are distinct characteristic values of $T_i$ for i = 1, 2, …, n and conditions (4) and (5) are satisfied.*

*Proof:* Suppose that T is n-diagonalizable i.e., each $T_i$ of T is diagonalizable with distinct characteristic values $(C_1^1, C_2^1, ..., C_{k_1}^1) \cup (C_1^2, C_2^2, ..., C_{k_2}^2) \cup \ldots \cup (C_1^n, C_2^n, ..., C_{k_n}^n)$, i.e., each set of $(C_1^i, C_2^i, ..., C_{k_i}^i)$ are distinct. Let $W_j^i$ be the space of characteristic vectors associated with the characteristic values $C_j^i$. As we have seen.

$$V = (W_1^1 \oplus \ldots \oplus W_{k_1}^1) \cup (W_1^2 \oplus \ldots \oplus W_{k_2}^2) \cup \ldots \cup (W_1^n \oplus \ldots \oplus W_{k_n}^n)$$

where each $V_i = W_1^i \oplus \ldots \oplus W_{k_i}^i$ for i = 1, 2, …, n.

Let $E_1^i, E_2^i, ..., E_{k_i}^i$ be the projections associated with this decomposition given in theorem. Then (2), (3), (4) and (5) are satisfied. To verify (1) we proceed as follows for each $\alpha = \alpha^1 \cup \alpha^2 \cup \ldots \cup \alpha^n$ in V; $\alpha^i \in V_i$; $\alpha^i = E_1^i\alpha + \ldots + E_{k_i}^i\alpha$ and so

$$T_i\alpha^i = TE_1^i\alpha^i + \ldots + TE_{k_i}^i\alpha^i$$

$$= C_1^iE_1^i\alpha^i + \ldots + C_{k_i}^iE_{k_i}^i\alpha^i.$$

In other words $T_i = C_1^iE_1^i + \ldots + C_{k_i}^iE_{k_i}^i$. Now suppose that we are given a n-linear operator $T = T_1 \cup T_2 \cup \ldots \cup T_n$ along with distinct n scalars $C^1 \cup C^2 \cup \ldots \cup C^n = C$ with scalar $C_j^i$ and non zero operator $E_j^i$ satisfying (1), (2) and (3) . This is true for each i = 1, 2, …, n and j = 1, 2, …, $k_i$. Since $E^i_j.E^i_k = 0$ when j ≠ k, we multiply both sides of

$$
\begin{aligned}
I &= I_1 \cup I_2 \cup \ldots \cup I_n \\
&= (E^1_1 + E^1_2 + \ldots + E^1_{k_1}) \cup (E^2_1 + E^2_2 + \ldots + E^2_{k_2}) \cup \ldots \cup \\
&\quad (E^n_1 + E^n_2 + \ldots + E^n_{k_n})
\end{aligned}
$$

by $E^1_t \cup E^2_t \cup \ldots \cup E^n_t$ and obtain immediately

$$E^1_t, E^2_t, \ldots, E^n_t = (E^1_t)^2 \cup (E^2_t)^2 \cup \ldots \cup (E^n_t)^2.$$

Multiplying

$$T = (C^1_1 E^1_1 + \ldots + C^1_{k_1} E^1_{k_1}) \cup \ldots \cup (C^n_1 E^n_1 + \ldots + C^n_{k_n} E^n_{k_n})$$

by $E^1_t \cup E^2_t \cup \ldots \cup E^n_t$ we have

$$T_1 E^1_t \cup T_2 E^2_t \cup \ldots \cup T_n E^n_t = C^1_t E^1_t \cup C^2_t E^2_t \cup \ldots \cup C^n_t E^n_t$$

which shows that any n-vector in the n range of $E^1_t \cup E^2_t \cup \ldots \cup E^n_t$ is in the n-null space of (T – CI) = $(T_1 - C^1_t I_1) \cup \ldots \cup (T_n - C^n_t I_n)$ where I = $I_1 \cup I_2 \cup \ldots \cup I_n$. Since we have assumed $E^1_t \cup E^2_t \cup \ldots \cup E^n_t \neq 0 \cup 0 \cup \ldots \cup 0$, this proves that there is a nonzero n-vector in the n-null space of

$$(T - CI) = (T_1 - C^1_t I_1) \cup \ldots \cup (T_n - C^n_t I_n)$$

i.e., that $C^i_t$ is a characteristic value of $T_i$ for each i, i = 1, 2, …, n; for if $C^i$ is any scalar then $(T_i - C^i I_i) = (C^i_1 - C^i)\, E^i_1 + \ldots + (C^i_{k_i} - C^i)\, E^i_{k_i}$ true for i = 1, 2, … , n so if $(T_i - C^i I_i)\alpha^i = 0$, we must have $(C^i_t - C^i)\, E^i_j \alpha^i = 0$. If $\alpha^i$ must be the zero vector then $E^i_j \alpha^i \neq 0$ for some j so that for this j we have $C^i_j - C^i = 0$.

Certainly $T_i$ is diagonalizable since we have shown that every non zero vector in the range of $E^i_j$ is a characteristic value of $T_i$ and the fact that $I_i = E^i_1 + \ldots + E^i_{k_i}$, shows that these characteristic vectors span $V_i$. This is true for each i, i = 1, 2, …, n. All that is to be shown is that the n-null space of (T – CI) = $(T_1 - C^1_k I_1) \cup (T_2 - C^2_k I_2) \cup \ldots \cup (T_n - C^n_k I_n)$ is exactly the n range of $E^1_k \cup ... \cup E^n_k$, but this is clear because T$\alpha$ = C$\alpha$ i.e., $T_i\alpha^i = C^i_j\alpha^i$, for each i, i = 1, 2, …, n. Thus

$$\bigcup_{i=1}^{k_i}\sum_{j=1}^{k_i}(C^i_j - C^i_k)E^i_j\alpha^i = 0 \text{ ; i.e.,}$$

$$\sum_{j=1}^{k_1}(C^1_j - C^1_k)E^1_j\alpha^1 \cup \sum_{j=1}^{k_2}(C^2_j - C^2_k)E^2_j\alpha^2 \cup ... \cup \sum_{j=1}^{k_n}(C^n_j - C^n_k)E^n_j\alpha^n$$

$$= 0 \cup 0 \cup \ldots \cup 0.$$

Hence $((C^i_j - C^i_k)E^i_j\alpha^i$ = 0 for each j; and each i = 1, 2, …, n and $E^i_j\alpha^i$ = 0, k ≠ j; for each i, i = 1, 2, …, n. Since $\alpha^i$ = $E^i_1\alpha^i + ... + E^i_{k_i}\alpha^i$ for each i and $E^i_j\alpha^i$ = 0 for j ≠ k we have $\alpha^i$ = $E^i_j\alpha^i$ which proves that $\alpha^i$ is the range of $E^i_j$. This is true for each i hence the claim.

We give the statement of the primary decomposition theorem for n-vector space V.

**THEOREM 2.31:** *Let T = $T_1 \cup T_2 \cup ... \cup T_n$ be a n-linear operator on the finite dimensional n-vector space V = $V_1 \cup V_2 \cup ... \cup V_n$ over the field F. Let p = $p^1 \cup p^2 \cup ... \cup p^n$ where $p^i$ = $p_1^{r_1^i} p_1^{r_2^i} ... p_{k_i}^{r_{k_i}^i}$, i = 1, 2, …, n. i.e.,*

$$p = p_{11}^{r_1^1} p_{12}^{r_2^1} ... p_{1k_1}^{r_{k_1}^1} \cup p_{21}^{r_1^2} p_{22}^{r_2^2} ... p_{2k_2}^{r_{k_2}^2} \cup ... \cup p_{11}^{r_1^n} p_{12}^{r_2^n} ... p_{nk_n}^{r_{k_n}^n}$$

*where $p_{ik}$ are distinct irreducible monic polynomials over F and the $r^i_j$ are positive integers. Let $W^i_j$ be the null space of $p_{ik}{}^{(T)^{r^i_k}}$, k = 1, 2, …, $k_i$; i = 1, 2, …, n then*

1. *V = $(W^1_1 \oplus ... \oplus W^1_{k_1})$ $\cup$ $(W^2_1 \oplus ... \oplus W^2_{k_2})$ $\cup$ … $\cup$ $(W^n_1 \oplus ... \oplus W^n_{k_n})$*
2. *each $W^i$ is invariant under $T_{ik}$, i = 1, 2, …, n, 1 $\leq r \leq k_i$.*

3. *if $T_{ij}$ is the operator induced on $W_j^i$ by $T_i$ then the minimal polynomial for $T_{ij}$ is $p_{ij}^{r_j^i}$, true for j = 1, 2, …, $k_i$ and i = 1, 2, …, n.*

Several interesting results can be found in this direction analogously.

Now we define the notion of n-diagonalizable part and n-nilpotent part of a n-linear operator T.

Given $V = V_1 \cup V_2 \cup \ldots \cup V_n$ is a n-vector space over the field F. $T = T_1 \cup T_2 \cup \ldots \cup T_n$ a n-linear operator on V. Suppose the n-minimal polynomial of T is the product of first degree polynomials, i.e., the case in which each $p_j^i$ is of the form $x - C_j^i$. Now range of $E_j^i$, for each $T_i$ in T is the null space $W_j^i$ of $((T_i - C_j^i I_i^i)^{r_j^i}$. This is true for each i, i = 1, 2, …, n. Put

$$\begin{aligned} D &= D_1 \cup D_2 \cup \ldots \cup D_n \\ &= (C_1^1E_1^1 + C_2^1E_2^1 + \ldots + C_{k_1}^1E_{k_1}^1) \cup \\ &\quad (C_1^2E_1^2 + C_2^2E_2^2 + \ldots + C_{k_2}^2E_{k_2}^2) \cup \ldots \cup \\ &\quad (C_1^nE_1^n + C_2^nE_2^n + \ldots + C_{k_n}^nE_{k_n}^n). \end{aligned}$$

Clearly D is n-diagonalizable operator which we define or call as the n-diagonalizable part of T. Let as consider N = T – D. Now

$$\begin{aligned} T &= [T_1E_1^1 + \ldots + T_1E_{k_1}^1] \cup [T_2E_1^2 + \ldots + T_2E_{k_2}^2] \cup \ldots \cup \\ &\quad [T_nE_1^n + \ldots + T_nE_{k_n}^n] \end{aligned}$$

$$\begin{aligned} D &= (C_1^1E_1^1 + C_2^1E_2^1 + \ldots + C_{k_1}^1E_{k_1}^1) \cup \\ &\quad (C_1^2E_1^2 + C_2^2E_2^2 + \ldots + C_{k_2}^2E_{k_2}^2) \cup \ldots \cup \\ &\quad (C_1^nE_1^n + C_2^nE_2^n + \ldots + C_{k_n}^nE_{k_n}^n) \end{aligned}$$

so

$$N = [(T_1 - C_1^1I_1^1)E_1^1 + (T_1 - C_2^1I_2^1)E_2^1 + \ldots + (T_1 - C_{k_1}^1I_{k_1}^1)E_{k_1}^1] \cup$$

$$[\,(T_2 - C_1^2 I_1^2)E_1^2 + (T_2 - C_2^2 I_2^2)E_2^2 + \ldots + (T_2 - C_{k_2}^2 I_{k_2}^2)E_{k_2}^2\,]$$
$$\cup \ldots \cup [\,(T_n - C_1^n I_1^n)E_1^n + \ldots + (T_n - C_{k_n}^n I_{k_n}^n)\,]\,.$$

Clearly

$$\begin{aligned} N^2 \;=\; & [\,(T_1 - C_1^1 I_1^1)^2 E_1^1 + \ldots + (T_1 - C_{k_1}^1 I_{k_1}^1)^2 E_{k_1}^1\,] \cup \\ & [\,(T_2 - C_1^2 I_1^2)^2 E_1^2 + \ldots + (T_2 - C_{k_2}^2 I_{k_2}^2)^2 E_{k_2}^2\,] \cup \ldots \cup \\ & [\,(T_n - C_1^n I_1^n)^2 E_1^n + \ldots + (T_n - C_{k_n}^n I_{k_n}^n)^2 E_{k_n}^n\,] \end{aligned}$$

and in general we have

$$\begin{aligned} N^r \;=\; & [\,(T_1 - C_1^1 I_1^1)^r E_1^1 + \ldots + (T_1 - C_{k_1}^1 I_{k_1}^1)^r E_{k_1}^1\,] \cup \ldots \cup \\ & [\,(T_n - C_1^n I_1^n)^r E_1^n + \ldots + (T_n - C_{k_1}^n I_{k_n}^n)^r E_{k_n}^n\,] \end{aligned}$$

where $r \geq (r^1, r^2, \ldots, r^n)$ i.e., $r > r^i$, $i = 1, 2, \ldots, n$ ( by misuse of notation) we have $N^r = 0$ because the n-operator $(T – CI)^r$ will be $(0 \cup 0 \cup \ldots \cup 0)$ i.e., each $(T_i - C_j^i I_i)\ r_j^i = 0$ where $r > r_j^i$ for $j = 1, 2, \ldots, k_i$ and $i = 1, 2, \ldots, n$.

Now we define a nilpotent n-linear operator T.

**DEFINITION 2.31:** *Let N be a n-linear operator on $V = V_1 \cup V_2 \cup \ldots \cup V_n$ we say N is n-nilpotent if there exists some positive integer r, $r > r^i$; $i = 1, 2, \ldots, n$ such that $N^r = 0$.*

*Note:* If $N = N_1 \cup N_2 \cup \ldots \cup N_n$ then $N_i: V_i \to V_i$ is of dimension $n_i$, $n_i \neq n_j$ if $i \neq j$ true for $i = 1, 2, \ldots, n$ so we may have $N_i^{r^i} = 0$, $i = 1, 2, \ldots, n$. We may not have $r^i = r^j$, if $i \neq j$; hence the claim.

Now we give only a sketch of the proof however the reader is expected to get the complete the proof using this sketch.

**THEOREM 2.32:** *Let $T = T_1 \cup T_2 \cup \ldots \cup T_n$ be a n-linear operator on a finite dimensional n-vector space $V = V_1 \cup V_2 \cup$*

*... $\cup V_n$ over the field F. Suppose the n-minimal polynomial for T decomposes over F in to product of n-linear polynomials, then there is a n-diagonalizable n-operator D on V and a n-nilpotent operator N on V such that*

*I. $T = D + N$*

*II. $DN = ND$.*

*The n-diagonalizable operator D and the n-nilpotent operator N are uniquely determined by (I) and (II) and each of them is a n-polynomial in T.*

*Proof:* We give only a sketch of the proof. However the interested reader can find a complete proof using this sketch.

Given $V = V_1 \cup V_2 \cup \ldots \cup V_n$ finite $(n_1, n_2, \ldots, n_k)$ dimensional a n-vector space. $T = T_1 \cup T_2 \cup \ldots \cup T_n$ a n-linear operator on T such that $T_i: V_i \to V_i$ for each i = 1, 2, …, n. We can write each $T_i = D_i + N_i$, a nilpotent part $N_i$ and a diagonalizable part $D_i$; i = 1, 2, …, n.

Thus

$$\begin{aligned} T &= T_1 \cup T_2 \cup \ldots \cup T_n \\ &= (N_1 + D_1) \cup (N_2 + D_2) \cup \ldots \cup (N_n + D_n) \\ &= (N_1 \cup N_2 \cup \ldots \cup N_n) + (D_1 \cup D_2 \cup \ldots \cup D_n) \end{aligned}$$

i.e., $T = N + D$ where $N = N_1 \cup N_2 \cup \ldots \cup N_n$ and $D = D_1 \cup D_2 \cup \ldots \cup D_n$. Since each $D_i$ and $N_i$ not only commute but are polynomials in $T_i$ we see D and N commute and are n-polynomials of T, as the result is true for each i, i = 1, 2, … , n. Suppose we have $T = D^1 + N^1$, i.e.,

$$\begin{aligned} T &= T_1 \cup T_2 \cup \ldots \cup T_n \\ &= (D_1^1 + N_1^1) \cup (D_2^1 + N_2^1) \cup \ldots \cup (D_n^1 + N_n^1) \\ &= D^1 + N^1 \end{aligned}$$

where $D^1$ is the n-diagonalizable part of T i.e., each $D_i^1$ is the diagonalizable part of $T_i$ for i = 1, 2, …, n and $N^1$ the n-nilpotent part of T i.e., each $N_i$ is the nilpotent part of $T_i$ for i = 1, 2, …, n.

Since each $D_i^1$ and $N_i^1$ commute for i = 1, 2, …, n we have $D^1$ and $N^1$ also n-commute with any n-polynomial in T. Hence in particular they commute with D and N.

Now we have $D + N = D^1 + N^1$ i.e., $D – D^1 = N^1 – N$ and these four n-operator commute with each other. Since D and $D^1$ are n-diagonalizable they commute and so $D – D^1$ is also n-diagonalizable.

Since both N and $N^1$ are n-nilpotent they n-commute and the operator $N – N^1$ is also n-nilpotent. Since $N – N^1 = D – D^1$ and $N – N^1$ is n-nilpotent we have $D – D^1$ the n-diagonalizable n-operator is also n-nilpotent.

Such an n-operator can only be a zero operator, for since it is n-nilpotent, the n-minimal polynomial for this n-operator is of the form $x^{r_1} \cup x^{r_2} \cup ... \cup x^{r_n}$ with $x^{r_i} = 0$ for appropriate $m_i \geq r_i$, i = 1, 2, …, n. But since the n-operator is n-diagonalizable the n-minimal polynomial cannot have repeated n-roots hence each $r_i = 1$ and the n-minimal polynomial is simple $x \cup x \cup \ldots \cup x$ which confirms the operator is zero. Thus we have $D = D^1$ and $N = N^1$.

The interested reader is expected to derive analogous results when F is the field of complex numbers.

Now we proceed on to work with n-characteristic values n-characteristic vectors of a special n-linear n-operator on V.

Given V is a n-vector space say of finite dimension, $V = V_1 \cup V_2 \cup \ldots \cup V_n$ of dimension $(n_1, n_1, \ldots, n_n)$ defined over the field F. Let $T = T_1 \cup T_2 \cup \ldots \cup T_n$ be a special n-linear operator on V; i.e., $T_i: V_i \to V_i$ for each i, i = 1, 2, …, n.

We say $C = (C_1 \cup C_2 \cup \ldots \cup C_n)$ is a n-characteristic value of T if some n-vector $\alpha = \alpha_1 \cup \alpha_2 \cup \ldots \cup \alpha_n$ we have $T\alpha = C\alpha$, i.e., $T\alpha = (T_1 \cup T_2 \cup \ldots \cup T_n)(\alpha_1 \cup \alpha_2 \cup \ldots \cup \alpha_n) = (C_1 \cup C_2 \cup \ldots \cup C_n)(\alpha_1 \cup \alpha_2 \cup \ldots \cup \alpha_n)$ i.e., $T_o T = T_1\alpha_1 \cup T_2\alpha_2 \cup \ldots \cup T_n\alpha_n = C_1\alpha_1 \cup C_2\alpha_2 \cup \ldots \cup C_n\alpha_n$, i.e., each $T_i\alpha_i = C_i\alpha_i$ for i = 1, 2, …, n.

Here $\alpha = \alpha_1 \cup \alpha_2 \cup \ldots \cup \alpha_n$ is defined to be the n-characteristic vector of T. The collection of all $\alpha$ such that $T\alpha = C\alpha$ is called the n-characteristic space associated with C.

We shall illustrate the working of the n characteristic values, n-characteristic vectors associated with aT.

***Example2.20:*** Let $V = V_1 \cup V_2 \cup V_3$ be 3 vector space over Q where $V_1 = Q \times Q \times Q$, $V_2 = Q \times Q$ and $V_3 = Q \times Q \times Q \times Q$ are vector spaces over Q of dimensions 3, 2 and 4 respectively i.e., V is of (3, 2, 4) dimension. Define T: V $\rightarrow$ V where the 3-matrix associated with T is given by

$$A = A_1 \cup A_2 \cup A_3$$

$$= \begin{bmatrix} 3 & 0 & 2 \\ 0 & 1 & 5 \\ 0 & 0 & 7 \end{bmatrix} \cup \begin{bmatrix} 1 & 2 \\ 0 & 3 \end{bmatrix} \cup \begin{bmatrix} 1 & 0 & 2 & 1 \\ 0 & 2 & 5 & 0 \\ 0 & 0 & 3 & 7 \\ 0 & 0 & 0 & 4 \end{bmatrix}.$$

Now we will determine the 3-characterstic values associated with T. The n-characteristic polynomial

$$p = \begin{bmatrix} x-3 & 0 & -2 \\ 0 & x-1 & -5 \\ 0 & 0 & x-7 \end{bmatrix} \cup \begin{bmatrix} x-1 & -2 \\ 0 & x-3 \end{bmatrix} \cup$$

$$\begin{bmatrix} x-1 & 0 & -2 & -1 \\ 0 & x-2 & -5 & 0 \\ 0 & 0 & x-3 & 0 \\ 0 & 0 & 0 & x-4 \end{bmatrix}$$

$$= (x-3)(x-1)(x-7) \cup (x-1)(x-3) \cup (x-1)(x-2)(x-3)(x-4).$$

Thus the 3- characteristic values of $A = A_1 \cup A_2 \cup A_3$ are {3, 1, 7} $\cup$ {1, 3} $\cup$ {1, 2, 3, 4}. One can find the 3-characteristic values as in case of usual vector spaces and their set theoretic union will give 3-row mixed vector, which will be 48 in number as we have 48 choices for the 3-characterstic values as {3} $\cup$ {1} $\cup$ {1}, {3} $\cup$ {1} $\cup$ {2}, {3} $\cup$ {1} $\cup$ {3}, {3} $\cup$ {1} $\cup$ {4} so on and {7} $\cup$ {3} $\cup$ {4}.

Now having seen the working of 3-characteristic values we just recall in case of matrices A we say A is orthogonal if $AA^t = I$. Further A is anti orthogonal if $AA^t = -I$.

Now we for the first time define the notion of n-orthogonal matrices and n-anti orthogonal matrices.

**DEFINITION 2.32:** *Let $A = (A_1 \cup A_2 \cup \ldots \cup A_n)$ be a n-matrix.*

$$A^t = (A_1 \cup \ldots \cup A_n)^t = A_1^t \cup A_2^t \cup \ldots \cup A_n^t .$$

$$AA^t = A_1A_1^t \cup A_2A_2^t \cup \ldots \cup A_nA_n^t .$$

*We say A is n-orthogonal if and only if $AA^t = I_1 \cup I_2 \cup \ldots \cup I_n$ where $I_j$ is the identity matrix, i.e., if $A = A_1 \cup A_2 \cup \ldots \cup A_n$ is $m_i \times n_i$ matrix i = 1, 2, …, n; then $AA^t = I_1 \cup I_2 \cup \ldots \cup I_n$ is such that $I_j$ is a $m_j \times m_j$ identity matrix, j = 1, 2, …, n. We say A is anti orthogonal if and only if $AA^t = (-I_1) \cup (-I_2) \cup \ldots \cup (-I_n)$ where $I_j$ is $m_j \times m_j$ identity matrix i.e., if*

$$I = \begin{bmatrix} 1 & 0 & 0 & 0 \\ 0 & 1 & 0 & 0 \\ 0 & 0 & 1 & 0 \\ 0 & 0 & 0 & 1 \end{bmatrix}$$

*then*

$$-I = \begin{bmatrix} -1 & 0 & 0 & 0 \\ 0 & -1 & 0 & 0 \\ 0 & 0 & -1 & 0 \\ 0 & 0 & 0 & -1 \end{bmatrix}.$$

*Now we say $AA^t$ is n-semi orthogonal if $AA^t = B_1 \cup B_2 \cup \ldots \cup B_n$ ; some of the $B_i$'s are identity matrices and some are not identity matrices on similar lines we define n-semi anti orthogonal if in $AA^t = C_1 \cup C_2 \cup \ldots \cup C_n$ some $C_i$'s are $-I_i$ and some are not $-I_j$.*

It is not a very difficult task for the reader can easily get examples of these 4 types of n-matrices.

Chapter Three

# APPLICATIONS OF n-LINEAR ALGEBRA OF TYPE I

In this chapter we just introduce the applications of the n-linear algebras of type I. We just recall the notion of Markov bichains and indicate the applications of vector bispaces and linear bialgebras in Markov bioprocess. For this we have to first define the notion of Markov biprocess and its implications to linear bialgebra / bivector spaces. We may call it as Markov biprocess or Markov bichains.

Suppose a physical or mathematical system is such that at any moment it occupies two of the finite number of states (Incase of one of the finite number of states we apply Markov chains or the Markov process). For example say about a individuals emotional states like happy, sad etc., suppose a system move with time from two states or a pair of states to another pair of states; let us construct a schedule of observation times and a record of states of the system at these times. If we find the transition from one pair of state to another pair of state is not predetermined but rather can only be specified in terms of certain probabilities depending on the previous history of the system then the biprocess is called a stochastic biprocess. If in addition these transition probabilities depend only on the immediate history of the system; that is if the state of the system at any observation is dependent only on its state at the immediately proceeding observations then the process is called Markov biprocess or Markov bichain.

The bitransition probability $p_{ij} = p^1_{i_1 j_1} \cup p^2_{i_2 j_2}$ (i, j = 1, 2,…, k) is the probabilities that if the system is in state j = $(j_1, j_2)$ at any observation, it will be in state i = $(i_1, i_2)$ at the next observation. A transition matrix

$$P = [p_{ij}] = \left[p^1_{i_1 j_1}\right] \cup \left[p^2_{i_2 j_2}\right]$$

is any square bimatrix with non negative entries for which the bicolumn sum is $1 \cup 1$. A probability bivector is a column bivector with non negative entries whose sum is $1 \cup 1$.

The probability bivectors are said to be the state bivectors of the Markov biprocess. If $P = P_1 \cup P_2$ is the transition bimatrix of the Markov biprocess and $x^n = x_1^n \cup x_2^n$ is the state bivector at the $n^{th}$ observation then $x^{(n+1)} = P\, x^{(n)}$ and thus $x_1^{(n+1)} \cup x_2^{(n+1)} = P_1\, x_1^{(n)} \cup P_2\, x_2^{(n)}$. Thus Markov bichains find all its applications in bivector spaces and linear bialgebras.

Now we proceed onto define the new notion of Markov n-chain $n \geq 2$. Suppose a physical or a mathematical system is such that at any time it can occupy a finite number of states; when we view them as stochastic biprocess or Markov bichains when we make an assumption that the system moves with time from one state to another so that a schedule of observation times keeps the states of the system at these times. But when we tackle real world problems, say even for simplicity; emotions of a person may be very unpredictable depending largely on the situation and the mood of the person and its relation with another so such study cannot come under Markov chains. Even more is the complicated situation when the mood of a boss with subordinates; where mood of a person with a n number of persons and with varying emotions at a time and in such cases more than one emotion is experienced by a person and such states cannot be included and given as a next set of observation.

These changes and several feelings say at least n at a time ($n \geq 2$) will largely affect the transition n-matrix

$$P = P_1 \cup \ldots \cup P_n = \left[p^1_{i_1 j_1}\right] \cup \ldots \cup \left[p^n_{i_n j_n}\right]$$

with non negative entries which we will explain shortly. We indicate how n-vector spaces and n-linear algebras are used in Markov n-process ($n \geq 2$), when $n = 2$ the study is termed as Markov bioprocess. We first define Markov n-process and its implications to linear n-algebra and n-vector spaces; which we may call as Markov n-process and Markov n-chains.

Suppose a physical or a mathematical system is such that at any moment it occupies two or more finite number of states (in case of one of the finite number of states we apply Markov chains or the Markov process; in case of two of the finite number of state we apply Markov bichains or Markov biprocess). For example individual emotional states; happy, sad, cold, angry etc. suppose a system move with time from n states or a n tuple of states to another n-tuple of states; let us construct a schedule of observation times and a record of states of the system at these times. If we find the transition from n-tuple of states to another n-tuple of states not predetermined but rather can only be specified in terms of certain probabilities depending on the previous history of the system then the n-process is called a stochastic n-process. If in addition these transition probabilities depend only on the immediate history of the system that is if the state of the system at any observation is dependent only on its state at immediately proceeding observations then the process is called Markov n-process or Markov n-chain.

The n-transition probability

$$p_{ij} = p^1_{i_1 j_1} \cup p^2_{i_2 j_2} \cup \ldots \cup p^n_{i_n j_n}$$

$i, j = 1, 2, \ldots, K$ is the probabilities that if the system is in state $j = (j_1, j_2, \ldots, j_n)$ at any observation it will be in state $i = (i_1, i_2, \ldots, i_n)$ at the next observation.

A transition matrix associated with it is

$$P = [p_{ij}] = [p^1_{i_1 j_1}] \cup \ldots \cup [p^n_{i_n j_n}]$$

is a square n-matrix with non negative entries for all of which n-column sum is $(1 \cup \ldots \cup 1)$. A probability n-vector is a column n-vector with non negative entries whose sum is $1 \cup \ldots \cup 1$.

The probability n-vectors are said to be the state n-vectors of the Markov n-process. If $P = P_1 \cup \ldots \cup P_n$ is the transition n-matrix of the Markov n-process and $x^m = x_1^m \cup \ldots \cup x_n^m$ is the state n-vector at the $m^{th}$ observation then $x^{(m+1)} = Px^{(m)}$ and thus $x_1^{m+1} \cup \ldots \cup x_n^{m+1} = P_1(x_1^{(m)}) \cup \ldots \cup P_n x_n^{(m)}$. Thus Markov n-chains find all its applications in n-vector spaces and linear n-algebras. (n-linear algebras).

***Example 3.1****:* (Random Walk): A random walk by n persons on the real lines i.e. lines parallel to x axis is a Markov n-chain such that $p^1_{j_1k_1} \cup \ldots \cup p^n_{j_nk_n} = 0 \cup \ldots \cup 0$ if $k_t = j_t - 1$ or $j_t + 1$, t = 1, 2, …, n. Transition is possible only to neighbouring states from j to j – 1 and j + 1. Here state n-space is $S = S_1 \cup \ldots \cup S_n$ where $S_i = \{ \ldots -3 -2 -1\ 0\ 1\ 2\ 3 \ldots\}$; i = 1, 2, …, n.

The following theorem is direct.

**THEOREM 3.1:** *The Markov n-chain* $\{X_{m_1}; m_1 \geq 0\} \cup \ldots \cup \{X_{m_n}; m_n \geq 0\}$ *is completely determined by the transition n-matrix* $P = P_1 \cup \ldots \cup P_n$ *and the initial n-distribution* $\{P^1_{K_1}\} \cup \ldots \cup \{P^n_{K_n}\}$ *defined as*

$$P_1[X_0^1 = K_1] \cup \ldots \cup P_n[X_0^n = K_n]$$
$$= p_{K_1} \cup \ldots \cup p_{K_n} \geq 0 \cup \ldots \cup 0$$

*and*

$$\sum_{K_1 \in S_1} p_{k_1} \cup \ldots \cup \sum_{K_n \in S_n} p_{k_n} = 1 \cup \ldots \cup 1.$$

The proof is similar to Markov chain.

The n vector $u = (u_1^1 \ldots u_{n_1}^1) \cup \ldots \cup (u_1^n \ldots u_{n_n}^n)$ is called a probability n-vector if the components are non negative and their sum is one.

The square n-matrix $P = P_1 \cup \ldots \cup P_n = (p^1_{i_1 j_1}) \cup \ldots \cup (p^n_{i_n j_n})$ is called a stochastic n-matrix if each of the n-row probability n-vector i.e. each element of $P_i$ is non negative and the sum of the elements in each row of $P_i$ is one for i = 1, 2, …, n.

We illustrate this by a simple example.

***Example 3.2**:* Let

$$P = \begin{bmatrix} 1 & 0 & 0 \\ 1/3 & 1/6 & 1/2 \\ 1/4 & 0 & 3/4 \end{bmatrix} \cup \begin{bmatrix} 1 & 0 \\ 3/7 & 4/7 \end{bmatrix} \cup \begin{bmatrix} 0 & 1/2 & 0 & 1/2 \\ 1 & 0 & 0 & 0 \\ 1/4 & 1/4 & 1/2 & 0 \\ 3/7 & 2/7 & 1/7 & 1/7 \end{bmatrix}$$

be a stochastic 3-matrix.

The transition n-matrix P of a Markov n-chain (m-n-C) is a stochastic n-matrix. A stochastic n-matrix $A = A_1 \cup \ldots \cup A_n$ is said to be n-regular if all the entries of some power of each $A_i$ i.e. $A_i^{m_i}$ is positive, $m_i$'s positive integer for every i, i = 1, 2 …, n; i.e. $(m_1, \ldots, m_n) > (1, 1, \ldots, 1)$. $A^m = A_1^{m_1} \cup \ldots \cup A_1^{m_n}$; m = $(m_1, \ldots, m_n)$ with $A_i^{m_i} > 0$ for each i so that we state $A^m > (0 \cup \ldots \cup 0)$. It is easily verified that if $P = P_1 \cup \ldots \cup P_n$ is a stochastic n-matrix then $P^m$ is also a stochastic n-matrix for all m > (1, 1, …, 1). Is P a stochastic n-matrix if $P^n$ is a stochastic n-matrix?

Prove (1, …, 1) is a n-eigen value of a stochastic n-matrix i.e. if $A = A_1 \cup \ldots \cup A_n$; $|\lambda I - A| = 0 \cup \ldots \cup 0 \Rightarrow \lambda = (1, \ldots, 1)$ if $|\lambda_1 I_1 - A_1| \cup \ldots \cup |\lambda_n I_n - A_n| = 0 \cup 0 \cup \ldots \cup 0 \cup 0$, implies $\lambda = (\lambda_1, \ldots, \lambda_n) = (1, \ldots, 1)$. We define n-independent trials analogous to independent trials if

$$P = P_1 \cup \ldots \cup P_n$$

and

$$P^m = P = P_1^m \cup \ldots \cup P_n^m$$
$$= P_1 \cup \ldots \cup P_n$$

for all m ≥ (1, …, 1) where $p^t_{i_t j_t} = p^t_{j_t}$ for t = 1, 2, …, n i.e. all the rows of each $P_t$ is the same then we say P is an n-independent trial.

We can also define the notion of Bernoulli n trials. We just depict n-random walk with absorbing barriers. Let the possible n-states be $(E^1_0, E^1, \ldots, E^1_{K_1}) \cup \ldots \cup (E^n_0, E^n_1, \ldots, E^n_{K_n})$.

Consider the n-matrix of transition n-possibilities

$$P = P^1_{i_1 j_1} \cup \ldots \cup P^n_{i_n j_n}$$

$$= \begin{bmatrix} 1 & 0 & 0 & \ldots & 0 \\ q_1 & 0 & p_1 & \ldots & 0 \\ 0 & q_1 & 0 & p_1 & 0 \\ 0 & & \ldots & & q_1 \\ 0 & \ldots & & 0 & 1 \end{bmatrix}_{K_1 \times K_1} \cup \ldots \cup \begin{bmatrix} 1 & 0 & 0 & \ldots & 0 \\ q_n & 0 & p_n & \ldots & 0 \\ 0 & q_n & 0 & p_n & \vdots \\ 0 & & \ldots & 0 & q_n \\ 0 & \ldots & & 0 & 1 \end{bmatrix}_{K_n \times K_n}$$

From each of the interior n states

$$\{E^1_1, \ldots, E^1_{K_{1-1}}\} \cup \ldots \cup \{E^n_1, \ldots, E^n_{K_{n-1}}\},$$

n-transmission are possible to the right and left neighbour with $(p_t)_{i_t, i_t+1} = p_t$, $(p_t)_{i_t, i_t-1} = q_t$; t = 1, 2, …, n.

However no n-transition is possible from either $E_0 = (E^1_0 \cup \ldots \cup E^n_0)$ and $E_K = [E^1_K \cup \ldots \cup E^n_K]$ to any other n-state.

This n-system may move from one n-state to another but once $E_0$ or $E_K$ is reached the n-system stays there permanently.

Now we describe random walk with reflecting barriers.

Let $P = P^1_{i_1 j_1} \cup \ldots \cup P^n_{i_n j_n}$ be a n-matrix with

$$P = \begin{bmatrix} q^1 & p^1 & 0 & \dots & 0 & 0 & 0 \\ q^1 & 0 & p^1 & \dots & 0 & 0 & 0 \\ 0 & q_1 & 0 & p^1 & 0 & 0 & 0 \\ 0 & 0 & 0 & \dots & q^1 & 0 & p^1 \\ 0 & 0 & 0 & \dots & 0 & q^1 & p^1 \end{bmatrix}$$

$$\cup \dots \cup \begin{bmatrix} q^n & p^n & 0 & \dots & 0 & 0 & 0 \\ q^n & 0 & p^n & \dots & 0 & 0 & 0 \\ 0 & q^n & 0 & p^n & 0 & 0 & 0 \\ 0 & 0 & 0 & \dots & q^n & 0 & p^n \\ 0 & 0 & 0 & \dots & 0 & q^n & p^n \end{bmatrix}$$

$p^t$ and $q^t$ for t = 1, 2, …, n is defined by

$$P^t_{ij} = P^t\,(X^t_{n_t} = j_t \,|X^t_{n_t-1} = i_t) = \begin{cases} p^t & \text{if } j_t = i_t + 1 \\ q^t & \text{if } j_t = 0 \\ 0 & \text{otherwise} \end{cases}$$

true for t = 1, 2, …, n.

It may be possible that $p^t_{i_t j_t} = 0$, $p^{t^{(2)}}_{i_t j_t} = 0$ but $p^{t^{(3)}}_{i_t j_t} > 0$. We say the state $j_t$ is accessible from state $i_t$ if $P^{t^{(n)}}_{i_t j_t} > 0$ for some n > 0.

In notation $i_t \to j_t$ i.e. $i_t$ leads to $j_t$. If $i_t \to j_t$ and $j_t \to i_t$ then $i_t$ and $j_t$ communicate and we denote it by $i_t \leftrightarrow j_t$, if this happens we say they n-communicate. If only some of them communicate and others do not communicate we say the n-system semi communicates.

Here

$$P^{t^{(n)}}_{i_t j_t} = \begin{cases} q^t p^{t^{j_t}} & \text{for } j_t = 0, 1, 2, \dots, i_t + n_t - 1 \\ j^t\ p^{j_t} & \text{for } j_t = j_t + n_t \\ 0 & \text{otherwise.} \end{cases}$$

The state $i_t$ is essential if $i_t \to j_t$ implies $i_t \leftarrow j_t$ i.e. if any state $j_t$ is accessible from $i_t$ then $i_t$ is accessible from that state, true for t

= 1, 2, …, n. Let $\mathfrak{I} = \mathfrak{I}_1 \cup \ldots \cup \mathfrak{I}_n$ denote the set of all essential n state i.e. each $\mathfrak{I}_t$ denotes the set of all essential states, t = 1, 2, …, n. States that not n-essential are called n-inessential. We have semi essential if a few of the $\mathfrak{I}_t$'s are essential. We have semi essential state as m-essential state where m < n and only m out of the n states are essential rest inessential or n – m inessential state.

A Markov n-chain is called n-irreducible (or n-ergodic) if there is only one n communicating class i.e. all states n-communicate with each other or every n-state can be reached from every other n-state.

A n-subset $c = c_1 \cup \ldots \cup c_n$ of $S = S_1 \cup \ldots \cup S_n$ is said to be closed (or n-transient) if it is impossible to leave c in one step i.e. $p_{ij} = 0 \cup \ldots \cup 0$, i.e. $p^1_{i_1 j_1} \cup \ldots \cup p^n_{i_n j_n} = 0 \cup \ldots \cup 0$ for all i $\in$ c i.e. $(i_1, \ldots, i_n) \in c_1 \cup \ldots \cup c_n$ and all $(j_1, \ldots, j_n) \notin c$ for all $i_t \in c_t$ and all $j_t \notin c_t$; t = 1, 2, …, n.

We say a n-subset $c = c_1 \cup \ldots \cup c_n$ of $S = S_1 \cup \ldots \cup S_n$ is semi n-closed (or semi n-transient) if it is impossible to leave (only m of the) $c_t$'s, $1 \leq t \leq n$, m < n in one state; i.e. $p^t_{i_t j_t} = 0$ for all $i_t \in c_t$, and for all $j_t \in c_t$. We call this also m-closed (m < n) or m-transient, m = 1, 2, …, n –1. If m = n – 1 we call c to be hyper n-closed (or hyper n-transient).

A Markov n-chain is n-irreducible if the only n-closed set in S is S itself i.e., there is no n-closed set other than the set all of n states.

We say a Markov n-chain is semi irreducible or m-irreducible (m < n) if the closed sets in $S = S_1 \cup \ldots \cup S_n$ are m in number from the n-states $\{S_1, \ldots, S_n\}$, m < n. If m = n – 1 then we say the Markov n-chain is hyper n irreducible.

A single n-state $\{K_1, \ldots, K_n\}$ forming a closed n-set is called n-absorbing (n-trapping) i.e., a n-state such that the n-system remains in that state once it enters there. Thus a n-state $\{K_1, \ldots, K_n\}$ is n absorbing if the $\{K_1^{th}, \ldots, K_n^{th}\}$ rows of the transition n-matrix $P = P_1 \cup \ldots \cup P_n$ has 1 on the main n-diagonal and 0 else where.

***Example 3.3:*** Let P = $P_1 \cup \ldots \cup P_4$ be a transition 4 matrix given by

$$P = \begin{bmatrix} 0 & 0 & 1/2 & 0 & 1/2 & 0 & 0 \\ 1 & 0 & 0 & 0 & 0 & 0 & 0 \\ 0 & 0 & 0 & 1/3 & 0 & 2/3 & 0 \\ 0 & 0 & 0 & 1 & 0 & 0 & 0 \\ 0 & 1/7 & 5/7 & 0 & 0 & 0 & 1/7 \\ 0 & 0 & 0 & 5/8 & 3/8 & 0 & 0 \\ 0 & 3/5 & 0 & 1/5 & 0 & 1/5 & 0 \end{bmatrix} \cup$$

$$\begin{bmatrix} 0 & 1/2 & 0 & 1/2 & 0 \\ 0 & 0 & 1/3 & 1/3 & 1/3 \\ 0 & 0 & 1 & 0 & 0 \\ 1/5 & 4/5 & 0 & 0 & 0 \\ 7/9 & 0 & 0 & 0 & 2/9 \end{bmatrix} \cup \begin{bmatrix} 1 & 0 & 0 & 0 \\ 0 & 1/9 & 0 & 8/9 \\ 1/8 & 0 & 7/8 & 0 \\ 1/4 & 0 & 1/2 & 1/4 \end{bmatrix} \cup$$

$$\begin{bmatrix} 0 & 0 & 0 & 1 & 0 & 0 \\ 1/2 & 0 & 1/2 & 0 & 0 & 0 \\ 1/4 & 0 & 0 & 3/4 & 0 & 0 \\ 0 & 0 & 0 & 0 & 7/8 & 1/8 \\ 0 & 1 & 0 & 0 & 0 & 0 \\ 0 & 0 & 0 & 0 & 0 & 1 \end{bmatrix}.$$

Clearly the n-absorbing state is (4, 3, 1, 6).

Several interesting results true in case of M C can be proved for M – n – C with appropriate changes and suitable modifications.

Now we briefly describe the method for spectral m-decomposition (m ≥ 2). Let $P = P_1 \cup \ldots \cup P_m$ be a $N \times N$ m-matrix with m set of latent roots $\lambda_1^1 \ldots \lambda_{N_1}^1, \lambda_1^2 \ldots \lambda_{N_2}^2, \ldots, \lambda_1^m \ldots \lambda_{N_m}^m$ all distinct and simple i.e. each set of latent roots $\{\lambda_1^t \ldots \lambda_{N_t}^t\}$ are all distinct and simple for t = 1, 2, …, m; then

$$(P_1 - \lambda_{i_1}^1 I_1)\, U_{i_1}^1 \cup \ldots \cup (P_m - \lambda_{i_m}^m I_m)\, U_{i_m}^m = 0 \cup \ldots \cup 0$$

for the n-column latent n-vector $U_{i_1}^1 \cup \ldots \cup U_{i_m}^m$ and

$$V_{i_1}^{1'}\,(P_1 - \lambda_{i_1}^1 I) \cup \ldots \cup V_{i_m}^{m'}\,(P_m - \lambda_{i_m}^m I) = 0 \cup \ldots \cup 0$$

for the row latent n-vector $V_{i_1}^1 \cup \ldots \cup V_{i_m}^m$.

$$A_{i_1}^1 \cup \ldots \cup A_{i_m}^m = U_{i_1}^1\, V_{i_1}^{1'} \cup \ldots \cup U_{i_m}^m\, V_{i_m}^{m'}$$

are called m latent or m-spectral m-matrix associated with $(\lambda_{i_1}^1, \ldots, \lambda_{i_m}^m)$; $i_t = 1, 2, \ldots, N_t$, t = 1, 2, …, m.

The following properties of $A_{i_1}^1 \cup \ldots \cup A_{i_m}^m$ are well known

(i) $A_{i_1}^1 \cup \ldots \cup A_{i_m}^m$ 's are m-idempotent i.e.

$$(A_{i_1}^1 \cup \ldots \cup A_{i_m}^m)^2 = A_{i_1}^1 \cup \ldots \cup A_{i_m}^m$$

i.e. each $\left(A_{i_t}^t\right)^2 = A_{i_t}^t$, t = 1, 2, …, m.

(ii) They are n-orthogonal i.e.

$$A_{i_1}^1 . A_{j_t}^t = 0,\ i_t \neq j_t\ ;\ t = 1, 2, \ldots, m.$$

(iii) They give a spectral n-decomposition

$$P_1 \cup \ldots \cup P_n = \sum_{i_1=1}^{N_1} \lambda_{i_1}^1 A_{i_1}^1 \cup \ldots \cup \sum_{i_m=1}^{N_m} \lambda_{i_m}^m A_{i_m}^m .$$

It follows from (i) to (iii), that

$$P^K = P_1^{K_1} \cup \ldots \cup P_m^{K_m} =$$

$$\left(\sum_{i_1=1}^{N_1}\lambda_{i_1}^1 \; A_{i_1}^1\right)^{K_1} \cup \ldots \cup \left(\sum_{i_m=1}^{N_m}\lambda_{i_m}^m \; A_{i_m}^m\right)^{K_m}$$

$$= \sum_{i_1=1}^{N_1}\lambda_{i_1}^{K_1} A_{i_1}^1 \cup \ldots \cup \sum_{i_m=1}^{N_m}\lambda_{i_m}^{K_m} A_{i_m}^m$$

$$= \sum_{i_1=1}^{N_1}\lambda_{i_1}^{K_1} U_{i_1}^1 \; V_{i_1}^{1'} \cup \ldots \cup \sum_{i_m=1}^{N_m}\lambda_{i_m}^{K_m} U_{i_m}^m \; V_{i_m}^{m'} \; .$$

Also we know that

$$P^K = UD^K \; U^{-1} = U^1 D_1^{K_1} (U^1)^{-1} \cup \ldots \cup U^m D_m^{K_m} (U^m)^{-1}$$

where $U = \{U_1^1, \ldots, U_{N_1}^1\} \cup \ldots \cup \{U_1^m, \ldots, U_{N_m}^m\}$ and

$$D = D_1 \cup D_2 \cup \ldots \cup D_m$$

$$= \begin{bmatrix} \lambda_1^1 & 0 & \ldots & 0 \\ 0 & \lambda_2^1 & \ldots & 0 \\ \vdots & \vdots & & \vdots \\ 0 & 0 & \cdots & \lambda_{N_1}^1 \end{bmatrix} \cup \ldots \cup \begin{bmatrix} \lambda_1^m & 0 & \ldots & 0 \\ 0 & \lambda_2^m & \ldots & 0 \\ \vdots & \vdots & & \vdots \\ 0 & 0 & \cdots & \lambda_{N_m}^m \end{bmatrix}.$$

Since the n-latent n-vectors are determined uniquely only upto a multiplicative constant, we have chosen them such that

$$U_{i_1}^{1'} V_{i_1} \cup \ldots \cup U_{i_m}^{m'} V_{i_m} = (1 \cup \ldots \cup 1).$$

One can work for any m-power of P to know $\lambda_{i_t}^t$ 's and $A_{i_t}^t$ 's; t = 1, 2, …, m. Now even if we say $P^K = P_1^{K_1} \cup \ldots \cup P_m^{K_m}$ we work for K = ($K_1$, …, $K_m$) and when the working with any $P_t$ is over that $t^{th}$ component remains as it is and calculations are performed for the rest of the components of P. With the advent of the appropriate programming using computers simultaneous working is easy; also one needs to know in the present technologically advanced age one cannot think of computing one by one and also things do not occur like that in many situations. So under these circumstances only the adaptation of n-matrices plays a vital role by saving both time and economy. Also stage by stage comparison of the simultaneous occurrence of n-events is possible.

Matrix theory has been very successful in describing the interrelations between prices, outputs and demands in an economic model. Here we just discuss some simple models based on the ideals of the Nobel-laureate Wassily Leontief. Two types of models discussed are the closed or input-output model and the open or production model each of which assumes some economic parameter which describe the inter relations between the industries in the economy under considerations. Using matrix theory we evaluate certain parameters.

The basic equations of the input-output model are the following:

$$\begin{bmatrix} a_{11} & a_{12} & \cdots & a_{1n} \\ a_{21} & a_{22} & \cdots & a_{2n} \\ \vdots & \vdots & & \vdots \\ a_{n1} & a_{n2} & \cdots & a_{nn} \end{bmatrix} \begin{bmatrix} p_1 \\ p_2 \\ \vdots \\ p_n \end{bmatrix} = \begin{bmatrix} p_1 \\ p_2 \\ \vdots \\ p_n \end{bmatrix}$$

each column sum of the coefficient matrix is one

i. $p_i \geq 0$, i = 1, 2, …, n.
ii. $a_{ij} \geq 0$, i , j = 1, 2, …, n.
iii. $a_{ij} + a_{2j} + \ldots + a_{nj} = 1$

for j = 1, 2 , …, n.

$$p = \begin{bmatrix} p_1 \\ p_2 \\ \vdots \\ p_n \end{bmatrix}$$

are the price vector. $A = (a_{ij})$ is called the input-output matrix

$$Ap = p \text{ that is, } (I - A)\, p = 0.$$

Thus A is an exchange matrix, then Ap = p always has a nontrivial solution p whose entries are nonnegative. Let A be an exchange matrix such that for some positive integer m, all of the entries of $A^m$ are positive. Then there is exactly only one linearly independent solution of (I – A) p = 0 and it may be chosen such that all of its entries are positive in Leontief open production model.

In contrast with the closed model in which the outputs of k industries are distributed only among themselves, the open model attempts to satisfy an outside demand for the outputs. Portions of these outputs may still be distributed among the industries themselves to keep them operating, but there is to be some excess some net production with which to satisfy the outside demand. In some closed model, the outputs of the industries were fixed and our objective was to determine the prices for these outputs so that the equilibrium condition that expenditures equal incomes was satisfied.

$x_i$ = monetary value of the total output of the $i^{th}$ industry.

$d_i$ = monetary value of the output of the $i^{th}$ industry needed to satisfy the outside demand.

$\sigma_{ij}$ = monetary value of the output of the $i^{th}$ industry needed by the $j^{th}$ industry to produce one unit of monetary value of its own output.

With these qualities we define the production vector.

$$x = \begin{bmatrix} x_1 \\ x_2 \\ \vdots \\ x_k \end{bmatrix}$$

the demand vector

$$d = \begin{bmatrix} d_1 \\ d_2 \\ \vdots \\ d_k \end{bmatrix}$$

and the consumption matrix,

$$C = \begin{bmatrix} \sigma_{11} & \sigma_{12} & \cdots & \sigma_{1k} \\ \sigma_{21} & \sigma_{22} & \cdots & \sigma_{2k} \\ \vdots & \vdots & & \vdots \\ \sigma_{k1} & \sigma_{k2} & \cdots & \sigma_{kk} \end{bmatrix}.$$

By their nature we have

$$x \geq 0,\ d \geq 0 \text{ and } C \geq 0.$$

From the definition of $\sigma_{ij}$ and $x_j$ it can be seen that the quantity

$$\sigma_{i1}\, x_1 + \sigma_{i2}\, x_2 + \ldots + \sigma_{ik}\, x_k$$

is the value of the output of the $i^{th}$ industry needed by all k industries to produce a total output specified by the production vector x.

Since this quantity is simply the $i^{th}$ entry of the column vector Cx, we can further say that the $i^{th}$ entry of the column vector x – Cx is the value of the excess output of the $i^{th}$ industry available to satisfy the outside demand. The value of the outside demand for the output of the $i^{th}$ industry is the $i^{th}$ entry of the demand vector d; consequently; we are led to the following equation:

$$x - Cx = d \text{ or}$$
$$(I - C)\, x = d$$

for the demand to be exactly met without any surpluses or shortages. Thus, given C and d, our objective is to find a production vector $x \geq 0$ which satisfies the equation $(I - C)x = d$.

A consumption matrix C is said to be productive if $(1 - C)^{-1}$ exists and $(1 - C)^{-1} \geq 0$.

A consumption matrix C is productive if and only if there is some production vector $x \geq 0$ such that $x > Cx$.

A consumption matrix is productive if each of its row sums is less than one. A consumption matrix is productive if each of its column sums is less than one.

Now we will formulate the Smarandache analogue for this, at the outset we will justify why we need an analogue for those two models.

Clearly, in the Leontief closed Input – Output model,
$p_i$ = price charged by the $i^{th}$ industry for its total output in reality need not be always a positive quantity for due to competition to capture the market the price may be fixed at a loss or the demand for that product might have fallen down so badly so that the industry may try to charge very less than its real value just to market it.

Similarly $a_{ij} \geq 0$ may not be always be true. Thus in the Smarandache Leontief closed (Input – Output) model (S-Leontief closed (Input-Output) model) we do not demand $p_i \geq 0$, $p_i$ can be negative; also in the matrix $A = (a_{ij})$,

$$a_{1j} + a_{2j} + \ldots + a_{kj} \neq 1$$

so that we permit $a_{ij}$'s to be both positive and negative, the only adjustment will be we may not have $(I - A)\, p = 0$, to have only one linearly independent solution, we may have more than one and we will have to choose only the best solution.

As in this complicated real world problems we may not have in practicality such nice situation. So we work only for the best solution.

On similar lines we formulate the Smarandache Leontief open model (S-Leontief open model) by permitting that $x \geq 0$ , $d \geq 0$ and $C \geq 0$ will be allowed to take $x \leq 0$ or $d \leq 0$ and or $C \leq 0$ . For in the opinion of the author we may not in reality have the monetary total output to be always a positive quality for all industries and similar arguments for di's and $C_{ij}$'s.

When we permit negative values the corresponding production vector will be redefined as Smarandache production vector (S-production vector) the demand vector as Smarandache demand vector (S-demand vector) and the consumption matrix as the Smarandache consumption matrix (S-consumption matrix). So when we work out under these assumptions we may have different sets of conditions

We say productive if $(1 - C)^{-1} \geq 0$, and non-productive or not up to satisfaction if $(1 - C)^{-1} < 0$.

The reader is expected to construct real models by taking data's from several industries. Thus one can develop several other properties in case of different models.

Matrix theory has been very successful in describing the interrelations between prices outputs and demands.

Now when we use n-matrices in the input – output model we can under the same set up study the price vectors of all the goods manufactured by that industry simultaneously. For in the present modernized world no industry thrives only in the production one goods. For instance take the Godrej industries it manufacturers several goods from simple locks to bureau. So if they want to study input output model to each and every goods it has to work several times with the exchange matrix; but with the introduction of n-mixed matrices we can use the n-matrix as the input output n-model to study interrelations between the prices outputs and demands of each and every goods manufactured by that industry. Suppose the industry manufactures n-goods, $n \geq 2$.

Thus $A = A_1 \cup \ldots \cup A_n$ is an exchange n-matrix where each $A_i$ is a $n_i \times n_i$ matrix $i = 1, 2, \ldots, n$. The basic n-equations of the input – output model is the following

$$\begin{bmatrix} a_{11}^1 & a_{12}^1 & \cdots & a_{1n_1}^1 \\ a_{21}^1 & a_{22}^2 & \cdots & a_{2n_1}^1 \\ \vdots & \vdots & & \vdots \\ a_{n_1 1}^1 & a_{n_1 2}^1 & \cdots & a_{n_1 n_1}^1 \end{bmatrix} \begin{bmatrix} p_1^1 \\ \vdots \\ p_{n_1}^1 \end{bmatrix} \cup$$

$$\begin{bmatrix} a_{11}^2 & a_{12}^2 & \cdots & a_{1n_2}^2 \\ a_{21}^2 & a_{22}^2 & \cdots & a_{2n_2}^2 \\ \vdots & \vdots & & \vdots \\ a_{n_2 1}^2 & a_{n_2 2}^2 & \cdots & a_{n_2 n_2}^2 \end{bmatrix} \begin{bmatrix} p_1^2 \\ \vdots \\ p_{n_2}^2 \end{bmatrix} \cup \ldots \cup$$

$$\begin{bmatrix} a^n_{11} & a^n_{12} & \dots & a^n_{1n_n} \\ a^n_{21} & a^n_{22} & \dots & a^n_{2n_n} \\ \vdots & \vdots & & \vdots \\ a^n_{n_n1} & a^n_{n_n2} & \dots & a^n_{n_nn_n} \end{bmatrix} \begin{bmatrix} p^n_1 \\ \vdots \\ p^n_{n_1} \end{bmatrix} =$$

$$\begin{bmatrix} p^1_1 \\ \vdots \\ p^1_{n_1} \end{bmatrix} \cup \begin{bmatrix} p^2_1 \\ \vdots \\ p^1_{n_2} \end{bmatrix} \cup \dots \cup \begin{bmatrix} p^n_1 \\ \vdots \\ p^n_{n_n} \end{bmatrix}$$

each n column sum of the coefficient n-matrix is $(1 \cup \dots \cup 1)$

(i) $p^t_i \geq 0$; t = 1, 2, …, n.

(ii) $a^t_{i_t j_t} \geq 0$; $i_t$, $j_t$ = 1, 2, …, $n_t$ and t = 1, 2, …, n.

(iii) $a^t_{ij_1} + a^t_{2j_t} + \dots + a^t_{n_t j_t} = 1$ for $j_t$ = 1, 2, …, $n_t$ and t = 1, 2, …, n.

$$p = p_1 \cup \dots \cup p_n = \begin{bmatrix} p^1_1 \\ p^1_2 \\ \vdots \\ p^1_{n_1} \end{bmatrix} \cup \begin{bmatrix} p^2_1 \\ p^2_2 \\ \vdots \\ p^2_{n_2} \end{bmatrix} \cup \dots \cup \begin{bmatrix} p^n_1 \\ p^n_2 \\ \vdots \\ p^n_{n_n} \end{bmatrix}$$

are the price n-vector of the n-goods.

$$A = A_1 \cup \dots \cup A_n = (a^1_{i_1 j_1}) \cup \dots \cup (a^n_{i_n j_n})$$

is called the input-output n-matrix.

$$Ap = p \text{ that is } (I - A)\, p = 0 \cup \dots \cup 0$$

i.e. $(I_1 - A_1)\, p_1 \cup \dots \cup (I_n - A_n)\, p_n = 0 \cup \dots \cup 0$.

Thus A is an exchange n-matrix then Ap = p always has a nontrivial n-solution $p = p_1 \cup \dots \cup p_n$, whose entries are nonnegative. Let A be the exchange n-matrix such that for some n-positive integers $(m_1, \dots, m_n)$ all the entries of $A^m = A_1^{m_1} \cup \dots \cup A_n^{m_n}$ are positive. Then there is exactly only one linearly n-independent solution of $(I - A)p = 0 \cup \dots \cup 0$ and that it may be chosen such that all of its entries are positive in Leontief open production n-model.

Thus the model provides at a time i.e. simultaneously the price n-vector i.e. the price vector of each of the n-goods. When n = 1 we see the structure corresponds to the Leontief open model. When n = 2 we get the Leontief economic bi models. This n-model is useful when the industry manufactures more than one goods and it not only saves time and economy but it renders itself stage by stage comparison of the price n-vector which is given by $p = p_1 \cup \ldots \cup p_n$.

Now we proceed onto describe the S-Leontief open n-model using n-matrices.

In reality we may not always have the exchange n-matrix $A = A_1 \cup \ldots \cup A_n = (a^1_{i_1 j_1}) \cup \ldots \cup (a^n_{i_n j_n}),\ a^t_{i_t j_t} \geq 0.$ For it can also be both positive or negative. Thus in S-Leontief closed (input - output) n-model we do not demand $p^t_{i_t} \geq 0$, $p^t_{i_t}$ can be negative also in the n-matrix $A = A_1 \cup \ldots \cup A_n = (a^1_{i_1 j_1}) \cup \ldots \cup (a^n_{i_n j_n})$ where $a^t_{ij_t} + \ldots + a^t_{K_t j_t} \neq 1$ for every t = 1, 2, …, n. i.e. we permit $a^t_{i_t j_t}$ to be both positive and negative, the only adjustment will be, we may not have $(I - A)p = 0 \cup \ldots \cup 0$ to have only one n-linearly independent solution, we may have more than one and we will have to choose only the best solution which will be helpful to the economy of the nation. The best by no means should favour in the interrelation high prices but a medium price with most satisfactory outputs and best catering to the demands as it is an economic n-model.

So n-matrices will be highly helpful and out of one set of solution which will have n-components associated with the exchange n-matrix $A = A_1 \cup \ldots \cup A_n$, we have to pick up from the nontrivial solution $p_1 = p_1 \cup \ldots \cup p_n$ the best suited $p_i$'s and once again find a $p' = p'_1 \cup \ldots \cup p'_n$ with the estimated $p_i$'s from the earlier p remain as zero and choose the best $p'_j$ for the solution $p'$ and so on. The final $p = p_1 \cup \ldots \cup p_n$ will be filled with the best $p_i$'s and $p_j$'s and so on.

Thus the solution would be the best suited solution of the economic model.

The difference between Leontief closed or input output n-model and the S-Leontief closed or input output economic n-model is that in the Leontief model there is only one independent solution where as in the S-Leontief closed input output economic n model we can choose the best solution from the set of solutions so that best solution also may vary from person to person for what is best for one may not be best for the other that too when it describes the interrelations between prices, outputs and demands in an economic n-model.

Now we briefly describe the Leontief open production n-model. In contrast with Leontief closed n-model here the n-set of or n-tuple of industries say ($K_1, \ldots, K_n$) where output of $K_i$ industries are distributed only among themselves the open n model attempts to satisfy an outside demand for the n-outputs, true for i = 1, 2, …, n. Portions of these n-outputs may still be distributed among the ($K_1, \ldots, K_n$) set of industries themselves to keep them operating, but there is to be some excess some net production with which to satisfy the outside demand.

In the closed n-model the n-outputs of the industries were fixed and the objective was to determine the n-prices for these n-outputs so that the equilibrium condition that expenditures equal income was satisfied.

$x^t_{i_t}$ = monetary value of the $i_t^{th}$ industry from the $t^{th}$ unit i.e. we have

$K_1$ = industries in the first unit denoted by $c_1$
$K_2$ = industries in the second unit denoted by $c_2$
$\vdots$
$K_t$ = industries in the $t^{th}$ unit denoted by $c_t$

and so on

$K_n$ – industries in the $n^{th}$ unit denoted by $c_n$.

$d^t_{i_t}$ – monetary value of the output of the $i_t^{th}$ industry need to satisfy the outside demand.

$\sigma^t_{i_t j_t}$ – monetary value of the output of the $i_t^{th}$ industry needed by the $j_t^{th}$ industry to produce one unit of monetary value of its own profit.

This is true for every t; t = 1, 2, …, n. With these qualities we define the n-production vector which is a n-vector.

$$x = x_1 \cup \ldots \cup x_n = \begin{bmatrix} x_1^1 \\ \vdots \\ x_{K_1}^1 \end{bmatrix} \cup \begin{bmatrix} x_1^2 \\ \vdots \\ x_{K_2}^2 \end{bmatrix} \cup \ldots \cup \begin{bmatrix} x_1^n \\ \vdots \\ x_{K_n}^n \end{bmatrix},$$

the n-demand vector which is a n-vector,

$$d = d_1 \cup \ldots \cup d_n \begin{bmatrix} d_1^1 \\ \vdots \\ d_{K_1}^1 \end{bmatrix} \cup \ldots \cup \begin{bmatrix} d_1^n \\ \vdots \\ d_{K_n}^n \end{bmatrix}$$

and the n-consumption matrix which is a n matrix

$$c = c_1 \cup \ldots \cup c_n$$

$$= \begin{bmatrix} \sigma_{11}^1 & \sigma_{12}^1 & \ldots & \sigma_{1K_1}^1 \\ \sigma_{21}^1 & \sigma_{22}^1 & \ldots & \sigma_{2K_1}^1 \\ \vdots & \vdots & & \vdots \\ \sigma_{K_1}^1 & \sigma_{K_1 2}^1 & \ldots & \sigma_{K_1 K_1}^1 \end{bmatrix} \cup \ldots \cup \begin{bmatrix} \sigma_{11}^n & \sigma_{12}^n & \ldots & \sigma_{1K_n}^n \\ \sigma_{21}^n & \sigma_{22}^n & \ldots & \sigma_{2Kn}^n \\ \vdots & \vdots & & \vdots \\ \sigma_{K_n 1}^n & \sigma_{K_n 2}^n & \ldots & \sigma_{K_n K_n}^n \end{bmatrix}.$$

We have $x \geq 0 \cup \ldots \cup 0$ i.e. $x = x_1 \cup \ldots \cup x_n \geq 0 \cup 0 \cup \ldots \cup 0$, $d \geq 0 \cup \ldots \cup 0$ i.e. $d = d_1 \cup \ldots \cup d_n \geq 0 \cup \ldots \cup 0$ and $c \geq 0 \cup \ldots \cup 0$ i.e. $c = c_1 \cup c_2 \cup \ldots \cup c_n \geq 0 \cup \ldots \cup 0$.

From the definition of $\sigma_{i_t j_t}^t$ and $x_{j_t}^t$ it can be seen that the quantity $\sigma_{i_t 1}^t x_1^t + \sigma_{i_t 2}^t x_2^t + \ldots + \sigma_{i_t K_t}^t$ is the value of the $i_t^{th}$ industry of the $t^{th}$ unit needed for all $K_t$ industries to produce a total output specified by the production component vector $x^t = x_1^t \cup \ldots x_{K_t}^t$ of the n-vector. $x = x_1 \cup x_2 \cup \ldots \cup x_n$. This is true for each t; t = 1, 2, …, n. Since the quantity is simply the $i_t^{th}$ entry of the $t^{th}$ unit column n vector $c_t x^t$ we can further say that the $i_t^{th}$ entry of the column vector $x^t - c_t x^t$ is the value of the excess output of the $i_t^{th}$ industry available to satisfy the outside demand for t = 1, 2, …, n. Thus the excess n-output of the ($i_1$, …, $i_n$) industry is given by the n-column vector $x - cx = x^1 - c_1 x^1 \cup \ldots \cup x^n - c_n x^n$. The value of the outside demand for the n output $(i_1, \ldots, i_n)$ is the $(c_1, \ldots, c_n)^{th}$ entry of the demand vector $d = d_1 \cup \ldots \cup d_n$. Consequently, we are led to the following equation.

$$x - cx = x^1 - c_1 x^1 \cup \ldots \cup x^n - c_n x^n = d_1 \cup \ldots \cup d_n$$
$$(I - c)\, x = d$$
$$(I_1 - c_1)\, x^1 \cup \ldots \cup (I_n - c_n) x^n = d_1 \cup \ldots \cup d_n,$$

for the demand to be exactly met without any surplus or shortages. Thus given c and d our objective is to find a production n-vector $x = x_1 \cup \ldots \cup x_n \geq 0 \cup \ldots \cup 0$ which satisfies the n-equation $(I - c)\, x = d$ $(I_1 - c_1)\, x_1 \cup \ldots \cup (I_n - c_n)\, x_n = d_1 \cup \ldots \cup d_n$.

The consumption n-matrix $c = c_1 \cup \ldots \cup c_n$ is said to be n-productive if $(1 - c)^{-1} = (1 - c_1)^{-1} \cup \ldots \cup (1 - c_n)^{-1}$ exists and $(1 - c)^{-1} \geq 0 \cup \ldots \cup 0$. A consumption n-matrix $c = c_1 \cup \ldots \cup c_n$ is productive if and only if there is some production n-vector $x = x_1 \cup \ldots \cup x_n \geq 0 \cup \ldots \cup 0$ such that $x > cx$; $x_1 \cup \ldots \cup x_n > c_1 x_1 \cup \ldots \cup c_n x_n$.

A consumption n-matrix is productive if each of the n-row sums is less than one. A consumption n-matrix is n-productive if each of its column sum is less than one.

Now we will formulate the Smarandache analogue for this, at the outset we will justify why we need an analogue for the open or production n-model.

In the Leontief open n-model we may assume also $x \leq 0$, or $d \leq 0$ and or $c \leq 0$. For in the opinion of the author we may not in reality have the monetary total output to be always a positive quantity for all industries and similar arguments for $d_i^t$'s and $c_{ij}^t$'s.

When we permit negative values the corresponding production n-vector will be redefined as S-production n-vector the demand n-vector as S-demand n-vector and the consumption n-matrix as S-consumption n-matrix. Under these assumptions we may have different sets of conditions.

We say n-productive if $(1 – c)^{-1} \geq 0$ and non n-productive or not upto satisfaction if $(1 – c)^{-1} < 0$.

Now we have given some application of these n-matrices to industrial problems.

Finally it has become pertinent here to mention that in the consumption n matrices a particular industry or many industries can be used in several or more than one consumption matrix. So in this situation only the open Leontief n-model will serve it purpose. Also we can study the performance such industries which is in several groups i.e. in several $c_i$'s. One can also simultaneously study the group in which an industry has the best performance also the group in which it has the worst performance. In such situation only this model is handy.

**Chapter Four**

# SUGGESTED PROBLEMS

In this chapter we suggest some problems for the readers. Solving these problems will be a great help to understand the notions given in this book.

1. Find all p-subspaces of the n-vector space $V = V_1 \cup V_2 \cup V_3 \cup V_4$ where n = 4 and $p \leq 4$ over Q.

$$V_1 = \left\{ \begin{pmatrix} a & b & e \\ c & d & f \end{pmatrix} \middle| a, b, c, d, e, f \in Q \right\},$$

   $V_2 = (Q \times Q \times Q \times Q)$ over Q,
   $V_3$ = {Q[x] contains only polynomials of degree less than or equal to 6 with coefficients from Q} and

$$V_4 = \left\{ \begin{pmatrix} a & b & c & d \\ e & f & g & h \end{pmatrix} \middle| a, b, \ldots, g, h \in Q \right\}.$$

   What is the 4-dimension of V? Find a 4 basis of V.

2. Let $V = V_1 \cup V_2 \cup V_3$ and $W = W_1 \cup W_2 \cup W_3 \cup W_4$ be 3 vector space and 4 vector space over the field Q of 3 dimension (3, 2, 4) and 4 dimension (5, 3, 4, 2) respectively.

Find a 3 linear transformation from V to W. Also find a shrinking 3 linear transformation from V into W.

3. Let $V = V_1 \cup V_2 \cup V_3$ and $W = W_1 \cup W_2 \cup W_3$ be 3-vector spaces of dimensions (4, 2, 3) and (3, 5, 4) respectively defined over Q. Find the 3 linear transformation from V to W. What is the 3 dimension of the 3-vector space of all 3 linear transformation from V into W?

4. Let $V = V_1 \cup V_2 \cup V_3 \cup V_4$ be a 4-vector space defined over Q of dimension (3, 4, 2, 1). Give a 4 linear operator T on V.
Verify: 4 rank T + 4 nullity T = n dim V = (3, 4, 2, 1).

5. Define $T: V \rightarrow W$ be a 4 linear operator where $V = V_1 \cup V_2 \cup V_3 \cup V_4$ and $W = W_1 \cup W_2 \cup W_3 \cup W_4$ with 4-dimension (3, 2, 4, 5) and (4, 3, 5, 2) respectively, such that 4 kerT is a 4-dimensional subspace of V. Verify 4 rank T + 4 nullily T = 4 dim V = (3, 2, 4, 5).

6. Explicitly describe the n-vector space of n-linear transformations $L^n$ (V,W) of $V = V_1 \cup V_2 \cup V_3$ into $W = W_1 \cup W_2 \cup W_3 \cup W_4$ over Q of 3-dimension (3, 2, 4) and 4-dimension (4, 3, 2, 5) respectively.

7. What is n-dimension of $L^n$ (V,W) given in the problem 6?

8. For $T = T_1 \cup T_2 \cup T_3$ defined for V and W given in problem 6; $T_1 : V_1 \rightarrow W_3$, $T_1$ (x y z) = (x + y, y + z) for all x, y, z $\in V_1$, $T_2 : V_2 \rightarrow W_2$ defined by $T_2 (x_1, y_1) = (x_1 + y_1, 2y_1, y_1)$ for all $x_1, y_1, \in V_2$ and $T_3 : V_3 \rightarrow W_4$ defined by $T_3$ (a, b, c, d) = (a + b, b + c, c + d, d + a, a + b + d) for all a, b, c, d $\in V_3$. Prove 3 rank T + 3 nullity T = dim V = (3, 2, 4).

9. Let $V = V_1 \cup V_2 \cup V_3 \cup V_4$ be a 4 vector space over Q, where $V_1 = Q \times Q \times Q$, $V_2 = Q \times Q \times Q \times Q$, $V_3 = Q \times Q$, $V_4 = Q \times Q \times Q \times Q \times Q$, $T^j = T_1^j \cup T_2^j \cup T_3^j \cup T_4^j : V \rightarrow V$

$T_i^j : V_i \to V_i$; i = 1, 2, 3, 4.

Define two distinct 4 transformations $T^1$ and $T^2$ and find $T^1 \circ T^2$ and $T^2 \circ T^1$.

10. Give an example of a special linear 6-transformation $T = T_1 \cup T_2 \cup \ldots \cup T_6$ of V into W where V and W are 6 vector space of same 6 dimension.

11. Let T : V → W where 3-dim V = (3, 7, 8) and 3-dim W = (8, 3, 7). Give an example of T and find $T^{-1}$. Define T only as a 3 linear transformation for which $T^{-1}$ cannot be found.

12. Derive for a n-vector space the Gram-Schmidt n-orthogonalization process.

13. Prove every finite n-dimensional inner product n-space has an n-orthonormal basis.

14. Give an example of a 4-orthogonal matrix.

15. Give an example of a 5-anitorthogonal matrix.

16. Give an example of a 7-semi orthogonal matrix.

17. Give an example of a 5-semi antiorthogonal matrix.

18. Is $A = \begin{bmatrix} 3 & 1 & 0 & 2 \\ 1 & 1 & 6 & 1 \\ 0 & 2 & 0 & 1 \\ 1 & 0 & 5 & 0 \end{bmatrix} \cup \begin{bmatrix} 1 & 1 & 1 & 1 & 0 \\ 2 & 0 & 0 & 2 & 1 \\ 0 & 1 & 2 & 3 & 4 \\ 5 & 6 & 7 & 0 & 12 \end{bmatrix} \cup \begin{bmatrix} 3 & 1 & 8 \\ 0 & 1 & 1 \\ 1 & 0 & 1 \\ 0 & 1 & 4 \\ -1 & 0 & 1 \\ 1 & 1 & 0 \end{bmatrix}$ a 3-semi orthogonal 3 matrix?

19. Find the 4-eigen values, 4-eigen vectors of $A = A_1 \cup A_2 \cup A_3 \cup A_4 =$

$$= \begin{bmatrix} 3 & 0 & 1 & 0 \\ 0 & 4 & 1 & 3 \\ 0 & 5 & 0 & 1 \\ 1 & 2 & 1 & 0 \end{bmatrix} \cup \begin{bmatrix} 3 & 0 & 1 \\ 0 & 1 & 4 \\ 0 & 0 & 5 \end{bmatrix} \cup \begin{bmatrix} 5 & 1 \\ 0 & 1 \end{bmatrix} \cup \begin{bmatrix} 0 & 1 & 2 & 3 & 6 \\ 0 & 2 & 1 & 0 & 2 \\ 0 & 0 & 2 & 1 & 1 \\ 0 & 0 & 1 & 0 & 0 \\ 0 & 0 & 0 & 0 & 5 \end{bmatrix}.$$

Find the 4-minimal polynomial and the 4-characteristic polynomial associated with A. Is A a diagonalizable transformation? Justify your claim.

20. Give an example of 5-linear transformation on $V = V_1 \cup V_2 \cup ... \cup V_5$ which is not a 5-linear operator on V.

21. Let $V = V_1 \cup V_2 \cup V_3$ be a 3-vector space over the field Q of finite (5, 3, 2) dimension over Q. Give a special 3 linear operator on V. Give a 3 linear transformation on V which is not a special linear operator on V.

22. Define a 3-innerproduct on V given in the above problem and construct a normal 3 linear operator T on V such that T*T = TT*.

23. Let $V = V_1 \cup V_2 \cup V_3 \cup V_4$ be a 4-vector space of (3, 5, 2, 4) dimension over Q. Find a 4-linear operator T on V so that the 4-minimal polynomial of T is the same as 4-characteristic polynomial of T. Give a 4-linear operator U on V so that the 4-minimal polynomial is different from the 4-characteristic polynomial.

24. Let $V = V_1 \cup V_2 \cup V_3 \cup V_4 \cup V_5$ be a 5-vector space over Q of (2, 3, 4, 5, 6) dimension over Q. Construct a linear operator T on V so that T is 5-diagonalizable.

25. Let $V = V_1 \cup V_2 \cup V_3 \cup V_4$ be a 4-vector space over Q. Define a suitable T and find the n-monic generator of the 4-ideals of the polynomials over Q which 4-annihilate T.

Prove or disprove every 4-linear operator T on V need not 4-annihulate T.

26. State and prove the Cayley Hamilton theorem for n-linear operator on a n-vector space V.

27. Let $V = V_1 \cup V_2 \cup V_3 \cup V_4 \cup V_5$ be a 5-vector space over Q of (2, 4, 6, 3, 5) dimension over Q. Give a 5-basis of V so that Cayley Hamilton Theorem is true. Is Cayley Hamilton Theorem true for every set of 5-basis of V? Justify your claim.

28. Given $V = V_1 \cup V_2 \cup V_3 \cup V_4$ is a 4-vector space over Q of dimension (3, 7, 4, 2). Construct a T, a 4 linear operator on V so that V has a 4-subspace 4-invariat under T. Does V have any 4-linear operator T and a non-trivial 4-subspace W so that W is 4-invariant under T? Justify your answer.

29. Let $V = V_1 \cup V_2 \cup V_3 \cup V_4 \cup V_5$ be a 5-vector space of (2, 4, 5, 3, 7) dimension over Q. Construct a 5-linear operator V on T so that the 5-minimal polynomial associated with T is linearly factorizable. Find a T on V so that the 5-minimal polynomial does not factor linearly over Q.

30. Let $V = V_1 \cup V_2 \cup V_3$ be a 3-vector space of (2, 4, 3) dimension over Q. Find $L^3$ (V, V) the set of all 3-linear transformations on V. Suppose $L_S^3(V, V)$ is the set of all special 3-linear transformations on V.

    a. Prove $L_S^3$ (V, V) $\subseteq L^3$(V, V).
    b. What is the 3-dimension of $L^3$ (V, V)?
    c. What is the 3-dimension of $L_S^3$ (V, V)?
    d. Find a set of 3-orthogonal 3 basis for $L_S^3$ (V, V).
    e. Find a set of 3-orthonormal 3-basis for $L^3$ (V, V)
    f. Find a $T : V \to V$, T only a 3-linear transformation which has a nontrivial 3-null space.
    g. Find the 3-rank T of that is given in (6)

h. Can any $T \in L_S^3$ (V, V) have nontrivial 3-null space? Justify your answer.
i. Define a 3-unitary operator on V.
j. Define a 3-normal operator on V which not 3-unitary.

31. Let V and W be two 6-inner product spaces of same dimension ($W \neq V$) defined over the same field F. Define a T linear operator from V into W which preserves inner products by taking (3, 4, 6, 2, 1, 5) to be the dimension of V and (6, 5, 4, 2, 3, 1) is the dimension of W.
Does every $T \in L^6$ (V, W) preserve inner product? Justify your claim.

32. Given $V = V_1 \cup V_2 \cup V_3$ is a (4, 5, 3) dimensional 3-vector space over Q. Give an example of a 3-linear operator T on V which is 3-diagonalizable. Does their exist a 3-linear operator $T'$ on V such that $T'$ is not 3 diagonalizable? Justify your answer.

33. Let $V = V_1 \cup V_2 \cup V_3 \cup V_4 \cup V_5$ be a (3, 4, 5, 2, 6) dimension 5-vector space over Q. Define a 5 linear operator T on V and decompose it into the 5-nilpotent operator and 5-diagonal operator.

a. Does there exist a 5-linear operator T on V such that the 5-diagonal part is zero, i.e., the operator T is nilpotent?
b. Does there exist a 5-linear operator P on V such that it is completely 5-diagonal and the 5-nilpotent part of it is zero.
c. Give examples of the above mentioned 5-operator in (1) and (2)
d. What is the form of the 5-minimal polynomial in case of (1) and (2)?

34. Define for a n-vector space V over a field F the notion of n-independent n-subspaces of V. Give an example when n = 4.

35. Let $V = V_1 \cup V_2 \cup \ldots \cup V_6$ be a 6-vector space over Q. Define a 6-linear operator E on V such that $E^2 = E$.

36. Let $V = V_1 \cup V_2 \cup V_3 \cup V_4$ be a 4-vector space over Q of (3, 4, 5, 2) dimension. Suppose $V = \left(W_1^1 \oplus W_2^1\right)$ $\left(W_1^2 \oplus W_2^2 \oplus W_3^2\right) \cup \left(W_1^3 \oplus W_2^3 \oplus W_3^3\right) \cup \left(W_1^4 \oplus W_2^4\right)$,
Define 4-linear operators, $E_i^1 \cup E_j^2 \cup E_k^3 \cup E_m^4$; i = 1, 2; j = 1, 2, 3; k = 1, 2, 3 and m = 1, 2 such that each $E_p^i$ is a projection, i = 1, 2, 3, 4 and

$$E_p^i E_j^i \begin{cases} = 0 \text{ if } p \neq j \\ = E_p^i \text{ if } p = j. \end{cases}.$$

37. Prove if T is any 4-linear operator on V then $TE_j^i = E_j^i T$ for i = 1,2, 3,4. j = 1, 2 or 1, 2, 3 or 1, 2, for the V given in the problem 36.

38. Given $V = V_1 \cup V_2 \cup V_3 \cup V_4 \cup V_5$ to be a 5-vector space over Q of (2, 3, 4, 5, 6) dimension. Define T a linear operator on V and find the 5 minimal polynomial for T. Is every 5-subspace of V related with the 5-minimal polynomials i.e. the 5-null space of the minimal polynomials invariant under T?
Obtain the 5-nilpotent and 5-diagonalizable operator N and D respectively so that T = N+D.
Verify ND = DN for the same N and D of T.

39. If T is a 7-linear operator on $V = V_1 \cup V_2 \cup \ldots \cup V_7$ of (3, 2, 5, 1, 6, 4, 7) dimension over Q. Is the generalized Cayley Hamilton Theorem true for T?

40. Prove for a 3-vector spaces $V = V_1 \cup V_2 \cup V_3$ of (3, 4, 2) dimension over Q and $W = W_1 \cup W_2 \cup W_3$ of dimension (4, 5, 3) over Q if T is any 3 linear transformation find the 3 matrix associated with T. Find the 3-adjoint of T.

41. For any n-linear transformation T of a n vector space V = $V_1 \cup V_2 \cup ... \cup V_n$ of dimension ($n_1$, $n_2$, ..., $n_n$) into a m-vector space W (m>n) of dimension ($m_1$, $m_2$, …, $m_m$) over Q. Prove there exists a n-matrix A = ($A_1 \cup A_2 \cup ... \cup A_n$) which is related to T. Prove $L^n$ (V, W) $\cong$ {set of all n-matrices $A_1 \cup A_2 \cup ... \cup A_n$ where each $A_i$ is a $n_i \times m_j$ matrix with entries from Q}.

42. If V = $V_1 \cup V_2 \cup ... \cup V_n$ is a n-vector space over the field F of ($n_1$, $n_2$,…, $n_n$) dimension. If T : V $\rightarrow$ V is such that $T_i$ : $V_i \rightarrow V_i$; i = 1, 2, …, n. Show $L_n^S$ (V, V) $\cong$ {All n-mixed square matrices A = ($A_1 \cup A_2 \cup ... \cup A_n$) where $A_i$ is a $n_i \times n_i$ matrix with entries from F}.

43. Define n-norm on V an inner product space and is it possible to prove the Cauchy Schwarz inequality?

44. Derive Gram-Schmidt orthogonalization process for a n-vector space V with an inner product for a n-set of n-independent vectors in V.

45. Let V be a n-inner product space over F. W a finite dimensional n-subspace of V. Suppose E is a n orthogonal projection of V on W, with E an n-idempotent n-linear transformation of V onto W. $W^{\perp}$ the n-null space of E.
Prove V = W $\oplus$ $W^{\perp}$.

# FURTHER READING

# INDEX

**I**

**L**

**M**

**N**

# ABOUT THE AUTHORS

**Dr.W.B.Vasantha Kandasamy** is an Associate Professor in the Department of Mathematics, Indian Institute of Technology Madras, Chennai. In the past decade she has guided 12 Ph.D. scholars in the different fields of non-associative algebras, algebraic coding theory, transportation theory, fuzzy groups, and applications of fuzzy theory of the problems faced in chemical industries and cement industries.

She has to her credit 646 research papers. She has guided over 68 M.Sc. and M.Tech. projects. She has worked in collaboration projects with the Indian Space Research Organization and with the Tamil Nadu State AIDS Control Society. This is her $37^{th}$ book.

On India's 60th Independence Day, Dr.Vasantha was conferred the Kalpana Chawla Award for Courage and Daring Enterprise by the State Government of Tamil Nadu in recognition of her sustained fight for social justice in the Indian Institute of Technology (IIT) Madras and for her contribution to mathematics. (The award, instituted in the memory of Indian-American astronaut Kalpana Chawla who died aboard Space Shuttle Columbia). The award carried a cash prize of five lakh rupees (the highest prize-money for any Indian award) and a gold medal.

She can be contacted at vasanthakandasamy@gmail.com
You can visit her on the web at: http://mat.iitm.ac.in/~wbv

---

**Dr. Florentin Smarandache** is a Professor of Mathematics and Chair of Math & Sciences Department at the University of New Mexico in USA. He published over 75 books and 150 articles and notes in mathematics, physics, philosophy, psychology, rebus, literature.

In mathematics his research is in number theory, non-Euclidean geometry, synthetic geometry, algebraic structures, statistics, neutrosophic logic and set (generalizations of fuzzy logic and set respectively), neutrosophic probability (generalization of classical and imprecise probability). Also, small contributions to nuclear and particle physics, information fusion, neutrosophy (a generalization of dialectics), law of sensations and stimuli, etc. He can be contacted at smarand@unm.edu